\documentclass[12pt]{amsart}
\usepackage{amsmath,amsfonts,amssymb,amscd}
\usepackage{latexsym}

\def\HH{\mbox{${\mathcal H}$\kern-5.2pt${\mathcal H}$}}

\newtheorem{theorem}{Theorem}[section]
\newtheorem{proposition}[theorem]{Proposition}

\newtheorem{lemma}[theorem]{Lemma}
\newtheorem{corollary}[theorem]{Corollary}

\newtheorem{theorem }{Theorem}[section]
\newtheorem{proposition }[theorem]{Proposition}
\newtheorem{definition }[theorem]{Definition}
\newtheorem{lemma }[theorem]{Lemma}
\newtheorem{corollary }[theorem]{Corollary}
\newtheorem{notation }[theorem]{Notation}
\newtheorem{remark }[theorem]{Remark}
\newtheorem{example }[theorem]{Example}

\newtheorem{ theorem}{Theorem}[section]
\newtheorem{ proposition}[theorem]{Proposition}
\newtheorem{ definition}[theorem]{Definition}
\newtheorem{ lemma}[theorem]{Lemma}
\newtheorem{ corollary}[theorem]{Corollary}
\newtheorem{ notation}[theorem]{Notation}
\newtheorem{ remark}[theorem]{Remark}
\newtheorem{ example}[theorem]{Example}

\def\for{\  \hbox{ for } \ }
\def\iif{ \ \hbox{ if } \ }

\def\where{\  \hbox{ where } \ }
\def\and{\  \hbox{ and } \ }

\def\equal{\stackrel{\,\mathbf{def}}{= \kern-3pt =}}

\def\la{\lambda}
\def\La{\Lambda}
\def\om{\omega}

\def\th{\theta}
\def\al{\alpha}
\def\be{\beta}
\def\ga{\gamma}
\def\ep{\epsilon}

\def\Up{\Upsilon}
\def\de{\delta}
\def\De{\Delta}
\def\ka{\kappa}
\def\si{\sigma}
\def\Si{\Sigma}
\def\Ga{\Gamma}
\def\ze{\zeta}

\def\vph{\varphi}

\def\vep{\varepsilon}

\def\vth{{\vartheta}}

\def\tal{\tilde{\alpha}}
\def\tbe{\tilde{\beta}}
\def\tde{\tilde{\delta}}

\def\tV{\tilde{V}}

\def\tla{\tilde{\lambda}}
\def\tga{\tilde{\gamma}}
\def\tGa{\tilde{\Gamma}}

\def\tU{\tilde{U}}
\def\tw{\tilde w}
\def\tW{\tilde W}
\def\tB{\tilde B}

\def\tV{\tilde V}
\def\tz{\tilde z}
\def\tb{\tilde b}

\def\tG{\tilde G}

\def\tR{\tilde R}

\def\hH{\hat{H}}

\def\hY{\hat{Y}}

\def\hT{\hat{T}}

\def\hw{\hat{w}}
\def\hW{\hat{W}}
\def\hu{\hat{u}}

\def\hv{\hat{v}}

\def\hE{\widehat{E}}

\def\C{\mathbf{C}}
\def\I{\mathbf{I}}
\def\P{\mathbf{P}}
\def\Q{\mathbf{Q}}

\def\R{\mathbf{R}}
\def\N{\mathbf{N}}
\def\Z{\mathbf{Z}}

\def\F{\mathbf{F}}

\def\0{\mathbf{0}}
\def\H{\mathbf{H}}

\def\o{\mathcal{O}}

\def\r{\mathcal{R}}

\def\h{\mathcal{H}}

\def\y{\mathcal{Y}}
\def\e{\mathcal{E}}
\def\v{\mathcal{V}}
\def\z{\mathcal{Z}}
\def\x{\mathcal{X}}

\def\w{\mathcal{W}}
\def\i{\mathcal{I}}
\def\j{\mathcal{J}}

\def\lan{\langle}
\def\llb{(\!(}
\def\ran{\rangle}
\def\rrb{)\!)}

\def\dim{\mathop{\hbox{\rm dim}\,}_{\mathbf C}\,}

\def\lng{\hbox{\tiny lng}}
\def\sht{\hbox{\tiny sht}}
\def\sph{\hbox{\tiny sph}}

\def\JJ{\mathfrak{J}}
\def\HH{\mathfrak{H}}
\def\FF{\mathfrak{F}}

\def\CC{\mathfrak{C}}
\def\LL{\mathfrak{L}}

\def\AA{\mathfrak{A}}

\def\HH{\hbox{${\mathcal H}$\kern-5.2pt${\mathcal H}$}}

\font\smm=msbm10 at 12pt 
\def\symbol#1{\hbox{\smm #1}}
\def\lsmash{{\symbol n}}

\def\#{\sharp}

\title [Doble Hecke Algebras]{Double Affine Hecke Algebras\\
and Difference Fourier Transforms}
\author[Ivan Cherednik]{Ivan Cherednik $^\dag$} 
\date{November 2001}

\thanks{$^\dag$ \ Partially supported by NSF grant 
DMS-9877048,}

\address[I. Cherednik]{Department of Mathematics, UNC 
Chapel Hill, North Carolina 27599, USA\\
chered@math.unc.edu}

\begin{document}
\maketitle
{\small
\tableofcontents
}

% {0} Introduction \\ 
% {1} Affine Weyl groups 
%({\small Reduction modulo $W$}) \\ 
% {2} Double Hecke algebras 
%({\small Automorphisms, Demazure-Lusztig operators}) \\ 
% {3} Macdonald polynomials 
%({\small Intertwining operators}) \\ 
% {4} Fourier transform on polynomials 
%({\small Basic transforms}) \\ 
% {5} Jackson integrals 
%({\small Macdonald's $\eta$-identities}) \\ 
% {6} Semisimple representations  
%({\small Main Theorem, $GL_n$ and other applications}) \\ 
% {7} Spherical representations
%{\small (Semisimple spherical representations)} \\ 
% {8} Gaussian and self-duality
%{\small (Perfect representations, Main examples)}
%\vfill
%\intro %(Optional, Introduction)
%\vfil

We study difference Fourier transforms in the
representations of double affine Hecke algebras
in Laurent polynomials, polynomials multiplied by the
Gaussian, and in various spaces of delta-functions,
including finite-dimensional representations.
The theory generalizes the Harish-Chandra theory of
spherical functions, its $p$-adic variant, and harmonic
analysis on Heisenberg and Weyl algebras, predecessors
of the double Hecke algebras. 
There are applications to the Harish-Chandra inversion, 
Verlinde algebras, Gauss-Selberg sums, and
Macdonald's $\eta$-type identities. 

We discuss two major directions, compact and non compact, 
generalizing the corresponding parts of the
classical theory of spherical functions. 
They are based on the imaginary and real integrations 
and their substitutes, namely, the constant term 
functional and the Jackson summation. 
Both directions exist in two different variants: 
real, as $q<1,$ and unimodular for $|q|=1,$
including roots of unity. At roots of unity,
compact and non compact directions coincide and the 
irreducible representations are all finite-dimensional.

It is important to extend the difference 
Fourier transform to arbitrary irreducible
semisimple representations of double Hecke algebras.
Their description  is reduced in the paper to a certain
combinatorial problem, which can be managed in several cases,
including the case of $GL_n$ for generic $q$ as
$k,$ the other structural constant of 
the double Hecke algebra, is a special rational number. 

For such algebras, 
semisimple and general irreducible representations 
are described in terms of periodic skew and generalized Young
diagrams using the technique of intertwiners. The classification
generalizes that due to Bernstein-Zelevinsky 
in the affine case and the papers [C9,C11,Na].
For generic $k,$ the irreducible representations are
induced from the affine Hecke algebra and the theory 
is not too promising. Recently Vasserot generalized 
the "geometric approach"
from [KL1] to the double affine case. His results are
for arbitrary root systems and are also the most fruitful
in the case of rational $k.$

In the paper, we focus on irreducible 
finite-dimensional semisimple unitary (or pseudo-unitary)
spherical self-dual (Fourier-invariant) representations, 
which are called perfect. 
They exist either as $q$ is a root of unity or
when $q$ is generic but $k$ is special rational.
Each double Hecke algebra (for an arbitrary
reduced root system) may have only one such 
representation up to isomorphisms and the choice
of the character of the affine Hecke algebra (necessary
in the definition of spherical representations).
It is a generalization of the uniqueness of the 
irreducible 
representations of the Weyl algebras at roots of unity.
In the $A$-case, such representations are also unique 
among
all finite-dimensional representations of the 
corresponding double Hecke algebra.

As $q$ is a root of unity, perfect representations
lead to a new class of Gauss-Selberg sums, the eigenvalues
of the Gaussian which is an eigenfunction of the 
difference Fourier transform. Similar sums for
generic $q$ and special rational $k$ are directly related
to Macdonald's $\eta$-type identities and to the classification
of degenerations of the multi-dimensional Bessel and 
hypergeometric functions. 

\vfil
Arbitrary reduced root systems are considered in the paper.
Modern methods are good enough to transfer smoothly
the one-dimensional theory to such setting.
It is worth mentioning that there are quite a few new 
results even in the rank one case, including new 
identities and new proofs of the classical formulas. 
It seems
that the Fourier transform of the considered type was
never studied systematically before [C3], even in the
one-dimensional case. Its kernel, which is essentially
the basic hypergeometric function, of course appeared in
the integrals. However, as far as I know, the self-duality
of the corresponding Fourier transform and the 
Fourier-invariance of the Gaussian remained unknown.  

Let me also mention recent paper [St] 
(see also ref. therein)
and the manuscript by Macdonald devoted to the so-called
$C^\vee C$ case, which has all important 
features of the case of 
reduced root systems considered in this paper, including the self-duality
and the existence of the Gaussians. The paper and the manuscript
contain reasonably complete algebraic theory of the 
double affine Hecke algebra of type $C^\vee C,$
nonsymmetric polynomials, and the Fourier transforms, 
including the prior results due to Koornwinder, Noumi, 
and Sahi.

Note a challenging parallelism between 
our cyclotomic Gaussian sums and modular 
Gauss-Selberg sums (see, for instance, [E]).
However the difference is dramatic. 
In the modular case, Selberg-type kernels are calculated 
in  finite
fields and are  embedded into roots of unity right 
before the summation. 
Our sums are defined entirely in cyclotomic fields.

\vfil
%\vskip 0.2cm
{\it Nonsymmetric polynomials.}
The key tool of the recent progress in the theory of orthogonal
polynomials, related combinatorics, and harmonic 
analysis is
the technique of nonsymmetric Opdam-Macdonald polynomials.
The main references are [O2,M4,C4]. Opdam mentions 
in [O2] that
a definition of the nonsymmetric polynomials 
(in the differential 
setup) was given in Heckman's unpublished lectures.  
These polynomials are expected to be, generally speaking,
beyond quantum groups and Kac-Moody algebras because of
the following metamathematical reason. 
Many special functions in the Lie and Kac-Moody
theory, including classical and affine characters,
spherical functions, conformal blocks and more 
are $W$-invariant. One can expect these functions 
to be $W$-symmetrizations of simpler and hopefully more
fundamental nonsymmetric functions, but it doesn't happen
in the traditional theory, with some reservation 
about the Demazure
character formula. One needs double Hecke algebras and 
nonsymmetric Macdonald's polynomials to manage this problem. 

However our considerations are not quite new in
the representation theory. We directly generalize the
harmonic analysis on Heisenberg and Weyl algebras, 
theory of metaplectic (Weil) representations, 
and borrow a lot of the affine Hecke algebra technique.
The $p$-adic counterparts of the nonsymmetric 
polynomials are the limit of our polynomials
as $q\to \infty$ and they are well known. The book
[Ma] is devoted to them. See [O4] about recent developments. 
Let me also mention the Eisenstein integrals and series, 
which are nonsymmetric. 
The classical number theory does require
nonsymmetric functions. For instance, we need all 
one-dimensional theta functions, not only even ones.

There are important applications of the new technique.
Let me mention first [O2] and its certain continuation [C8],
where the Harish-Chandra transforms of the coordinates
treated as multiplication operators were calculated,
which was an old open problem. This technique was
successfully applied in [KnS] and further papers to the 
combinatorial conjectures about the coefficients of 
Macdonald's polynomials. The third and, I think, very convincing 
demonstration is the theory of nonsymmetric 
Verlinde algebras.  

\vfil
%\vskip 0.2cm
{\it Nonsymmetric Verlinde algebras.}
Double Hecke algebras generalize Hei\-sen\-berg and
Weyl algebras, so it is not surprising that
they are very helpful to study Fou\-rier transforms.
In a sense, perfect representations constitute all 
commutative algebras with "perfect" Fourier transforms,
i.e. enjoying all properties of the classical one. 
Spherical representations appear here because 
they have natural structures of commutative  algebras. 

The subalgebras of symmetric ($W$-invariant) elements
generalize the Verlinde algebras. The latter are formed by
integrable representations of Kac-Moody algebras of fixed level
with the celebrated  fusion procedure as the multiplication.  
They can be also introduced as restricted
categories of representations of quantum groups at roots of
unity due to Kazhdan and Lusztig [KL2] and Finkelberg.

The perfect representations "triple-generalize" the 
Verlinde algebras, and what is important, dramatically 
simplify the theory. First, the Macdonald
symmetric polynomials replaced the characters [Ki,C3].
Technically, the parameter $k,$ which is 
$1$ in the Verlinde theory [Ve], 
became an arbitrary positive integer (satisfying
certain inequalities).  Second and very
important, the nonsymmetric polynomials replaced the 
symmetric ones [C4]. It is much more comfortable 
to deal with the Fourier transform when 
all functions are available, not only symmetric.
It is exactly what perfect representations provide.
Third, fractional $k$ appeared.
Generally speaking,  $hk\in \Z$ for the
Coxeter number $h,$ which is one of the results of
the paper (Sections 5,8).
 
The final construction is rather distant from 
the original Verlinde algebra (and from the Lie theory). 
The structural constants of the multiplication are not
positive integers anymore. However all other important
features, including the action of $PGL_2(\Z)$ and
the positivity of the inner product, were saved. The main 
hermitian form considered in the paper 
(and previous works) extends that due 
to Verlinde. It makes the generators $X,T,Y$ of the 
double Hecke algebra unitary and sends $q\mapsto q^{-1}.$
It is positive as $|q|=1$ and the angle of $q$ is sufficiently small.

%\vfil
%\vskip 0.2cm 
There is another important inner product which does not 
require the conjugation of $q$ and leads to a more traditional 
variant of the double Hecke harmonic analysis. It is also 
considered in the paper. In the standard representations 
in $\de$-functions, the operators $Y$ are not normal with 
respect to this pairing, which diminish the role 
of the $Y$-semisimplicity in the theory.
It is similar to the regular representation of the affine 
Hecke algebra, where these operators are due 
to Bernstein and
Zelevinsky. This representation can be obtained as 
a limit of the $\de$-representation as $q\to \infty.$ 
The $X$-operators are selfadjoint with respect to 
this pairing and are very helpful in the double Hecke theory.
However they do not survive under this limit. 
In the $p$-adic harmonic
analysis, one needs the complete classification and
advanced methods like $K$-theory [KL1]. It is what can
be expected for the pairing without $q\to q^{-1}$ for the
double Hecke algebras. 

In the unimodular
case with the pairing involving the $q$-conjugation, 
reasonable unitary representations are $Y$-semisimple 
and their eigenfunctions, for instance the nonsymmetric
polynomials, play a very important role. Generally
speaking, the theory becomes
more elementary, and there is a hope to manage it using 
classical tools like the technique of the
highest vectors. The paper is based on this approach. 
Note that in the Verlinde case,  
the $W$-symmetrizations of the $Y$-eigenfunctions
are the restrictions of the classical characters to roots
of unity.

\vfil
%\vskip 0.2cm
{\it Macdonald-Mehta integrals.}
Selberg's integrals were actually the starting
point. By the way, they also gave birth to
the theory of symmetric Macdonald polynomials,
which appeared for the first time in Kadell's
unpublished work. The history is as follows. 
Mehta suggested a formula for the integral of 
the product $\prod_{1\le i< j\le n} (x_i-x_j)^{2k}$
with respect to the Gaussian measure. 
The formula was readily deduced
from Selberg's integral by Bombieri.
Macdonald extended Mehta's integral in [M1]
to arbitrary root systems and verified his conjectures 
for the classical systems. He employed Selberg's 
integrals too.
Later Opdam found  a uniform proof for all root systems
using the shift operators in [O1]. Note that neither
alternative proofs without shift operators 
are known for special root systems.
In [C5], $q$ was added to this theory. 

The  (classical) Macdonald-Mehta  integral is the 
normalization constant
for the generalized Hankel transform introduced
by Dunkl [D]. 
The generalized Bessel functions [O3]
multiplied by the Gaussian are
eigenfunctions of this transform. The eigenvalues are 
given in terms of this constant [D,J]. It is the same 
in the $q$-case.

The Hankel transform 
is a "rational" degeneration of the Harish-Chandra 
spherical transform. The symmetric space $G/K$
is replaced by its tangent space $T_e(G/K)$ with the adjoint
action of $G$
(see [H]). Here $G$ is a semi-simple Lie group, 
$K$ its maximal compact
subgroup. Only very special  $k,$ the root multiplicities,
may appear in the Harish-Chandra theory. For instance, 
$k=1$ is
the so-called group case. When the Harish-Chandra transform
is replaced by its generalization in terms
of the hypergeometric functions (see [HO1]), the 
parameter $k$ becomes arbitrary (and extra parameters appear).

\vfil
\vskip 0.2cm
{\it From Gauss integrals to Gaussian sums.}
Study of the $q$-counterparts of the Hankel and 
Harish-Chandra transforms and the passage to
the roots of unity are the goals  of the paper.
There are four main steps:  \\ 
 {a)} the $q$-deformation of the Macdonald-Mehta 
integrals,  \\ 
 {b)} the transfer from the $q$-integrals to the Jackson sums,  \\ 
 {c)} their counterparts/limits at roots of unity, \\ 
 {d)} the product formulas for the Gauss-Selberg sums.

Let us demonstrate the complete procedure restricting 
ourselves with the classical one-dimensional Gauss integral.

We want to go from the celebrated formula
\begin{align}
&\int_{-\infty}^\infty e^{-x^2}x^{2k}\hbox{d}x\ =\ 
\Gamma (k+1/2), \
\Re k>-1/2, 
\label{(1)}
\end{align}
to the almost equally famous
Gauss formula
\begin{align}
&\sum_{m=0}^{2N-1} e^{\frac{\pi m^2 }{ 2N}i}\ =\ 
(1+i)\sqrt{N},\ \ N\in {\mathbf N}.
\label{(2)}
\end{align}
The actual objective is to incorporate $k$ into the last formula.
It seems that this natural problem was never 
considered in the classical works on Gaussian sums.
Modular Gauss-Selberg sums give some solution,
but it is not what one may expect from the viewpoint
of the consistent $q$-theory.

Obviously a trigonometric counterpart of 
(\ref{(1)}) is needed. 
A natural candidate is the Harish-Chandra spherical transform,
where $\sinh(x)^{2k}$ substitutes for  $x^{2k}$. 
However such choice creates problems. 
The Harish-Chandra transform is not self-dual anymore,
the Gaussian looses its Fourier-invariance, and
the formula (\ref{(1)}) has no $\sinh$-counterpart 
(at least for generic $k$).

\vfil
a) {\it Difference setup.} It was demonstrated recently that these 
important
features of the classical Fourier transform are restored for the 
kernel
\begin{align}
&\delta_k(x;q)\equal\prod_{j=0}^\infty \frac{(1-q^{j+2x})(1-q^{j-2x})
}{
(1-q^{j+k+2x})(1-q^{j+k-2x})},\ 0<q<1, \ k\in \mathbf{C}.
\label{(3)}
\end{align}

Actually the self-duality of the corresponding transform
can be expected a priori because the Macdonald truncated 
theta-function $\delta$ is a unification of
$\sinh(x)^{2k}$ and the Harish-Chandra function ($A_1$) serving 
the inverse spherical transform.

Setting $q=\exp(-1/a),\ a>0,$
\begin{align}
&(-i)\int_{-\infty i}^{\infty i}q^{-x^2}
\delta_k\hbox{\, d}x=
{2\sqrt{a\pi}}\prod_{j=0}^\infty
\frac{1-q^{j+k}}{
 1-q^{j+2k}},\ \ \Re k>0.
\label{(4)}
\end{align}
The limit of (\ref{(4)})  multiplied by $(a/4)^{k-1/2}$
as $a\to \infty$ is (\ref{(1)}) in the imaginary variant.

\vfil
%\vskip 0.2cm
b) {\it Jackson sums.}
Special  $q$-functions have many interesting properties 
which have no classical counterparts. The most promising
feature is a possibility to replace the integrals by sums,
the Jackson integrals. 

Let $\int_\#$ be the integration for the 
 path which  begins at $z=\epsilon i+\infty$, moves
to the left till $\epsilon i$, then down
through  the origin to  $-\epsilon i$, and finally
returns down the positive real axis to $-\epsilon 
i+\infty$
(for small $\epsilon>0$).

Then for $|\Im k|<2\epsilon, \Re k>0,$
\begin{align}
&\frac{1}{ 2i}\int_{\#} q^{x^2} \delta_k\hbox{\, d}x\ =
-\frac{a\pi}{ 2}\prod_{j=0}^\infty
\frac{(1-q^{j+k})(1-q^{j-k})}{
 (1-q^{j+2k})(1-q^{j+1})}\times\lan q^{x^2}\ran_\#,
\notag\\ 
&\lan q^{x^2}\ran_\#\equal
\sum_{j=0}^\infty q^{\frac{(k-j)^2}{ 4}}
\frac{1-q^{j+k}}{
1-q^{k}} \prod_{l=1}^j 
\frac{1-q^{l+2k-1}}{ 
 1-q^{l}}\notag
\end{align}
\begin{align} 
& =\ \prod_{j=1}^\infty
\frac{(1-q^{j+k})}{(1-q^{j})} 
\sum_{j=-\infty}^{\infty} q^{\frac{(k-j)^2
}{ 4}}\ =
\label{(5)}\\ 
& q^{\frac{k^2}{ 4}}\prod_{j=1}^\infty
\frac{(1-q^{j/2})(1-q^{j+k})(1+q^{j/2-1/4+k/2})
(1+q^{j/2-1/4-k/2})
}{ (1-q^j)}.
\label{(6)}
\end{align}
The sum for $\lan q^{x^2}\ran_\#$
is the Jackson integral for a special choice ($-k/2$) 
of the starting point. Its convergence is for all $k$. 
The transfer (\ref{(5)}) $\Rightarrow$ (\ref{(6)}) holds 
for arbitrary root systems 
and, generally speaking, requires a variation of the 
starting point (and representations of the double Hecke
algebra with highest weights). 
However in the one-dimensional case we can simply use
the classical $\eta-$type identity.

\vfil
%\vskip 0.2cm
c) {\it Gaussian sums.} 
When $q=\exp(2\pi i/N)$ and $k$ is a 
positive integer $\le N/2$ we come to the 
simplest Gauss-Selberg cyclotomic sum:
\begin{align}
&\sum_{j=0}^{N-2k} q^{\frac{(k-j)^2}{ 4}}
\frac{1-q^{j+k}}{
1-q^{k}} \prod_{l=1}^j
\frac{1-q^{l+2k-1}}{
 1-q^{l}}=
\prod_{j=1}^k
(1-q^{j})^{-1}\sum_{m=0}^{2N-1} q^{m^2/4}.
\label{(7)}
\end{align}
Its modular counterpart is (1,2b) from  [E].

Formula (\ref{(7)}) can be deduced directly 
from (\ref{(5)}) following the
classical limiting  procedure from [Ch]. 
Formulas (\ref{(4)})-(\ref{(6)}) 
can be verified by elementary
methods too. However the generalizations of 
these formulas involving 
$q$-ultraspherical Rogers' polynomials
do require double Hecke algebras and lead to new 
one-dimensional identities. 
Their nonsymmetric variants
and multidimensional generalizations are based entirely
on Hecke algebras.

There are other examples of Gaussian sums, for instance,
those used in the classical  quadratic reciprocity. 
All of them  were interpreted using ``the smallest''
representations of the rank one double Hecke algebra
in [C7]. It seems that the formulas with integers 
$k>1$ did not appear in arithmetic. 
There are also interesting new formulas for 
half-integral $k$. 

\vfil
%\vskip 0.2cm
d) {\it New proof of the Gauss formula.} 
Substituting $k=[N/2]$ (not  $k=0,$ as one may expect,
which makes (\ref{(7)}) a trivial
identity) we get the product formula for  
$\sum_{m=0}^{2N-1} q^{m^2/4},$ which can be readily calculated
as $q=e^{2\pi i/N}$ and quickly  results in (\ref{(2)}).

Let us consider  the case of $N=2k$ only 
(odd $N=2k+1$ are quite similar):
\begin{align}
&\sum_{m=0}^{2N-1} q^\frac{m^2}{ 4}\ =\ 
q^{k^2}{ 4}\,
\Pi \for
\Pi=\prod_{j=1}^k (1-q^{j}).
\notag \end{align}

First, $\Pi\bar{\Pi}=(X^{N}-1)(X+1)(X-1)^{-1}(1)=2N.$
Second, $\arg(1-e^{i\phi})=\phi/2-\pi/2$ when 
$0< \phi< 2\pi,$ and therefore 
$$\arg\Pi = \frac{\pi}{ N} {k(k+1)}{ 2}-{\pi k}{ 2} =
\frac{\pi (1-k)}{ 4}.$$

Here we can switch to arbitrary primitive
$q=\exp(2\pi i l/N)$ as $(l,N)=1.$ It is necessary 
to control somehow the set of $\arg(q^{j})$ for 
$1\le j\le k.$ This leads
to the quadratic reciprocity (see [CO]). We do not discuss
this direction in this paper.

Finally, 
$$\arg( q^\frac{k^2}{ 4}\Pi)={\pi (1-k)}{ 4}+
\frac{\pi k}{ 4}={\pi}{ 4},\and 
 q^\frac{k^2}{ 4}\Pi=\sqrt{N}(1+i).
$$

\vfil
\vskip 0.2cm
{\it Macdonald's $\eta$-identities.}
The approach via double Hecke algebras provides  
interesting deformations 
of (some) cyclotomic Gauss-Selberg sums 
and the corresponding unitary self-dual representations. 
Surprisingly, the roots of unity
can be replaced by arbitrary $q$
such that $|q|=1.$ The
numbers of terms in the sums and 
respectively the dimensions of the representations 
remain the same. The Fourier transform and all other 
related structures and formulas
can be deformed as well. For instance, as $k=1$ we
deform the Verlinde algebras together with the hermitian 
forms and the projective $PGL_2(\Z)$-action.
%These deformations might reflect some properties
%(filtrations?) of the integrable representations 
%of the Kac-Moody algebras.  

The deformation is as follows.
We start with finite-dimensional irreducible
spherical representations of double affine Hecke algebras
where $q$ is generic and $k$ is 
special negative rational. Such representations
are studied in the paper in detail 
(note that an example was considered in [DS]). 
They remain irreducible  as $q$ becomes a
root of unity under certain conditions.
Since positive and negative $k$ 
cannot be distinguished at roots of unity,
it gives the desired deformation.

%This relation resembles to a certain 
%extent the correspondence 
%between representations of Kac-Moody algebras of negative 
%integral levels and quantum groups at roots of unity.

The representation theory of double Hecke
algebras at such special negative
$k$ is directly connected with the 
Macdonald identities [M6].
They are product formulas for the Gauss-type sums over 
the root lattices upon the shift by $-(1/h)\rho,$ where 
$\rho$ is the half-sum of all positive roots, $h$ the
Coxeter number. The Macdonald identities correspond 
exactly to one-dimensional representations of 
the  double Hecke algebra.
Replacing $-(1/h)$ by $-(m/h)\rho$ for integral  $m>0$ relatively
prime to  $h,$ we come to multidimensional
representations.

\vfil
\vskip 0.2cm
{\it Relations to the classical theories}.
There are two variants of difference Fourier transforms
based on the pairing a) involving or b) not involving the
formal conjugation $q\mapsto q^{-1},\, t\mapsto t^{-1}.$
Case b) is new. The previous papers [M4,C4] and 
others devoted to the nonsymmetric polynomials were based
on the pairing of type a). This conjugation is quite natural 
as $|q|=1,$ but creates certain 
problems in the real theory. In the paper, we 
remove the conjugation, for instance,  calculate the
scalar products of the nonsymmetric polynomials without
it. The final theory is actually very close to that from
[O2]. See also [HO1]. In a way, we  replace $w_0$ in Opdam's 
paper by $T_{w_0}.$ 

To be more exact, the conjugation $q\to q^{-1}$ in the
main scalar product is not a problem and can be eliminated
when $W$-invariant functions are considered. However it
plays an important role in the nonsymmetric setup.

This development clarifies completely the relations
to the Harish-Chandra theory of spherical functions and
the $p$-adic theory due to Macdonald, Matsumoto
and others. 
In papers [O4,O5], Opdam developed the Matsumoto
theory of "nonsymmetric" spherical functions towards the
theory of nonsymmetric polynomials. The operator $T_{w_0}$
appears there in the inverse transform in a way 
similar to that of the present paper 
([O4], Proposition 1.12). 
The $p$-adic theory corresponds to
the limit $q\to \infty,$ with $t$ being $p$ or its proper 
power (in the rigorous $p$-adic setup). Technically, 
the case $q=\infty$ is much simpler than the $q,t$-case, 
however quite a few features of the general case can be seen 
under such degeneration. 
Note that the self-duality of the Fourier transform 
and the Gaussian do not survive 
in this limit. See the end of [CO]. 

It is worse mentioning  that in
the theory of double affine Hecke algebras 
the Fourier transform has a topological interpretation
as the transposition of the periods of the elliptic
curve via the elliptic braid group. See [CO] and
also [Bi].

As $q\to 1,\, t=q^k,\ $ we come to 
[O2] and, via the symmetrization, to the Harish-Chandra theory 
(with $k$ being the root multiplicity). Algebraically, this 
limiting case is similar to the 
general theory, however the analytic problems
are quite different. The nonsymmetric Fourier transform in 
the space of compactly supported $C^\infty$-functions is analyzed 
in [O2] in detail. Analytically, the $q,t$-case seems closer to the 
$p$-adic theory and to the theory of generalized Hankel transform.
We do not discuss the analytic problems
in this paper. There are some results in this direction
in [C5] (mainly the construction of the general spherical functions), 
[C10] (analytic continuations in terms of 
$k$ with applications to $q$-counterparts of Riemann's zeta 
function), [KS] and in other papers by Koelink and Stokman.
The latter papers are devoted to explicit 
analytic properties of the one-dimensional
Fourier transforms in terms of the basic hypergeometric 
function, closely connected with those considered in [C3,C5]. 

\vfil
\vskip 0.2cm
{\it Plan of the paper.}
Let us describe the structure of the paper.
The first three sections mainly
contain the results from the previous papers,
adjusted and extended to include
a) the scalar products
without the conjugation $q\mapsto q^{-1},$
b) the roots of unity. There are also some 
new combinatorial results about affine Weyl groups
necessary for the classification technique from the
next sections, which is essentially the technique
of intertwining operators.

Sections 4,5 are about general theory
of the Fourier transforms in the compact and non compact
cases. They act on polynomials with and without 
multiplication by the Gaussian and in spaces of  
delta-functions (called the functional representations).
The application to the $\eta$-identities concludes
Section 5. 

In Section 6, we find necessary and
sufficient conditions for a representation with the 
cyclic vector to be semisimple and pseudo-unitary,
including a special consideration of the case of $GL_n.$
It is the main tool for the next 
Sections 7,8  devoted to the classification 
of spherical  and  self-dual representations.
We consider there main examples of 
finite-dimensional prefect representations, which provide
generalizations of the Macdonald $\eta$-identities 
and the Verlinde algebras.  

The exposition is self-contained. The 
one-dimen\-si\-onal papers  
[CM], [C7], and [CO]
combined with the second part of [C6]
could be a reasonable introduction.
The papers [C5] (Macdonald-Mehta integrals,
general $q$-spherical
functions) and  [C1] (induced and spherical
representations) also contain additional results and
examples.

%\vfil
\vskip 0.2cm
{\it Acknowledgments.}
The author thanks V.~Kac and I.~Macdonald for useful 
comments, P.~Etingof and E.~Vasserot for discussing the
case of $GL_n,$ and acknowledges his indebtedness
to E.~Opdam for reading the paper and 
making important helpful suggestions which especially 
influenced Introduction and Section 7.
Significant part of the paper
was written while delivering a course of lectures
on double Hecke algebras 
at Harvard University (spring 2001).
I am grateful to D.~Kazhdan,
P.~Etingof, and V.~Ostrik, who made the course possible,
participated in it, and stimulated the paper a great deal. 
I also would like to thank R.~Rentschler
for the kind invitation (summer 2001), and IHES for 
hospitality (summer 2002). 

%\vfill
\vskip 0.2cm
\section{Affine Weyl groups}
\setcounter{equation}{0}

Let $R=\{\al\}   \subset \R^n$ be a root system of type $A,B,...,F,G$
with respect to a euclidean form $(z,z')$ on $\R^n 
\ni z,z'$,
$W$ the Weyl group  generated by the reflections $s_\al$,
$R_{+}$ the set of positive  roots ($R_-=-R_+$), 
corresponding to (fixed) simple roots  
roots $\al_1,...,\al_n,$ 
$\Ga$ the Dynkin diagram  
with $\{\al_i, 1 \le i \le n\}$ as the vertices, 
$R^\vee=\{\al^\vee =2\al/(\al,\al)\}$ the dual root system,
\begin{align}
& Q=\oplus^n_{i=1}\Z \al_i \subset P=\oplus^n_{i=1}\Z \om_i, 
\notag \end{align}
where $\{\om_i\}$ are fundamental weights:
$ (\om_i,\al_j^\vee)=\de_{ij}$ for the 
simple coroots $\al_i^\vee.$

Replacing $\Z$ by $\Z_{\pm}=\{m\in\Z, \pm m\ge 0\}$ we get
$Q_\pm, P_\pm.$
Note that $Q\cap P_+\subset Q_+.$ Moreover, each $\om_j$ has all 
nonzero coefficients (sometimes rational) when expressed 
in terms of 
$\{\al_i\}.$
Here and further see [B].  

The form will be normalized
by the condition  $(\al,\al)=2$ for the 
short roots. This normalization coincides with that
from the tables in [B]  for $A,C,D,E,G.$
Hence $\nu_\al\equal (\al,\al)/2$ can be either $1,2$ or $1,3$ and 
$Q\subset Q^\vee,  P\subset P^\vee,$ where $P^\vee$ is generated by
the fundamental coweights $\om_i^\vee.$ 
Sometimes we write $\nu_{\sht}$ for short roots (it is always $1$)
and  $\nu_{\lng}$ for long ones.

%\vfil
Let  $\vth\in R^\vee $ be the maximal positive 
coroot. All simple coroots appear in its
decomposition. 
See [B] to check that  $\vth$ considered
as a root 
(it belongs to $R$ because of the choice of normalization)
is maximal among all short positive roots of $R.$
For the sake of completeness, let us prove another defining property 
of $\vth$ (see Proposition \ref{BSTIL} below for a  uniform proof).

\begin{proposition}\label{LEAST}
The least nonzero element in $Q_+^{+}=Q_+\cap 
P_+=Q\cap P_+$ 
with respect to $Q_+$ is
$\vth,$ i.e. $b-\vth\in Q_+$ for all $b\in Q_+^{+}.$  
\end{proposition}

{\it Proof.} If $a\in Q_+^{+}$ then all coefficients of the 
decomposition of $a$ in terms of $\al_i$ are nonzero. Indeed,
if $\al_i$ does not appear in this decomposition then its
neighbors  
(in the Dynkin diagram) have negative scalar products with $a.$
This readily gives the claim when 
$\vth=\sum_{i=1}^n\al_i,$
i.e. for the systems $A,B.$ In these cases,  $\vth$ is 
$\om_1+\om_n\for  A$ and $\om_1 \for B.$ Otherwise, it is $\om_2\ 
(C,D,E_6)$,
$\om_1\ (E_7),$ $\om_8\ (E_8),$ $\om_4\ (F_4),$ and
$\om_1\ (G_2)$ respectively (in the notation from [B]). 
The corresponding subscripts of $\om$ will be denoted by $\tilde{o}.$
They are the indices of the simple roots
neighboring the root $-\vth$ added to the Dynkin diagram $\Ga$.
 
The proposition holds for $E_8,F_4,G_2,$ because 
in these cases
$P=Q$ and $ \om_j-\om_{\tilde{o}} \in Q_+$ 
(use the tables from
[B]). So we need to examine $C,D,E_{6,7}.$

As to $C,$ one verifies that $ \om_j-\om_{2} \in 
(1/2)Q_+$ for 
$j\ge 2$ and 
$Q_+^{+}$ is generated by $\om_{2i}$ and 
$\om_{2i+1}+\om_{2j+1}$
for $1\le i\le j\le [n/2].$ Thus the relation  
$2\om_1-\om_2\in Q_+$
proves the claim. 
In the $D$-case, $a-\om_2$ for $a\in Q_+^{+}$
can be apart from $(1/4)Q_+$  
only if it is a linear combination of
$\om_1,\om_{n-1},$ and $\om_n.$ However such a 
combination either
``contains'' $\om_{n-1}+\om_n$ or $2\om_1,$ which is 
sufficient to
make it greater than  $\om_2$
(with respect to $Q_+$). 

For $R=E_6,$  $ \ \om_j-\om_{2} \in (1/3)Q_+$ when 
$j\ge 2=\tilde{o}$ except for $j=6.$ The intersection
of $Q$ and $\Z_+\om_1+\Z_+\om_6$ is generated
by $3\om_{1},$ $3\om_{6},$ $\om_1+\om_6.$ Either weight
is greater than $\om_2.$
In the case of $E_7,$ $\tilde{o}=1,$
$ \om_j-\om_{1} \in (1/2)Q_+$ if $j\le 6,$ and  
$2\om_7$ is
greater than $\om_1.$ 
\qed 

Setting
$\nu_i\ =\ \nu_{\al_i}, \ 
\nu_R\ = \{\nu_{\al}, \al\in R\},$ one has 
\begin{align}
&\rho_\nu\equal (1/2)\sum_{\nu_{\al}=\nu} \al \ =
\ \sum_{\nu_i=\nu}  \om_i, \hbox{\ where\ } \al\in R_+, 
\ \nu\in\nu_R.
\end{align}

\vskip 0.2cm
{\it Affine roots.}
The vectors $\ \tal=[\al,\nu_\al j] \in 
\R^n\times \R \subset \R^{n+1}$ 
for $\al \in R, j \in \Z $ form the 
{\it affine root system} 
$\tR \supset R$ ( $z\in \R^n$ are identified with $ [z,0]$).  
We add $\al_0 \equal [-\vth,1]$ to the simple
 roots for the 
maximal short root $\vth$.
The corresponding set $\tR$ of positive roots coincides with
$R_+\cup \{[\al,\nu_\al j],\ \al\in R, \ j > 0\}$.

We complete the Dynkin diagram $\Ga$ of $R$  
by $\al_0$ (by $-\vth$ to be more
exact). The notation is $\tGa$. One can get it from the
completed Dynkin diagram for $R^\vee$ [B]
reversing all the arrows.

The set of
the indices of the images of $\al_0$ by all 
the automorphisms of $\tGa$ will be denoted by $O$ 
($O=\{0\} \for E_8,F_4,G_2$). Let $O'={r\in O, r\neq 0}$.
The elements $\om_r$ for $r\in O'$ are the so-called minuscule
weights: $(\om_r,\al^\vee)\le 1$ for
$\al \in R_+$.

Given $\tal=[\al,\nu_\al j]\in \tR,  \ b \in B$, let  
\begin{align}
&s_{\tal}(\tz)\ =\  \tz-(z,\al^\vee)\tal,\ 
\ b'(\tz)\ =\ [z,\ze-(z,b)]
\label{ondon}
\end{align}
for $\tz=[z,\ze] \in \R^{n+1}$.

The {\it affine Weyl group} $\tW$ is generated by all $s_{\tal}$
(we write $\tW = <s_{\tal}, \tal\in \tR_+>)$. One can take
the simple reflections $s_i=s_{\al_i}\ (0 \le i \le n)$ 
as its
generators and introduce the corresponding notion of the  
length. This group is
the semidirect product $W\lsmash Q'$ of 
its subgroups $W=<s_\al,
\al \in R_+>$ and $Q'=\{a', a\in Q\}$, where
\begin{align}
& \al'=\ s_{\al}s_{[\al,\nu_{\al}]}=\ 
s_{[-\al,\nu_\al]}s_{\al}\for 
\al\in R.
\label{ondtwo}
\end{align}

The {\it extended Weyl group} $ \hW$ generated by $W\and P'$
(instead of $Q'$) is isomorphic to $W\lsmash P'$:
\begin{align}
&(wb')([z,\ze])\ =\ [w(z),\ze-(z,b)] \for w\in W, b\in B.
\label{ondthr}
\end{align}
Later in this paper,  $b$ and $b'$ will be identified.

Given $b\in P_+$, let $w^b_0$ be the longest element
in the subgroup $W_0^{b}\subset W$ of the elements
preserving $b$. This subgroup is generated by simple 
reflections. We set
\begin{align}
&u_{b} = w_0w^b_0  \in  W,\ \pi_{b} =
b( u_{b})^{-1}
\ \in \ \hW, \  u_i= u_{\om_i},\pi_i=\pi_{\om_i},
\label{wo}
\end{align}
where $w_0$ is the longest element in $W,$
$1\le i\le n.$

The elements $\pi_r\equal\pi_{\om_r}, r \in O'$ and
$\pi_0=\hbox{id}$ leave $\tGa$ invariant 
and form a group denoted by $\Pi$, 
 which is isomorphic to $P/Q$ by the natural 
projection $\{\om_r \mapsto \pi_r\}$. As to $\{ u_r\}$,
they preserve the set $\{-\vth,\al_i, i>0\}$.
The relations $\pi_r(\al_0)= \al_r= ( u_r)^{-1}(-\vth)$ 
distinguish the
indices $r \in O'$. Moreover (see e.g. [C2]):
\begin{align}
& \hW  = \Pi \lsmash \tW, \where
  \pi_rs_i\pi_r^{-1}  =  s_j \iif \pi_r(\al_i)=\al_j,\ 
 0\le j\le n.
\end{align}

\vskip 0.2cm
{\it The length on} $\hW.$
Setting
$\hw = \pi\tw \in \hW,\ \pi\in \Pi, \tw\in \tW,$
the length $l(\hw)$ 
is by definition the length of the reduced decomposition 
$\tw\ =\ s_{i_l}...s_{i_2} s_{i_1} $
in terms of the simple reflections 
$s_i, 0\le i\le n.$ The number of  $s_{i}$
in this decomposition 
such that $\nu_i=\nu$ is denoted by   $l_\nu(\hw).$

The length can be also introduced as the 
cardinality $|\la(\hw)|$
of  
$$
\la(\hw)\equal\tR_+\cap \hw^{-1}(\tR_-),\ \hw\in \hW.
$$ 
This set is 
a disjoint union of $\la_\nu(\hw):$
\begin{align}
&\la_\nu(\hw)\  
=\ \{\tal\in \la(\hw),\nu({\tal})=\nu \}\ =\ 
\{\tal, \ l_\nu( \hw s_{\tal}) < l_\nu(\hw) \}.
\label{lambda}
\end{align}

The coincidence with the previous definition 
is based on the equivalence of the relation
$l_\nu(\hw\hu)=
l_\nu(\hw)+l_\nu(\hu)$  
for
$\hw,\hu\in\hW$
and
\begin{align}
&  \la_\nu(\hw\hu) = \la_\nu(\hu) \cup
\hu^{-1}(\la_\nu(\hw)), 
\label{ltutw} 
\end{align}
which, in its turn, is equivalent to  
the positivity condition $\hu^{-1}(\la_\nu(\hw))
\subset \tR_+.$

Formula (\ref{ltutw}) includes obviously 
the positivity condition. It also
implies that
$$
\la_\nu(\hu) \cap \hu^{-1}(\la_\nu(\hw))\ =\ 
\hu^{-1}\bigr(\hu(\la_\nu(\hu)) \cap \la_\nu(\hw)\bigl)\ 
=\ \emptyset,
$$ thanks to the formula  
$$
\la_\nu(\hw^{-1}) = -\hw(\la_\nu(\hw)).
$$
Thus it results in the equality $l_\nu(\hw\hu)=
l_\nu(\hw)+l_\nu(\hu).$

The other  equivalences are based on the
following simple general fact:
\begin{align}
&  \la_\nu(\hw\hu)\setminus \{ \la_\nu(\hw\hu)\cap 
\la_\nu(\hu)\}
=\hu^{-1}(\la_\nu(\hw))\cap \tR_+ \hbox{\ for \ any \ }
\hu,\hw.
\label{latutw} 
\end{align}
For instance, the condition for the lengths readily 
implies
(\ref{ltutw}). For the sake of completeness,  let us 
deduce (\ref{latutw}) from the positivity condition above.
We follow [C2].

It suffices to check that  $\la_\nu(\hw\hu)\supset 
\la_\nu(\hu).$
If there exists a positive  $\tal\in  \la_\nu(\hu)$ 
such that 
$(\hw\hu)(\tal)\in \tR_+,$ then 
$$
\hw(-\hu(\tal))\in \tR_-\ \Rightarrow\ 
-\hu(\tal)\in \la_\nu(\hw)\ \Rightarrow\ 
-\tal\in \hu^{-1}(\la_\nu(\hw)).
$$ 
We come to a contradiction with the positivity.

Applying (\ref{ltutw}) to the reduced decomposition above,
\begin{align}
\la(\hw) = &\{ \tal^1=\al_{i_1},\
\tal^2=s_{i_1}(\al_{i_2}),\ 
\tal^3=s_{i_1}s_{i_2}(\al_{i_3}),\ldots \notag\\ 
&\ldots,\tal^l=\tw^{-1}s_{i_l}(\al_{i_l}) \}.\label{tal}
\end{align}
The cardinality $l$ of the set $\la(\hw)$ equals $l(\hw).$

This set can be introduced for nonreduced decompositions
as well. Let us denote it by $\tla(\hw)$ to differ from 
$\la(\hw).$
It always contains $\la(\hw)$ and, moreover,  
can be represented in the form 
\begin{align}
&\tla(\hw)\ =\ \la(\hw)\, \cup\, \tla^+(\hw)\, 
\cup\,-\tla^+(\hw),\label{tlaw}
\\ &\where \tla^+(\hw)\ =\ 
(\tR_+\cap\ \tla(\hw))\setminus \la(\hw). 
\notag \end{align}
The coincidence  with $\la(\hw)$ is  for
reduced decompositions only.

\vskip 0.2cm
{\bf Reduction modulo $W$.} The following proposition is from [C4].
It generalizes the construction of the elements 
$\pi_{b}$ for $b\in P_+.$ 

\begin{proposition} \label{PIOM}
 Given $ b\in P$, there exists a unique decomposition 
$b= \pi_b  u_b,$
$ u_b \in W$ satisfying one of the following equivalent conditions:

{i) \ \ } $l(\pi_b)+l( u_b)\ =\ l(b)$ and 
$l( u_b)$ is the greatest possible,

{ii)\  }
$ \la(\pi_b)\cap R\ =\ \emptyset$.

The latter condition implies that 
$l(\pi_b)+l(w)\ =\ l(\pi_b w)$
for any $w\in W.$ Besides, the relation $ u_b(b)
\equal b_-\in P_-=-P_+,$
holds, which, in its turn,
determines $ u_b$ uniquely if one of the following equivalent 
conditions is imposed:

{iii) }
$l( u_b)$ is the smallest possible, 

{iv)\ }
if $\al\in \la( u_b)$ then $(\al,b)\neq 0$.
\end{proposition}
\qed

Since $\pi_b = b u_b^{-1} = u_b^{-1} b_-, $ 
the set $\pi_P=\{\pi_b, b\in P\}$
can be described in terms of $P_-$:
\begin{align}
&\pi_P = \{u^{-1} b_- \for b_-\in P_-,\ u\in W\notag\\ 
&\hbox{\ such\ that\ }
\al\in \la( u^{-1}) \Rightarrow (\al,b_-)\neq 0\}.
\label{lapiom} 
\end{align}
Using the longest element $w'_0$ in the centralizer
$W'_0$ of $b_-,$ such $u$ constitute the set
$$\{u,\, \mid \, l(u^{-1}w'_0)\ =\ l(w'_0)+l(u^{-1})\}.$$
Their number is $|W|/|W'_0|.$

For $\tal=[\al,\nu_\al j]\in \tR_+,$ one has:
\begin{align}
\la(b) = \{ \tal,\  &( b, \al^\vee )>j\ge 0 \iif \al\in R_+,
\label{lambi}\\ 
&( b, \al^\vee )\ge j> 0 \iif \al\in R_-\},
\notag \\  
\la(\pi_b) = \{ \tal,\ \al\in R_-,\ 
&( b_-, \al^\vee )>j> 0 
\iif  u_b^{-1}(\al)\in R_+,
\label{lambpi} \\    
&( b_-, \al^\vee )\ge j > 0 \iif   
u_b^{-1}(\al)\in R_- \}, \notag \\ 
\la(\pi_b^{-1}) = \{ \tal,\  -&(b,\al^\vee)>j\ge 0 \}.
\label{lapimin}\\ 
\la(u_b) = \{ \al\in R_+,\  &(b,\al^\vee)> 0 \}.
\label{laumin}  
\end{align}

Let us introduce 
the following {\it affine 
action} of $\hW$ on $z \in \R^n$:
\begin{align}
& (wb)\llb z \rrb \ =\ w(b+z),\ w\in W, b\in P,\notag\\ 
& s_{\tal}\llb z\rrb\ =\ z - ((z,\al)+j)\al^\vee,
\ \tal=[\al,\nu_\al j]\in \tR.
\label{afaction}
\end{align}
For instance, $(b w)\llb 0\rrb=b$ for any $w\in W.$
The relation to the above action is given in terms of 
the {\it affine pairing} $([z,l], z'+d)
\equal (z,z')+l:$
\begin{align}
& (\hw([z,l]),\hw\llb z' \rrb+d) \ =\ 
([z,l], z'+d) \for \hw\in \hW, 
\label{dform}
\end{align}
where we treat $d$ formally.

Introducing the {\it affine Weyl chamber}
\begin{align}
&\CC_a\ =\ \bigcap_{i=0}^n \LL_{\al_i},\ 
\LL_{[\al,\nu_\al j]}=\{z\in \R^n,\ 
(z,\al)+j>0 \},
\notag \end{align}
we come to another interpretation of the $\la$-sets:
\begin{align}
&\la_\nu(\hw)\ =\  \{\tal\in R^a_+, \,   \CC_a  
\not\subset \hw\llb \LL_{\tal}\rrb, \, \nu_{\al}=\nu \}.
\label{lamaff}
\end{align}

For instance, $\Pi$ is the group of all elements of $\hW$
preserving $\CC_a$ with respect to the affine action.
Geometrically, the elements $\pi_b$ are exactly those
sending the negative $-\CC_a$ of $\CC_a$ to 
the negative $-\CC$ of the nonaffine Weyl chamber
$\CC\equal\{z\in \R^n,$$ (z,\al_i)>0$ as $i>0\}.$
More generally, we have the following proposition.

\begin{proposition}\label{PIOMG}
Given two finite sets of positive 
affine roots $\{\tbe=[\be,\nu_\be i]\}$ and  
$\{\tga=[\ga,\nu_\ga j]\},$
let $L_{\be,\ga}$ be the closure of
the union of $\hw\llb-\CC_a\rrb$
over $\hw\in \hW$ such that
$\tbe\not\in \la(\hw)\ni \tga.$ Then
$$L_{\be,\ga}=
\{z\in \R^n,\  (z,\be)+i\le 0,\ (z,\ga)+j\ge 0
\ \hbox{\ for\ all\ } \tbe,\tga\}.$$
The same holds in the nonaffine variant for
$\CC$ in place of $\CC_a.$
\end{proposition}
\qed

We will later need the
following "affine" variant of Proposition \ref{PIOM}.
Given $z\in \R^n,$
there exists a unique element 
$\tilde{w}$ $=u_z a_z$ with $a_z\in Q$ and
$u_z \in W$  satisfying
the relations
\begin{align}
& z_-\equal \tilde{w} \llb z \rrb \in -\bar{\CC_a}, 
\label{toaff}
\\ 
&(z_-,\vth)=-1  \ \Rightarrow\   u_z^{-1}(\vth) \in R_-,
\hbox{\ and\ }\notag\\ 
&(\al_i,z_-)= 0 \ \Rightarrow\   u_z^{-1}(\al_i)\in R_+,
\ i>0,  
\notag \end{align}
where $-\bar{\CC_a}$ is the negative of
the closure $\bar{\CC_a}$ of ${\CC_a}.$

The element $b_{-}= u_b(b)$ is a unique element
from $P_{-}$ which belongs to the orbit $W(b)$. 
So the equality   $c_-=b_- $ means that $b,c$
belong to the same orbit. We will also use 
$b_{+} \equal w_0(b_+),$ a unique element in $W(b)\cap P_{+}.$
In terms of $\pi_b,$
$$u_b\pi_b\ =\ b_-,\ \pi_b u_b\ =\ b_+.$$

Note that $l(\pi_b w)=l(\pi_b)+l(w)$ for all $b\in P,\ 
 w\in W.$
For instance, 
\begin{align}
&l(b_- w)=l(b_-)+l(w),\ l(wb_+)=l(b_+)+l(w),
\label{lupiw}
\\ 
&l(u_b\pi_b w)=l(u_b)+l(\pi_b)+l(w) \for b\in P,\,
 w\in W.\notag
\end{align}
We will use these relations together with the following 
proposition when calculating the
conjugations of the nonsymmetric Macdonald polynomials.

\begin{proposition}\label{LAPOSI}
The set $\la(\hw)$ for $\hw=w b$
consists of $\tal=[\al,\nu_\al j]$ with
positive $\al$ if and only if
\begin{align}
&(b,\al^\vee)\ge -1 \and (b,\al^\vee)=-1 
\Rightarrow \al\in\la(w)
\label{lupiww}
\end{align}
for all $\al \in R_{+},\ b\in P.$
The elements $\hw\in W\cdot P_+$ are of this type.
\end{proposition}
{\it Proof.} 
We use (\ref{tlaw}):
\begin{align}
&\la(wb)\ =\ \la(b) \cup ({}^{-b}\la(w))\,\setminus\,
 (\tla^+ \cup -\tla^+),\notag\\ 
&\where \tla^+\ =\ \tR_-\cap ({}^{-b}\la(w)). 
\label{lawb}
\end{align}
The latter set is 
$$\{[\al,\nu_\al (b,\al^\vee)]\} \hbox{\ such\ that\ }
(b,\al^\vee)<0,\, \al\in\la(w).
$$
Let us calculate $\la(wb)\cap \la(b).$ We need to remove from
$\la(b)$ the roots
$$-[\al,\nu_\al (b,\al^\vee)]
\for \al>0 \hbox{\ such\ that\ }
(b,\al^\vee)<0,\ \al\in \la(w).
$$
The roots in the form  
$[-\al,\nu_\al j]$ belong to  $\la(b)$ exactly 
for $0<j\le -(b,\al^\vee).$ Therefore $\la(wb)\cap \la(b)=$ 
\begin{align}
&=\ \{\ [\al,\nu_\al j],\ 0\le j<(b,\al^\vee)
\hbox{\ as\ } (b,\al^\vee)>0, \notag \\ 
&[-\al,\nu_\al j],\ 0< j<-(b,\al^\vee)
\hbox{\ as\ } (b,\al^\vee)<0 \hbox{\ and\ } \al\in \la(w),
\label{lawbb}
\\  &[-\al,\nu_\al j],\ 0< j\le -(b,\al^\vee)
\hbox{\ as\ } (b,\al^\vee)<0 \hbox{\ and\ } \al\not\in 
\la(w)\ \}, \notag 
\end{align}
for $\al\in R_+.$

Now let us assume that
nonaffine components of the roots from $\la(wb)$ are 
all positive.
Then $(b,\al^\vee)\ge -1$ for all $\al>0$
and $\la(w)$ contains all $\al$ making $-1$ in the
the scalar product with $b$ above.
These conditions are also sufficient.
\qed

%\vfil
\vskip 0.2cm
{\it Partial ordering on $P.$} 
It will be necessary in the theory of nonsymmetric polynomials.
See [O2,M4]. This ordering was also used in [C2] in the 
process of calculating the coefficients of $Y$-operators:
\begin{align}
&b \le c, c\ge b \for b, c\in P \iif c-b \in A_+,
\label{order}
\\ &b \preceq c, c\succeq b \iif b_-< c_- \hbox{\ or\  }
\{b_-=c_- \hbox{\ and\ } b\le c\}.
\label{succ}
\end{align} 
Recall that $b_-=c_- $ means that $b,c$
belong to the same $W$-orbit.
We  write  $<,>,\prec, \succ$ respectively if $b \neq c$.

The following sets  
\begin{align}
&\si(b)\equal \{c\in P, c\succeq b\},\ 
\si_*(b)\equal \{c\in P, c\succ b\}, \notag\\   
&\si_-(b)\equal \si(b_-),\ 
\si_+(b)\equal \si_*(b_+)= \{c\in P, c_->b_-\}. 
\label{cones}
\end{align}
are convex. 
By {\it convex}, we mean that if
$ c, d= c+r\al\in \si$  
for $\al\in R_+, r\in \Z_+$, then
\begin{align}
&\{c,\ c+\al,...,c+(r-1)\al,\ d\}\subset \si.  
\label{convex}
\end{align}

The convexity of the intersections 
$\si(b)\cap W(b), \si_*(b)\cap W(b)$
is by construction. For the sake of completeness, 
let us check the convexity
of the sets $\si_{\pm}(b).$ 

Both sets are $W$-invariant. Indeed,
$c_- > b_-$ if and only if  $b_+>w(c)>b_-$ for all 
$w\in W.$
The set $\si_-(b)$ is the union of  $\si_+$ and the orbit $W(b).$
Here we use that $b_+\and b_-$ are the greatest and the least
elements of  $W(b)$ with respect to ``$>$''. 
This is known (and can
be readily checked by the induction with respect to the 
length - see e.g. [C2]).

If the endpoints of (\ref{convex}) are between $b_+$ and $b_-$ then
it is true for the orbits of all inner points even if  $w\in W$
changes the sign of $\al$ (and the order of the endpoints). 
Also the elements from $\si(b)$
strictly between $c$ and $d$ (i.e.
$c+q\al,\ $ $0<q<r$) belong to $\si_+(b).$ This gives the required. 

The next two propositions are essentially from [C1].

\begin{proposition }\label{BSTAL}
i) Let $c=\hu\llb 0\rrb,$ where
$\hu$ is obtained by striking out any number of 
$\{s_{j}\}$
from a reduced decomposition of $\pi_b (b\in P).$ 
Then $c\succ  b$. Generally speaking, the
converse is not true even if $c\in W(b)$ (I. Macdonald).
In other words, 
the Bruhat order of $\hW$ when restricted to
$\{\pi_b, b\in P\}$  is stronger than $\succ.$

ii) Letting  $b=s_i\llb c\rrb$
 for $0\le i\le n,$ the element $s_i\pi_c$ can be represented
in the form $\pi_b$ for some $b\in P$ if and only if
$(\al_j,c+d)\neq 0.$ More exactly,
the following three conditions are equivalent: 
\begin{align}
& \{c\succ b\}\Leftrightarrow  \{(\al_i^\vee,c+d)>0\}
\Leftrightarrow \{s_i\pi_c=\pi_b,\ l(\pi_b)=l(\pi_c)+1\}.  
\label{aljb}
\end{align}
\qed
\end{proposition }
%\vfil
\begin{proposition}\label{BSTIL}
i) Assuming (\ref{aljb}), let $i>0$. Then $b=s_i(c),\ $
$b_-=c_-,$ and
$u_b=u_c s_i.$ The set $\la(\pi_b)$ is obtained from
$\la(\pi_c)$ by replacing the strict inequality 
$( c_-, \al^\vee )>j> 0$ for $\al=u_c(\al_i)$
(see (\ref{lapimin})) 
with $( c_-, \al^\vee )\ge j> 0.$ 
Here $\al\in R_-$, and $(c_-,\al^\vee)=(c,\al_i^\vee)>0.$

ii) In the case $i=0,\ $ the following holds: 
$b=\vth+s_{\vth}(c)$, the element $c$ is from $\si_+(b)$,
$b_-=c_- - u_c(\vth)\in P_-,$ and $u_b=u_c s_{\vth}$.
For $\al= u_c(-\vth)=\al^\vee,$ the $\la$-inequality  
$( c_-, \al^\vee )\ge j> 0$ is replaced
with the strict inequality $( c_-, \al^\vee )+2> j> 0.$ 
Here $\al\in R_-$, and $( c_-, \al^\vee)=-(c,\vth)\ge 0.$

iii) For any $c\in P, r\in O',\ $ 
$\pi_r\pi_c=\pi_b$ where $b=\pi_r\llb c\rrb.$
Respectively, $u_b=u_c u_r,\ b= \om_r+u_r^{-1}(c),\
b_-=c_- +u_c w_0(\om_{r}).$ In particular,
the latter weight always belongs to $P_-.$
\end{proposition}
{\it Proof.} Let us check i). First, $\pi_b=s_i\pi_c=
s_i c u_c^{-1}= s_i(c) (u_c s_i)^{-1}.$ The uniqueness
of the latter decomposition gives the coincidence
$u_b=u_c s_i.$ Second, $\la(\pi_b)$ is the union of $\la(\pi_c)$
and $(u_c c^{-1})(\al_i)=[\al,(c,\al_i)]$
for $\al=u_c(\al_i)\in R_-$
(see (\ref{ltutw})).
Third, the inequality with $\al^\vee$ in
(\ref{lapimin}) is strict for $c$ because 
$u_c^{-1}(\al)=\al_i   \in R_+$ 
and becomes non-strict for $b$ since 
$u_b^{-1}(\al)=-\al_i\in R_-.$

For ii), it is the other way round. Namely, the extra affine
root from $\la(\pi_b)\setminus \la(\pi_c)$
is $\pi_c^{-1}(\al_0)=[u_c(-\vth),1-(c,\vth)].$ Therefore
$\al=u_c(-\vth)\in R_-$ and the $\la$-inequality 
for $\pi_c$ is 
non-strict.
As for $\pi_b,$ $\ u_b^{-1}(\al)=
(u_c s_{\vth})^{-1}(\al)=\vth\in R_+,$
and the inequality becomes strict. Explicitly,
$(b_-,\al)=(b,-\vth)=(s_{\vth}(c-\vth),-\vth)=2-(c,\vth).$
\qed

\vfil
{\it Arrows in} $P.$
We write $c\rightarrow b$ or $b\leftarrow c$ in the cases i),ii)
from the proposition above and
use the left-right arrow  $c\leftrightarrow b$ for iii)
or when
$b$ and $c$ coincide.
By $c\rightarrow\!\rightarrow b,$  we mean that 
$b$ can be obtained from $c$ by
a chain of (simple) right arrows. Respectively,
 $c\leftrightarrow\!\rightarrow b$ indicates that 
$\leftrightarrow$
can be used in the chain. Actually no more than one 
left-right arrow
is always sufficient and it can be placed right after $c.$
If such arrows are not involved then the
ordering is obviously stronger than the Bruhat order
given by the procedure from i), Proposition \ref{BSTAL}, 
which in its turn is stronger than ``$\succ$''.

If $l(\pi_b)=l(\pi_c)+l(\pi_b\pi_c^{-1})$ then the reduced
decomposition of $\hw=\pi_b\pi_c^{-1}$ readily 
produces a chain
of simple arrows from $c$ to $b.$ 
The number of right arrows is precisely
$l(\hw)$ (see(\ref{aljb})). 
Only transforms of type ii),iii) will change the
$W$-orbits, adding negative short roots to the corresponding
$c_-$ for ii) and  
the weights in the form $w(\om_r)\ (w\in W)$ for iii).

For instance,
let $c,b,b-c\in P_-.$ Then $\pi_c=c,\ \pi_b=b,\ $  
$\hw=b-c$ is
of length $l(b)-l(c),$ and we  need to decompose $b-c.$ 
If $c=0,\ $ then $b-\om_r\in Q,$ the reduced decomposition of $b$ 
begins with $\pi_r\ (r\in O'),$ and 
the first new $c_-$ is $w_0(\om_r).$
When there is no
$\pi_r$ and $b\in Q\cap P_-,$ the chain always
starts with  $-\vth.$
Note that it gives another proof of 
Proposition \ref{LEAST}.

\vskip 0.2cm
Let us now examine the arrows  
$c\Rightarrow b$ from the
viewpoint of the $W$-orbits. The following proposition is
useful for the classification and description of the
perfect representation (the end of the paper).

%\vfil
\begin{proposition}\label{WORB}
i) Given $c_-\in P_-$, any element in the form 
$c_-+u(\vth)\in P_-$
for $u\in W$ such that $u(\vth)\in R_-$
can be represented as $b_-$
for proper $c$ such that $W(c_-)\ni c\rightarrow b.$ Respectively,
 $c_-+u(\om_r)\in P_- $ can be represented as $b_-$ for proper
$c\leftrightarrow b$.

ii) In the case of $A,D,E,$
any element $b_-$ such that 
$b_-\prec c_-$ (both belong to  $P_-$) can be obtained from
$c_-\in P_-$ using consecutive arrows  
$W(c_-)\ni c\rightarrow b\in W(b_-).$ This cannot be true
for all root systems because only short roots 
may be added to
$c_-$ using such a construction.
\end{proposition}

%\vfil
{\it Proof.}
Let $c'=c_-+u(\om_r)\in P_-.$ Given $a=u(\om_r),$
one can assume that $u$ is the greatest
possible, i.e.
$u(\be)\in R_-$ if $\be\in R_+$ and $(\be^\vee,\om_r)=0.$
Explicitly, $u=u_a^{-1}w_0.$
We use the inequalities $(c',\al)\le 0\ (\al\in R_+)$ 
for $\al$ such that $(\al,c_-)=0$ (if such $\al$  exist).
This gives that $(c',\al)=(\om_r,u^{-1}(\al))$
is either $0$ or negative (actually only $-1$
may appear). In the latter case, $u^{-1}(\al)$
contains $\al_r$ with a negative coefficient, and therefore
 $u^{-1}(\al)\in R_-.$ If the scalar product is zero then
$u^{-1}(\al)$ is negative too because of the maximality of $u.$
Thus $w_0u^{-1}$ leaves such $\al$ in $R_+$ and can be represented
as $u_c^{-1}$ for proper $c\in W(c_-)$ (see Proposition \ref{PIOM}
and (\ref{lapiom})). 

If $c'=c_-+u(\vth)\in P_-,$ we use the same reasoning:
$(\vth,u^{-1}(\al))<0$ implies that $ u^{-1}(\al)$ is negative.
So  $u$ can be represented as $u_c w_0$ and
$c'=c_- -u_c(\vth)$ for proper $c.$ Then we  observe
that $( c_-, u_c(\vth))=(c,\vth)\le 0$ because 
$ u(\vth)$ was assumed to be negative. 

\vfil
It is worth mentioning that the direct statement
(that $c'\in P_-$ for $u$ from
Proposition \ref{BSTIL}) also holds and can
be readily checked. One needs to use that
the weights $\{\om_r\}$ are minuscule and 
$(\vth,\al^\vee)\ge 2$
$(\al\in R),$ with the equality  
exactly for $\al=\vth$ ($\al\in R$).

Now let us assume that $c,b\in P_-$ and 
$0\neq x=b-c\in Q_-$
in the simply-laced case $A,D,E.$ We pick a 
a connected component $I_t$ of
$I=\{i\ |\ (x,\al_i)\le 0\}$ in the Dynkin 
diagram $\Ga$ and
set $x_t=\sum_{i\in I_t} l_i\al_i$ for 
$x=\sum_{i=1}^n l_i\al_i.$
Then $(x_t,\al_i)\le 0$ for $i\in I_t,$
since $(\al_i,\al_j)<0$ as $i\in I_t\not\ni j.$ Obviously, $x_t\ge x.$

Let $\vth_t$ be the maximal short root for the root subsystem 
$R_t\subset R$ 
generated by $\{\al_i, i\in I_t\}$ as simple roots,
Then $-\vth_t\ge x_t$ thanks to Proposition \ref{LEAST}. Moreover,
$(\vth_t,\al_j)\ge -1$ for $j\not\in I_t.$ Here we use 
that the coefficient of $\al_i$ in $\vth_c$ is $1$ if 
$\al_i$ is an
end point of $I_t$ and at the same time
an inner point in $\Ga.$ This can be readily
checked using the tables of [B]. 
Note that one may not assume that  
$\al_i$ is inner in $\Ga$
if $I_t$ is of type $A,D.$ 

Here the value  $-1$ is reached precisely for the
neighbors of the endpoints of $I_t$ in $\Ga\setminus I_t$
(otherwise the scalar product is zero). 
Such points $\al_j$
do not belong to $I$ because $I_t$ is a connected 
component of $I.$
Therefore $(c,\al_j)<0$ since
$(x,\al_j)> 0$ and $(b+x,\al_j)\le 0.$ We see that 
$c>c'=c-\vth_t\ge b$
(all three weights are from $P_-$)
and can continue by induction.
\qed

\vskip 0.2cm
\section{Double Hecke algebras}
\setcounter{equation}{0}

We keep the notations from the previous section. 
For the sake of uniformity,
let 
 $$
([b,l],[b',l'])=(b,b'),\
[\al,\nu_\al j]^\vee=[\al^\vee,j], \ X_0=X_{\al_0}=
q_{\vth}X_{\vth}^{-1}.
$$  
By  $m,$ we denote the least natural number 
such that  $(P,P)=(1/m)\Z.$  So
$m=2 \for D_{2k},\ m=1 \for B_{2k}, C_{k},$
otherwise $m=|\Pi|$.

The definition of double affine Hecke algebra will depend 
on the parameters
$q, \{t_\nu, \nu\in \nu_R\}$. 
The basic field will be $\Q_{q,t}\equal$
$\Q[q^{\pm 1/(2m)},t^{\pm 1/2}]$ formed by 
polynomials in terms of $q^{\pm 1/(2m)}$ and  
$\{t_\nu^{\pm 1/2} \}.$
If $q,t$ are regarded to be complex numbers 
(e.g. when $q$ is a root of unity) 
then the fractional powers of $q,t$ have to be somehow
fixed. We set
\begin{align}
&   t_{\tal} = t_{\al}=t_{\nu_\al},\ t_i = t_{\al_i},\ 
q_{\tal}=q^{\nu_\al},\ q_i=q^{\nu_{\al_i}},\notag\\ 
&\where \tal=[\al,\nu_\al j] \in \tR,\ 0\le i\le n.\label{talj}
\end{align}

To simplify formulas we will use the parameters
$\{k_\nu\}$ together with  $\{t_\nu \}$ setting
$$
t_\al=t_\nu=q_\al^{k_\nu} \for \nu=\nu_\al \and
\rho_k=\sum_\nu k_\nu \rho_\nu.
$$
For instance, by $q^{(\rho_k,\al)}$ we mean 
$\prod_{\nu\in\nu_R}t_\nu^{(\rho_\nu,\al^\vee)}.$  
Here $\al\in R$ and therefore it is a product of
integral powers of $t_{\sht}$ and $t_{\lng}.$ 

Let $X_1,\ldots,X_n$ be pairwise
commutative and algebraically independent.
We set   
\begin{align}
& X_{\tb}\ =\ \prod_{i=1}^nX_i^{l_i} q^{ j} 
\iif \tb=[b,j],
\label{Xde}
\\  
&\hbox{where\ } b=\sum_{i=1}^n l_i \om_i\in P,\ j \in 
\frac{1}{ m}\Z.
\notag \end{align}
The elements $\hw \in \hW$ act in 
the ring $\Q[q^{\pm 1/m}][X]$ of polynomials in 
terms of $X_b  (b\in P)$ and $\ q^{\pm 1/m}$
by the formulas:
\begin{align}
&\hw(X_{\tb})\ =\ X_{\hw(\tb)}. 
\end{align}
In particular,
\begin{align}
&\pi_r(X_{b})\ =\ X_{ u^{-1}_r(b)} q^{(\om_{r^*},b)} 
\for \al_{r^*}\equal \pi_r^{-1}(\al_0), \ r\in O'.
\label{pi} \end{align}

Note that the involution $r\mapsto r^*$ of the set $O'$
satisfies the relations $ u_r u_{r^*}=1=\pi_r\pi_{r^*}$.
Moreover, $w_0(\om_r)=-\om_{r^*}$ for the longest 
element $w_0\in W.$
Thus $\al_r\mapsto \al_{r^*}$ is nothing else 
but the automorphism
of the nonaffine Dynkin diagram (preserving $\al_0$). 
This can be 
readily seen from the tables of [B], 
where only the case of $D_n$
requires some consideration.

Recall the notation $(d,[\al,j])=j.$ For instance,
$(\al_0^\vee,b+d)=1-(b,\vth).$

\begin{definition }
The  double  affine Hecke algebra $\HH\ $
(see [C2]) 
is generated over the field $ \Q_{ q,t}$ by 
the elements $\{ T_i,\ 0\le i\le n\}$, 
pairwise commutative $\{X_b, \ b\in P\}$ satisfying 
(\ref{Xde}),
and the group $\Pi,$ where the following relations are imposed:

(o)\ \  $ (T_i-t_i^{1/2})(T_i+t_i^{-1/2})\ =\ 
0,\ 0\ \le\ i\ \le\ n$;

(i)\ \ \ $ T_iT_jT_i...\ =\ T_jT_iT_j...,\ m_{ij}$ 
factors on each side;

(ii)\ \   $ \pi_rT_i\pi_r^{-1}\ =\ T_j \iif 
\pi_r(\al_i)=\al_j$; 

(iii)\  $T_iX_b T_i\ =\ X_b X_{\al_i}^{-1} \iif 
(b+d,\al^\vee_i)=1,\
0 \le i\le  n$;

(iv)\ $T_iX_b\ =\ X_b T_i$ if $(b+d,\al^\vee_i)=0 
\for 0 \le i\le  n$;

(v)\ \ $\pi_rX_b \pi_r^{-1}\ =\ X_{\pi_r(b)}\ =\ 
X_{ u^{-1}_r(b)}
 q^{(\om_{r^*},b)},\  r\in O'$.
\label{double}
\end{definition }
\qed

Given $\tw \in W^a, r\in O,\ $ the product
\begin{align}
&T_{\pi_r\tw}\equal \pi_r\prod_{k=1}^l T_{i_k},\where 
\tw=\prod_{k=1}^l s_{i_k},
l=l(\tw),
\label{Tw}
\end{align}
does not depend on the choice of the reduced decomposition
(because $\{T\}$ satisfy the same ``braid'' relations 
as $\{s\}$ do).
Moreover,
\begin{align}
&T_{\hv}T_{\hw}\ =\ T_{\hv\hw}\  \hbox{ whenever}\ 
 l(\hv\hw)=l(\hv)+l(\hw) \for
\hv,\hw \in \hW. \label{TT}
\end{align}
In particular, we arrive at the pairwise 
commutative elements 
\begin{align}
& Y_{b}\ =\  \prod_{i=1}^nY_i^{k_i} \iif  
b=\sum_{i=1}^nk_i\om_i\in P,\where  
 Y_i\equal T_{\om_i},
\label{Yb}
\end{align}
satisfying the relations
\begin{align}
&T^{-1}_iY_b T^{-1}_i\ =\ Y_b Y_{a_i}^{-1} \iif 
(b,\al^\vee_i)=1,
\notag\\ 
& T_iY_b\ =\ Y_b T_i \iif (b,\al^\vee_i)=0,
 \ 1 \le i\le  n.
\end{align}

\vskip 0.2cm
{\bf Automorphisms.}
The following maps can be uniquely extended to automorphisms of
$\HH\ $(see [C1,C4]):
\begin{align}
\vep:\ &X_i \mapsto Y_i,\ \  Y_i \mapsto X_i,\  
\ T_i \mapsto T_i^{-1},\ 
t_\nu \mapsto t_\nu^{-1},
 q\mapsto  q^{-1},
\label{vep}
\\ 
 \tau_+: \ X_b \mapsto X_b,\ \ &Y_r \mapsto 
X_rY_r q^{-(\om_r,\om_r)/2},\ 
T_i\mapsto T_i,\ t_\nu \mapsto t_\nu,\
 q\mapsto  q,
\notag\\ 
\tau_+:\ &Y_0 \mapsto Y_0 T_0^{2} X_0=
q^{-1}T_{s_\vth}^{-1}T_0X_0,\
\ Y_0=Y_{{\al_0}}\equal q_{\vth}^{-1}Y_{\vth}^{-1},
\label{tau}
\end{align}
where $1\le i\le n,\ r\in O'.$
The formulas for the images of $\{Y_r,Y_0\}$ readily give:
\begin{align}
\vep(T_0)\ =\ &(Y_0 T_0 X_0)^{-1}\ =\ X_\vth T_{s_\vth},\ 
\tau_+(T_0)\ =\ X_0^{-1}T_0^{-1},\notag\\ 
\vep(\pi_r)\ =\ &X_rT_{ u_r^{-1}},\ \tau_+(\pi_r)\ =
\ q^{-(\om_r,\om_r)/2}X_r\pi_r= q^{(\om_r,\om_r)/2}
\pi_rX_{r^*}^{-1},\notag\\ 
&\for \pi_r X_{r^*}\pi_r^{-1}\ =\ q^{(\om_r,\om_r)}
X_{r}^{-1},\ 
X_{r^*}T_{ u_r}X_r\ =\ T_{ u_{r^*}}^{-1}.
\label{vphto}
\end{align}

Theorem 2.3 from [C4] states that the mapping
\begin{align}
& \Bigl(\begin{matrix}  0 &-1\\ -1& 0\end{matrix} \Bigr) 
\mapsto \vep,\ 
  \Bigl(\begin{matrix} 1& 1\\  0& 1 \end{matrix} \Bigr) 
\mapsto \tau_+
\label{glz}
\end{align}
can be extended to a homomorphism of $GL_2(\Z)$ to the
group of automorphisms of \HH \ 
modulo  conjugations 
by the central elements from the group generated
by $T_1,\ldots,T_n$. In this statement,
the quadratic relation (o) from Definition 
\ref{double} may be omitted. Only group relations matter. 

We will also need the following automorphisms
of $\HH\ :$
\begin{align}
&\tau_- \equal  \vep\tau_+\vep,\and  
\si\equal \tau_+\tau_-^{-1}\tau_+ =
\tau_-^{-1}\tau_+\tau_-^{-1}= \vep\si^{-1}\vep. 
\label{taumina}
\end{align}
They preserve
$T_1,\ldots,T_n,\ t,q$ and are uniquely defined from the
following relations:

\begin{align}
& \tau_-: Y_b \mapsto Y_b, \ X_r \mapsto 
Y_r X_r q^\frac{(\om_r,\om_r)}{ 2},\notag\\ 
& \tau_-: X_0 \mapsto (T_0X_0Y_0 T_0)^{-1}=qT_{s_{\vth}}
X_0^{-1}T_0^{-1},
\label{tauom} \\ 
& \si: X_b \mapsto Y_b^{-1},\ Y_r \mapsto 
Y_r^{-1}X_r Y_r q^{-(\om_r,\om_r)},\ 
Y_0 \mapsto T_{s_\vth}^{-1}X_0^{-1}T_{s_\vth},
\notag \end{align}
where $b\in P,\ r\in O'$. Besides,
\begin{align}
&\tau_-(T_0)\ =\ T_0,\ \tau_-(\pi_r)\ =\ \pi_r\ (r\in O'),
\notag\\ 
&\si(T_0)=T_{s_\vth}^{-1}X_\vth^{-1},\ 
\si(\pi_r)=Y_r^{-1}X_r\pi_r q^{-(\om_r,\om_r)}=
T_{u_r}^{-1}X_{r^*}^{-1}. 
\label{taumin}
\end{align}

Thus $\tau_-$ corresponds to  $(\begin{matrix} 1& 0\\  
1& 1 \end{matrix})$,
$\si$ to  $(\begin{matrix} 0& 1\\  -1& 0 \end{matrix})$.

The relation $ \tau_+\tau_-^{-1}\tau_+ =$
$\tau_-^{-1}\tau_+\tau_-^{-1}$ is exactly the definition of 
the projective
$PSL_2(\Z)$ due to Steinberg. It formally gives that
$\si^2$ 
commutes with $\tau_\pm.$ Generally speaking,
$\si^2$ is not inner in $\HH\ $. It is the 
conjugation by $T_{w_0}^{-1}$ if $w_0=-1.$
See below. Always, $\si^4$ is
the conjugation by $T_{w_0}^{-2}$, which is central 
in the nonaffine Hecke algebra $\lan T_i,\, 
1\le i\le n\ran.$

All automorphisms introduced above are unitary. 
An automorphism
$\om$ of \HH\ is called {\it unitary} if 
$\star\,\om\,\star=\om^{-1}$ 
for the main anti-involution $\star$ from [C2]:
\begin{align}
  & X_i^\star\ =\  X_i^{-1},\ \  Y_i^\star\ =\ 
 Y_i^{-1},\  \
 T_i^\star \ =\  T_i^{-1}, \notag\\ 
&t_\nu \mapsto t_\nu^{-1},\
 q\mapsto  q^{-1},\ 0\le i\le n,\ 
(AB)^\star=B^\star A^\star.
\label{star}
\end{align}
The commutativity with $\star$ is obvious because it
is the inversion with respect to
the multiplicative structure of $\HH\ .$

The following automorphism of $\HH^{\flat}\ $
will play an important role in the
paper: $\ \eta\equal \vep\si.$ It is uniquely defined
from the relations
\begin{align}
\eta:\ &T_i\mapsto T_i^{-1},\ X_b\mapsto X_b^{-1},\ 
\pi_r\mapsto \pi_r,\notag\\  
&\where 0\le i\le n,\ b\in P,\ r\in O'.
\label{etatxpi}
\end{align}
It 
conjugates $q,t,$ and extends the Kazhdan-Lusztig involution
on the affine Hecke algebra generated by 
$\{T_i,\, i\ge 0\}.$ 

Here we have used (\ref{tauom}) to calculate the images of
$X, T_i\ (i>0).$ Applying (\ref{taumin}):
\begin{align}
&\vep\si(T_0)=\vep(T_{s_\vth}^{-1}X_0)= 
T_{s_\vth}Y_0=T_0^{-1},\ 
\vep\si(\pi_r)=\notag\\ 
&\vep(T_{u_r}^{-1}X_{r^*}^{-1})= 
T_{u_r^{-1}}Y_{r^*}^{-1}=
T_{u_{r^*}}Y_{r^*}^{-1}= \pi_{r^*}^{-1}=\pi_{r}.
\notag \end{align}

We remind the reader that  $\pi_{r^*}=\pi_r^{-1}$ and
$u_{r^*}=u_r^{-1},$ where $r\in O,$ 
$\al_{r^*}= \pi_r^{-1}(\al_0).$ The map $r\mapsto r^*$ 
corresponds to the automorphism of the nonaffine Dynkin
diagram. It is straightforward to check that
$$\eta(Y_r)=T_{w_0}Y_{r^*}^{-1}T_{w_0}^{-1} \for r\in O'.
$$
We are going to generalize it to arbitrary $Y_b.$

We will use the same symbol $\varsigma$ for two
different maps, namely, for the following one: 
\begin{align}
&b\mapsto b^\varsigma\ =\ -w_0(b),\, b\in P,
\ \ s_i\mapsto s_{\varsigma(i)},\where
\label{varsigma}
\\ 
\varsigma(\al_i)=& 
\al_{\varsigma(i)} = -w_0(\al_i) \for 1\le i\le n,\ 
\varsigma(0)=0,\ \varsigma(\al_0)=\al_0,
\notag \end{align}
which is extended naturally to $\hW,$
and for the corresponding automorphism of $\HH\ :$
\begin{align}
&\varsigma:\ X_b\mapsto X_{\varsigma(b)}, \
Y_b\mapsto Y_{\varsigma(b)},\ 
T_i\mapsto T_{\varsigma(i)},\
\pi_r\mapsto \pi_{r^*}.
\label{varsigmah}
\end{align}
The one from (\ref{varsigmah})
commutes with all previous automorphisms.

\begin{proposition}\label{ETAY}
\begin{align}
&\eta(Y_b)\ =\ T_{w_0}Y_{w_0(b)}T_{w_0}^{-1},\
\si(Y_b)= T_{w_0}^{-1}X_{w_0(b)}T_{w_0}, 
\label{etay}
\\ 
&\si^2(X_b)=T_{w_0}^{-1}X_{w_0(b)}^{-1}T_{w_0},\
 \si^2(Y_b)=T_{w_0}^{-1}Y_{w_0(b)}^{-1}T_{w_0},\notag\\ 
&\si^2(T_0)=T_{w_0}^{-1} T_0 T_{w_0},\
T_{w_0}^{-1} Y_\vth T_{w_0}=T_{s_\vth}^{-1} Y_\vth 
T_{s_\vth},\notag\\ 
&T_{\varsigma(i)} = T_{w_0}^{-1}T_iT_{w_0} \for i>0,\
T_{w_0}^{-1}T_{s_\vth}T_{w_0}=T_{s_\vth},\notag\\ 
&\si^2(H)=T_{w_0}^{-1}\varsigma(H)T_{w_0},\
 \si^4(H)=T_{w_0}^{-2}HT_{w_0}^2\for H\in \HH.
\notag \end{align}
\end{proposition}

{\it Proof.} The formulas for $\si$ easily follow from those
for $\eta,$ so it suffices to check that 
$\eta(Y_b)=T_{-b}^{-1}$ for $b\in P_+$
coincides with  
$$T_{w_0}Y_{w_0(b)}T_{w_0}^{-1}=
T_{w_0}Y_{\varsigma(b)}^{-1}T_{w_0}^{-1}=
T_{w_0}T_{\varsigma(b)}^{-1}T_{w_0}^{-1},
$$
i.e. we need the relation $T_{-b}=
T_{w_0}T_{\varsigma(b)}T_{w_0}^{-1}.$

Since $-b\in P_-\ $,$\ T_{-b}T_{w_0}=T_{(-b)\cdot w_0}$.
Here we apply (\ref{lupiw}). 
Similarly, $\varsigma(b)\in P_+$
and 
$$T_{w_0}T_{\varsigma(b)}\ =\ T_{w_0\cdot(\varsigma(b))}.
$$
However $w_0\cdot(\varsigma(b))=w_0\cdot w_0\cdot 
(-b)\cdot 
w_0=(-b)\cdot w_0.$
\qed

%\comment
%We see that $\eta$ coincides with the conjugation
%$^*$ on the standard
%generators $\{X,T,\pi\}.$ However it is a homomorphism in
%contrast to $^\star,$ which is an anti-homomorphism.
%Note that $(\hY_a)^\star=\eta(\hY_a)$ in the polynomial
%representation, but $\eta(Y_a)$ does not 
%coincide with $Y_a^\star,$ 
%and cannot be expressed in terms of $Y$ only.
%For instance, 
%$$
%\eta(Y_r)=\eta(\pi_rT_{u_r})=
%\pi_{r^*}T_{u_{r^{*}}}^{-1}=
%T_{u_r}Y_{r}^{-1}T_{u_{r^*}}^{-1}
%$$
%is different from $(Y_r)^\star=Y_r^{-1}.$
%\endcomment

\vskip 0.2cm
{\bf Demazure-Lusztig operators.} Following  
[L, KK, C2], 
\begin{align}
&\hT_i\  = \  t_i ^{1/2} s_i\ +\ 
(t_i^{1/2}-t_i^{-1/2})(X_{\al_i}-1)^{-1}(s_i-1),
\ 0\le i\le n,
\label{Demaz}
\end{align}
preserve $\Q[q,t^{\pm 1/2}][X]$.
We note that only $T_0$ involves $q$: 
\begin{align}
&\hT_0\  =  t_0^{1/2}s_0\ +\ (t_0^{1/2}-t_0^{-1/2})
( q X_{\vth}^{-1} -1)^{-1}(s_0-1),\notag\\ 
&\where
s_0(X_i)\ =\ X_iX_{\vth}^{-(\om_i,\th)}
 q^{(\om_i,\vth)},\ 
\al_0=[-\vth,1].
\end{align}

\begin{theorem }
i) The map $ T_j\mapsto \hT_j,\ X_b \mapsto X_b$ 
(see (\ref{Xde})),
$\pi_r\mapsto \pi_r$  (see (\ref{pi})) induces a 
$ \Q_{ q,t}$-linear 
homomorphism from \HH\ to the algebra of linear endomorphisms 
of $\Q_{ q,t}[X]$.
This representation, which will be called polynomial, 
is faithful 
and remains faithful when   $q,t$ take  
any nonzero complex values assuming that
$q$ is not a root of unity. 

ii) For arbitrary nonzero $q,t$ any element $H \in \HH\ $  
have
a unique decomposition in the form
\begin{align}
&H =\sum_{w\in W }  g_{w}  T_{w} f_w,\ 
g_{w} \in \Q_{ q,t}[X],\ f_{w} \in \Q_{ q,t}[Y]  
\label{hatdec}
\end{align}
and five more analogous decompositions corresponding 
to the other orders
of $\{T,X,Y\}.$

iii) The image $\hH$ of $H\in \HH\ $  
is uniquely determined from the following condition:
\begin{align}
&\hH(f(X))\ =\ g(X)\for H\in \HH\ ,\iif Hf(X)-g(X)\ 
\in\notag\\ 
 & \{\sum_{i=0}^n H_i(T_i-t_i^{1/2})+
\sum_{r\in O'} H_r(\pi_r-1), \where H_i,H_r\in \HH\ \}.
\label{hat}
\end{align}
The automorphism $\tau_-$ preserves (\ref{hat}) and
therefore acts in the polynomial representation.
\label{FAITH}
\end{theorem }

{\it Proof}. This theorem is 
essentially from [C2] and [C3]. One only needs to
extend ii) to roots of unity (it was not formulated in
[C2] in the complete generality). In the first place, 
the existence of such a decomposition is true for all 
$q,$ which follows directly from the defining relations. 
Secondly, the uniqueness holds for generic
$q$ because the polynomial representation is faithful 
(a similar argument was used in [C2]). Finally, the number 
of linearly independent expressions in the form 
(\ref{hatdec}) 
(when the degrees of $f,g$ are  bounded) may not become 
smaller for special $q$ than for generic $q$ 
(here by special 
we mean the roots of unity).   
\qed

We will also use the above statements for the 
{\it intermediate subalgebras} of \HH\
with $P$ replaced by its
subgroups containing $Q$. 
Let $B$ be any lattice between $Q$ and $P,$ respectively, 
$\Pi^{\flat}$
the preimage of $B/Q$ in $\Pi,$ and 
$\hW^{\flat}=\Pi^{\flat}\cdot \tW=B\cdot W.$
By $\HH^{\flat}\ ,$ we denote the subalgebra of $\HH\ $ generated
by $X_b (b\in B), \pi\in \Pi^{\flat}$ and $T_i 
(0\le i\le n).$
Thus $Y_b\in \HH^{\flat}$ when $b\in B.$ Actually the lattices
for $X_b$ and $Y_b$ can be different, but when discussing
the Fourier transforms it is convenient to impose the coincidence.
Here $m$ can be replaced by $\tilde{m}$ such that 
$\tilde{m}(B,B)\subset \Z$ in the definition of 
$\Q_{q,t}.$

\begin{proposition}\label{interfa}
The algebra  $\HH^{\flat}$ satisfies all claims of 
Theorem \ref{FAITH} for $B$ instead of $P$ and the polynomial
representations in $\Q_{q,t}[X_b, b\in B].$
The automorphisms $\vep,\tau_{\pm},\si,\eta$
and the anti-involution $\star$ preserve this subalgebra. 
\end{proposition}

{\it Proof}. The  compatibility of the definition of 
 $\HH^{\flat}$
with the automorphisms $\tau_{\pm}$ follows directly
from the formulas for their action on $X_r,Y_r (r\in O')$ and
$T_0$ (see (\ref{vphto}, \ref{tauom})). Since  
$\HH^{\flat}$ is
generated by $X_b,Y_b (b\in\tB)$ and 
$\{T_1,...,T_n\}$ this holds
for $\vep$ as well. Claim i) from Theorem \ref{FAITH} remains
true for the polynomials in $X_b (b\in B)$
because the formulas for $\hT_i$ involve $X_\al$ only.

To check ii) we need a more complete version of the 
relations (iii,iv) from Definition \ref{double} (cf. 
[L] and formula (\ref{Demaz}) above).
Namely, for all $b\in P,$
\begin{align}
&T_iX_b -X_{s_i(b)}T_i\ =\ 
(t_i^{1/2}-t_i^{-1/2})\frac{s_i(X_b)-X_b}{
X_{\al_i}-1},\ 0 \le i\le  n.
\label{tixi}
\end{align}
This relation and its dual counterpart for 
$Y$ instead of $X$
ensure the existence of the decompositions (\ref{hatdec}).
Using the fact that the polynomial representation is faithful 
(for generic $q$) we establish the uniqueness (for all 
$q\neq 0$).
Claim iii) formally results from ii).
\qed

The following Proposition is essentially from [C4]
(Proposition 3.3). In the notation from 
(\ref{cones}), we set
\begin{align}
&\Si(b)\equal \oplus_{c\in \si(b)}\Q_{q,t} x_{c} \ ,
\Si_*(b)\equal \oplus_{c\in \si_*(b)}\Q_{q,t} x_{c}
\label{sibx}
\\ 
&\Si_-(b)\ =\ \Si(b_-),\ \Si_+(b)\ =\ \Si_*(b_+), 
\ b\in B.
\notag \end{align}
 
Let us denote the maps $b\mapsto -b$ by $\imath$.
 Recall that
$b\mapsto -w_0(b)$ is denoted by  $\varsigma.$ Respectively,
$$
\imath: X_b\mapsto X_{-b}=X_b^{-1},\
\varsigma: X_b\mapsto X_b^\varsigma=\varsigma(X_b)=
X_{-w_0(b)}.
$$
Then  
$\varsigma$ 
preserves the
ordering and sends 
$\Si(b)$ and $\Si_*(b)$ to $\Si(\varsigma(b))$ and 
$\Si_*(\varsigma(b)).$

\begin{proposition}
i)Given $b\in B,\ \tal=[\al,\nu_\al j]\in \tR$ with 
$\al>0,$
the operator 
\begin{align}
&\r_{\tal}\  = \  t_\al ^{1/2} +\ 
(t_\al^{1/2}-t_\al^{-1/2})(X_{\tal}^{-1}-1)^{-1}
(1-s_{\tal})
\label{rtal}
\end{align}
preserves $\Si(b),\Si_*(b),\Si_{\pm}(b).$ Moreover,
\begin{align}
\r_{\tal}(X_b)\,\hbox{mod}\, \Si_+(b)\  &= 
\  t_\al ^{1/2}X_b +\ 
(t_\al^{1/2}-t_\al^{-1/2})s_{\tal}(X_b) \iif (b,\al)<0,
\notag\\ 
&= \  t_\al^{-1/2}X_b \iif (b,\al)>0,\notag\\ 
&= \  t_\al^{1/2}X_b \iif (b,\al)=0.
\label{rtalx}
\end{align}

ii)The operators $\hw^{-1}T_{\hw}$ preserve 
$\Si(b),\Si_*(b),$ and
$\Si_{\pm}(b)$ 
if all 
$\tal=[\al,\nu_\al j]\in\la(\hw)$
have positive $\al.$ The elements $\hw\in W\cdot B_+$
have this property. For instance, $Y_c$ leave
the $\Si$-sets invariant for all $c\in B$
and 
$$\imath\cdot T_{w_0}:
\Si(b)\to \Si(\varsigma(b)),\ 
\Si_*(b)\to \Si_*(\varsigma(b)).
$$

iii) Assuming ii),
\begin{align}
&\hw^{-1}T_{\hw}(X_b)\,\hbox{mod}\, \Si_*(b)\  
 = \  \prod _{(\tal,b)\le 0} t_\al ^{1/2}\, 
\prod_{(\tal,b)> 0} t_\al ^{-1/2}\, X_b\notag\\ 
&\hbox{multiplied\ over\ } 
\tal=[\al,\nu_\al j]\in \la(\hw).
\label{rtalt}
\end{align}
In particular,
\begin{align}
Y_c(X_b)\,\hbox{mod}\, \Si_*(b)\ &=\ 
q^{(c,u_b^{-1}(\rho_k)-b)}\, X_b,
\label{rtalp}
\\  
\imath\cdot T_{w_0}(X_b)\,\hbox{mod}\, \Si_*(b) &=
\prod_{\nu\in \nu_R} t_\nu ^{l_\nu(w_0)/2 -l_\nu(u_b)}\, 
X_{\varsigma(b)}\,\hbox{mod}\, \Si_*(\varsigma(b)). 
\notag \end{align}
\label{RTAL}
\end{proposition}
\qed

The operators
$\{T_j, 0\le j\le n\}$ may not fix the spaces 
$\Si(b), \Si_*(b).$ However they  
leave $\Si_-(b)$ and $\Si_+(b)$
invariant. Hence 
the {\it affine Hecke algebra} $\h_Y^{\flat}$
generated by $T_i (0\le i\le n)$ and the group 
$ \Pi^{\flat}$
act in the space
\begin{align} 
V(b)\equal\Si_-(b)/\Si_+(b)\ \cong\ \Q_{q,t}[X_c,
 c\in W(b)].
\label{vbmin}
\end{align} 

This representation is the (parabolically) induced 
$\h_Y^{\flat}$-module 
generated by the image $X_+$ of $X_{b_+}$ subject to   
\begin{align}
&T_i(X_+)\ =\ t_i^{1/2}X_+\iif (b_+,\al_i)=0,\ 
1\le i\le n,\notag\\ 
&Y_a(X_+)\ =\ q^{(a,u_{b_+}^{-1}(\rho_k)-b_+)}X_+,
\ a\in B.
\label{xplus}
\end{align}
It is true for arbitrary nonzero $q,t.$

Note the following explicit formulas: 
\begin{align}
s_i(X_c)\ &=\ (t_i^{1/2}T_i(X_c))\iif (c,\al_i)>0,
\notag\\ 
&=\ (t_i^{1/2}T_i(X_c))^{-1} \iif (c,\al_i)<0,
\notag\\ 
T_i(X_c)\ &= t_i^{1/2} X_c\iif (c,\al_i)=0,\notag\\ 
\pi_r(X_c)&=X_{\pi_r(c)},\ c\in W(b),\ r\in O'.
\label{sonx}
\end{align}
Here $0\le i\le n,\ (c,\al_0)=-(c,\vth)$.

\vskip 0.2cm
\section {Macdonald polynomials}
\setcounter{equation}{0}

Continuing the previous section, we set 
\begin{align}
&\mu\ = \mu^{(k)}=\prod_{\al \in R_+}
\prod_{i=0}^\infty \frac{(1-X_\al q_\al^{i}) 
(1-X_\al^{-1}q_\al^{i+1})
}{
(1-X_\al t_\al q_\al^{i}) 
(1-X_\al^{-1}t_\al^{}q_\al^{i+1})}.\
\label{mu}
\end{align}
It is considered as a Laurent series with the
coefficients in  $\Q[t_\nu][[q_\nu]]$ for
$\nu\in \nu_R.$ 

We denote  the constant term
of $f(X)$ 
by $\langle  f \rangle.$ 
Let $\mu_\circ\equal \mu/\langle \mu \rangle$, 
where 
\begin{align}
&\langle\mu\rangle\ =\ \prod_{\al \in R_+}
\prod_{i=1}^{\infty} \frac{ (1- q^{(\rho_k,\al)+i})^2
}{
(1-t_\al q^{(\rho_k,\al)+i}) 
(1-t_\al^{-1}q^{(\rho_k,\al)+i})
}.
\label{consterm}
\end{align}
Recall that $q^{(\rho_k,\al)}=q_\al^{(\rho_k,\al^\vee)},
\ t_\al=q_\al^{k_\al}.$
This formula is from [C2]. It is nothing else but the
Macdonald constant term conjecture from [M1].

The coefficients of the Laurent series $\mu_\circ$ are actually 
from the field of rationals 
$\Q(q,t)\equal\Q(q_\nu,t_\nu),$ where $\nu\in \nu_R.$ One 
also has $\mu_\circ^*\ =\ \mu_\circ$ for the involution  
$$
 X_b^*\ =\  X_{-b},\ t^*\ =\ t^{-1},\ q^*\ =\ q^{-1}.
$$
This involution is the restriction
of the anti-involution $\star$
from (\ref{star}) to $X$-polynomilas (and Laurent series).  
These two properties of $\mu_\circ$
can be directly seen from the difference relations
for $\mu.$ We will prove the first in the next proposition.

Setting 
\begin{align}
&\langle f,g\rangle_\circ\ \equal\ \langle 
\mu_\circ f\ {g}^*\rangle\ =\ 
\langle g,f\rangle_\circ^* \for
f,g \in \Q(q,t)[X], 
\label{innerpro}  
\end{align}
we  introduce the {\it nonsymmetric Macdonald 
polynomials} $\ E_b(x)\ =\ E_b^{(k)}\in$
$\Q(q,t)[X]$ for $b \in P$ by means of 
the conditions
\begin{align}
&E_b-X_b\ \in\ \oplus_{c\succ b}\Q(q,t) X_c,\
\langle E_b, X_{c}\rangle_\circ = 0 \for P \ni c\succ b.
\label{macd}
\end{align}
They are well defined
because the  pairing   
is nondegenerate and form a 
basis in $\Q(q,t)[X]$. 

This definition is due to Macdonald (for
$k_{\sht}=k_{\lng}\in \Z_+ $),
who extended
Opdam's nonsymmetric polynomials introduced
in the differential case in [O2]
(Opdam mentions Heckman's contribution in [O2]). 
The general case was considered in [C4]. 

Another approach is based on the $Y$-operators.
See formulas (\ref{xplus}) and (\ref{rtalp}) above,
and also Proposition \ref{RTAL}.

\begin{proposition}
The polynomial representation is $\star$-unitary:
\begin{align}
&\langle H(f),g\rangle_\circ\  = \ 
\langle f, H^\star(g)\rangle_\circ \for H\in \HH, 
\ f\in \Q_{q,t}[X].
\label{staru}
\end{align}
The polynomials $\{E_b, b\in P\}$
are unique (up to proportionality) eigenfunctions of
the operators $\{L_f\equal f(Y_1,\cdots, Y_n), 
f\in \Q[X]\}$
acting in $\Q_{q,t}[X]:$
\begin{align}
&L_{f}(E_b)\ =\ f(q^{-b_\#})E_b, \where
b_\#\equal b- u_b^{-1}(\rho_k),
\label{Yone}
\\ 
& X_a(q^{b})\ =\
q^{(a,b)}\for a,b\in P,\ 
 u_b=\pi_b^{-1}b \hbox{\ is\ from\ Section\ 1.}
\label{xaonb}
\end{align}\label{YONE}
\end{proposition}
\qed

In the previous section we denoted $Y,T$ acting in the 
polynomial representation by $\hY,\hT.$ 
Here and further we will mainly omit the hat 
(if it does not lead to confusion). 
Thus we have two equivalent definitions of the
nonsymmetric polynomials. Both are compatible
with the transfer to the intermediate subalgebras
$\HH^{\flat}$ and subspaces 
$\Q_{q,t}[X_b]\equal\Q_{q,t}[X_b,
 b\in B]\subset \Q_{q,t}[X].$ 

Let us check that the coefficients of 
$\mu_\circ$ and the Macdonald
polynomials 
are rational functions in terms of $q_\nu,t_\nu.$
We need to make the construction of 
$\lan\ ,\ \ran_\circ$ more abstract.
A form $( f\, ,\, g )$ will be called 
 {\it $\ast$-bilinear} if 
\begin{align}
& ( r f\, ,\, g )\ =\  r\,(f\, ,\, g )\ =\ 
 (f\, ,\, r^*\,g )\for
r\in \Q_{q,t}.
\label{bilin}
\end{align}

\begin{proposition}
i) Forms $\lan f\, ,\, g \ran$ on $\Q_{q,t}[X]$ satisfying 
(\ref{staru}) and which are $\ast$-bilinear 
are in one-to-one correspondence with  
$\Q_{q,t}$-linear maps
$\varpi:\ $ $\Q_{q,t}[X]\to$ $ \Q_{q,t}$ such that
\begin{align}
&\varpi(T_i(u)) = t_i\, \varpi(u),\ 
\varpi(Y_a(u)) = q^{(a\, ,\, \rho_k)}\varpi(u),\ 
i\ge 0,\ a\in P.
\label{tyva}
\end{align}
Such a form is
$\ast$-hermitian in the sense of (\ref{innerpro}) 
if and only
if $\varpi(f^*)=$ $ (\varpi(f))^*.$
Given  $\varpi,$ the corresponding form is 
$\lan f\, ,\, g \ran=$  
$\varpi(fg^*).$

ii) Replacing $\Q_{q,t}$ by the field of rationals
$\widetilde{\Q}(q,t)=\Q(q^{1/(2m)},t^{1/2})$ 
in the definition  of $\HH\, $ and 
in the polynomial representation,
there is a unique nonzero $\varpi$
up to proportionality. Namely, it is the restriction to 
$$\widetilde{\Q}(q,t)[X]\subset 
\HH\otimes_{\Q} \widetilde{\Q}(q,t)
$$
of the linear map
$\varpi_{\hbox{\small ext}}:\ \HH\,\to 
\widetilde{\Q}(q,t)$ uniquely 
determined from the
relation $H -\varpi_{\hbox{\small ext}}(H)\i\subset $ 
$\JJ\equal$
\begin{align}
&\{\sum_{r\in O'} \bigl(H_r(\pi_r-1)+ 
(\pi_r-1)H_r'\bigr)+
  \sum_{i=0}^n 
\bigl(H_i(T_i-t_i^{1/2})+(T_i-t_i^{1/2})H_i'\bigr)\},
\notag
\end{align}
where $ H_i,H_i',H_r,H_r'\in \HH\,$.
Equivalently, 
\begin{align}
&f(X)\  =\ \varpi(f)+ 
\sum_{r\in O'} (\pi_r-1)(h_r)+ 
\sum_{i=0}^n (\hT_i-t_i^{1/2})(h_i)
\label{hatf}
\end{align}
for $ h_i,\, h_r\, \in\, \widetilde{\Q}(q,t)[X].$  
\label{UNIABS}
\end{proposition}
{\it Proof.} In the first place,
given a form $\lan f,g\ran$ satisfying
(\ref{staru}), the corresponding linear form is of course
$\varpi(f)=$ $\lan f , 1 \ran.$ It obviously satisfies
(\ref{hatf}). Its extension to  $\HH\,$
is  $\varpi_{\hbox{\small ext}}(H)=$ 
$\varpi(\hH(1)).$ 
 
Now let $\varpi_{\hbox{\small unv}}$ be
the projection  $\HH\, \to \HH\,/\JJ.$
We set
$$
\lan A\, ,\, B\ran_{\hbox{\small unv}}
 =\ \varpi_{\hbox{\small unv}} (AB^\star)\ =
\ \lan B\, ,\, A\ran_{\hbox{\small unv}}^* ,
$$ 
where the  anti-involution acts naturally on   
$\HH\, /\JJ$
thanks to the $\star$-invariance of $\JJ.$

By construction, 
$$
\lan HA\, ,\, B\ran_{\hbox{\small unv}}\ =\ 
\varpi_{\hbox{\small unv}} (HAB^\star)\ =\ 
\varpi_{\hbox{\small unv}} (AB^\star H)\ =\ 
\lan A\, ,\, H^\star B\ran_{\hbox{\small unv}}
$$ 
for $H=T_i,\pi_r.$ Upon the restriction to $\Q_{q,t}[X],$
\begin{align}
&\lan \hH(f)\, ,\, g\ran\ =\
 \lan f\, ,\,H^\star(g)\ran \for
\lan f\, ,\, g\ran\equal\varpi_{\hbox{\small unv}} 
(f(X)g(X)^*)
\label{vrpf}
\end{align}
and such $H.$ Obviously (\ref{vrpf}) holds when $H$ is 
the multiplication by a polynomial. So it is true
for all $H\in \HH\, .$
Thus the form $\lan f\, ,\, g\ran$ satisfies the same 
relations as $\lan f\, ,\, g\ran_\circ.$ Namely, it is
hermitian with respect to $^*$ 
and $\HH\,$-invariant with respect to $\star.$
Recall that its  values are in the vector space 
$\HH\,/\JJ$
with a natural action of $\star.$ 

The form $\lan f\, ,\, g\ran$ is universal among such forms.  
To be more exact,
an arbitrary  $\Q_{q,t}$-valued linear form on 
$\Q_{q,t}[X]$
obeying the same
relations as $\varpi_{\hbox{\small unv}} $ 
is the composition of $\varpi_{\hbox{\small unv}}$
and a homomorphism  $\om:\HH\,/\JJ\to\Q_{q,t}.$
The latter has to satisfy the $\star$-invariance 
relations $\om(H^\star)=\om(H)^*$
for $H\in \HH\,$ to make the corresponding bi-form
$\ast$-hermitian.
 
Let us switch from $\Q_{q,t}$ to the field of rationals  
$\Q(q,t).$
We already know that at least one
form $\lan\ ,\ \ran$ exists for generic $q,t$ 
with the values in a proper completion of $\Q(q,t).$
It is given by (\ref{innerpro}) and is
unique up to proportionality. Indeed,
the polynomial representation is irreducible, 
the $Y$-operators
are diagonalizable there, and the $Y$-spectrum is simple. 
Thus the space  $\HH\,/\JJ$ is one-dimensional upon
proper completion, so it has to be one-dimensional
over $\Q(q,t)$ as well.
We get the rationality of the coefficients of
$\mu_\circ$ and the coefficients of the
nonsymmetric Macdonald polynomials.
\qed

\vskip 0.2cm
Following Proposition \ref{YONE}, the
{\it symmetric 
Macdonald polynomials} $\ P_b=P_b^{(k)}\ $ can be introduced as
eigenfunctions of the $W$-invariant
operators $L_f=$ $f(Y_1,\cdots,Y_n)$ 
defined for symmetric, i.e.  $W$-invariant, polynomials $f$ as follows:
\begin{align}
&L_{f}(P_b)\ =\ f(q^{-b+\rho_k}) P_b,\ \
b\in P_-,\notag\\ 
&P_b\ =\ \sum_{c\in W(b)}X_c \mod \Si_+(b).
\label{Lf}
\end{align}
Here it suffices to take the 
{\it monomial symmetric
functions}, namely,
$f_b=\sum_{c\in W(b)}X_c$ for $b\in P_-.$

The $P$-polynomials
are pairwise orthogonal with respect to 
$\langle\ ,\ \rangle_\circ$
as well as $\{E\}.$
Since they are $W$-invariant,  $\mu$ can be replaced
by Macdonald's truncated theta function:
\begin{align}
&\de\ = \de^{(k)}=\prod_{\al \in R_+}
\prod_{i=0}^\infty \frac{(1-X_\al q_\al^{i}) 
(1-X_\al^{-1}q_\al^{i})
}{
(1-X_\al t_\al q_\al^{i}) 
(1-X_\al^{-1}t_\al^{}q_\al^{i})}.\
\label{delp}
\end{align}
The corresponding pairing remains $\ast$-hermitian because
$\de_\circ$ is $\ast$-invariant. 

These polynomials were introduced in [M2,M3]. 
Actually they
appeared for the first
time in Kadell's work (for some root systems). 
In one-dimensional case, they are due to Rogers (see [AI]).  

The connection between $E$ and $P$ is as follows: 
\begin{align}
&P_{b_-}\ =\ \P_{b_+} E_{b_+}, \ b_-\in P_-,\ 
b_+=w_0(b_-),\notag\\  
&\P_{b_+}\equal\sum_{c\in W(b_+)}
\prod_\nu t_\nu^{l_{\nu}(w_c)/2} \hT_{w_c},
\label{symmetr}
\end{align}
where $w_c\in W$ is the element of the least length such that
$c=w_c(b_+).$ See [O2,M4,C4]. 

\vskip 0.2cm
{\it Spherical polynomials}.
Mainly we will use the following renormalization 
of the $E$-polynomials (see [C4]):  
\begin{align}
\e_b\ \equal&\ E_b(X)(E_b(q^{-\rho_k}))^{-1},\where 
\notag\\ 
E_{b}(q^{-\rho_k}) \ =&\ q^{(\rho_k,b_-)} 
\prod_{[\al,j]\in \la'(\pi_b)}
\Bigl(
\frac{ 
1- q_\al^{j}t_\al X_\al(q^{\rho_k})
 }{
1- q_\al^{j}X_\al(q^{\rho_k})
}
\Bigr),
\label{ebebs}\\ 
\la'(\pi_b)\ =&\ 
\{[\al,j]\ |\  [-\al,\nu_\al j]\in \la(\pi_b)\}.
\label{jbseto}
\end{align}
We call them {\it spherical polynomials}. 
Explicitly
(see (\ref{lambpi})),
\begin{align}
\la'(\pi_b)\ =& \{[\al,j]\, \mid\,\al\in R_+,
\label{jbset}  
\\  
&-( b_-, \al^\vee )>j> 0 
\iif  u_b^{-1}(\al)\in R_-,\notag\\ 
&-( b_-, \al^\vee )\ge j > 0 \iif   
u_b^{-1}(\al)\in R_+ \}. 
\notag \end{align}
Formula (\ref{ebebs}) is the Macdonald evaluation conjecture
in the nonsymmetric variant.

Note that one has to consider only long $\al$ (resp., short)
if $k_{\sht}=0,$ i.e. $t_{\sht}=1$ (resp., 
$k_{\lng}=0$) in this set.
All formulas below involving $\la$ or $\la'$ have to be modified
correspondingly in such case.  

We have the following duality relations:
\begin{align}
&\e_b(q^{c_{\#}})\ =\ \e_c(q^{b_{\#}}) \for b,c\in P,\ 
b_\# = b- u_b^{-1}(\rho_k),
\label{ebdual}  
\end{align}
justifying the definition above.

\vskip 0.2cm
{\it Conjugations.}
We now come to the formulas for
the conjugations of the nonsymmetric polynomials, which
are important in the theory of the Fourier transform.
The automorphism $^*$ is well defined on $E_b$ because
their coefficients are rational functions in terms of
$q,t.$

\begin{proposition}
The conjugation $\ast$ in the polynomial representation
is induced by the involution $\eta$ of $\HH^{\flat}.$
The paring $\lan fg\mu_0\ran$ induces on $\HH\ $ 
the anti-involution
$\star\cdot \eta=\eta\cdot \star.$ 
For $b\in B,$
\begin{align}
&E_b^*\ =\ 
\prod_{\nu\in \nu_R} t_\nu ^{l_\nu(u_b)-l_\nu(w_0)/2}\, 
T_{w_0}(E_{\varsigma(b)}),\where \varsigma(b)=-w_0(b).
\label{ebast}
\end{align}
\label{PCONJ}
\end{proposition}
{\it Proof.}
The conjugation $\ast$ sends $\hT_i$ to $\hT_i^{-1}$
for $0\le i\le n$ i.e. coincides with the action 
of $\star$ on $T_i.$ This can be checked
by a simple direct calculation. However $\ast$ is an
involution in contrast to $\star$ which is anti-involution.
By definition,
the conjugation takes $X_b$ to $X_b^{-1}$ and fixes 
$\pi_r$
for $r\in O.$ Using (\ref{etatxpi}), we conclude that
$\ast$ is induced by the involution $\eta$ of 
$\HH^{\flat}.$

Applying (\ref{rtalp}) to $\varsigma(b),$  
\begin{align}
&T_{w_0}(E_{\varsigma(b)})\,\hbox{mod}\, \Si_*(b)=
\prod_{\nu\in \nu_R} 
t_\nu ^{l_\nu(w_0)/2 -l_\nu(u_{\varsigma(b)})}\, 
X_{-b}\,\hbox{mod}\, \imath(\Si_*(b)). 
\notag \end{align}
Here $\imath$ sends $X_b\mapsto X_{-b}$ and fixes $q,t.$
Obviously $\imath(\Si_*(b))$ coincides with  
$(\Si_*(b))^\ast$ and $l_\nu(u_{\varsigma(b)})=
l_\nu(u_b).$
This gives the coefficient of proportionality in 
(\ref{ebast}).

Conjugating,
$$\eta(Y_c)(E_b^*)=q^{(c,b_\#)}E_b^*=
q^{(w_0(c),w_0(b_\#))}E_b^*=q^{-(w_0(c),w_0(-b)_\#)}
E_b^*.
$$
The last transformation is possible because 
the automorphism $\varsigma=-w_0$
is compatible with the representations $b=\pi_b u_b.$
Finally,
$\eta(Y_c)=$ $T_{w_0}Y_{w_0(c)}T_{w_0}^{-1}$ due to
Proposition \ref{ETAY} (formula (\ref{etay})), and
$T_{w_0}^{-1}(E_b^*)$ has to be proportional to
$E_{w_0(-b)}.$
\qed

For the symmetric polynomials, 
$\ P_b^*=P_{\varsigma(b)},\ b\in P_-,\ $ which readily
results from the orthogonality property and 
the second relation from 
(\ref{Lf}).

The conjugation formula
(\ref{ebast}) 
gets simpler in the spherical normalization. Indeed,
$$q^{(\rho_k,b_-)}=\prod_{\nu\in \nu_R} 
t_\nu^{-l_\nu(b)/2},$$
and (\ref{ebebs}) reads
\begin{align}
&E_b(q^{-\rho_k})\ =\ \prod_\nu t_\nu^{-l_\nu(u_b)/2} 
M_b,
\label{ebubl} 
 \\ 
&M_b\ =\ \prod_{[\al,j]\in \la'(\pi_b)}
\Bigl(
\frac{ 
t_\al^{-1/2}- q_\al^{j}t_\al^{1/2} X_\al(q^{\rho_k})
 }{
1- q_\al^{j}X_\al(q^{\rho_k})
}
\Bigr),
\notag \end{align}
where the factor $M_b$ is "real", i.e. $\ast$-invariant.
Hence,
\begin{align}
\e_b^*\ &=\ M_b^{-1}(E_b\prod_\nu 
t_\nu^{l_\nu(u_b)/2})^*\notag\\  
&=\prod_{\nu} t_\nu ^{l_\nu(u_b)/2-l_\nu(w_0)/2}\, 
T_{w_0}(E_{\varsigma(b)})\notag\\ 
&=\prod_{\nu} t_\nu ^{-l_\nu(w_0)/2}\, 
T_{w_0}(\e_{\varsigma(b)}).
\label{ebeconj}
\end{align}
We will give later a better proof of this relation 
without
using (\ref{ebebs}). It will be based on 
the intertwining operators and the definition of the Fourier
transform. 

\vskip 0.2cm
{\bf Intertwining operators.}
The {\it $X$-intertwiners} (see e.g. [C1])
are introduced as follows:
\begin{align}
&\Phi_i\ =\ 
T_i + (t_i^{1/2}-t_i^{-1/2})(X_{\al_i}-1)^{-1},
\notag\\ 
& S_i=(\phi_i)^{-1}\Phi_i,\ 
G_i=\Phi_i (\phi_i)^{-1},\
\notag\\  
&\phi_i\ =\  t_i^{1/2} + 
(t_i^{1/2} -t_i^{-1/2})(X_{\al_i}-1)^{-1},
\label{Phi}
\end{align}
for $0\le i\le n$. They belong to \HH\ extended by the field 
$\Q_{q,t}(X)$
of rational functions in $\{X\}$. The elements
 $S_i$ and $G_i$ 
satisfy the same relations
as $\{s_i,\pi_r\}$ do.  Hence
the map 
\begin{align} 
\hw\mapsto S_{\hw}\ =\ \pi_r S_{i_l}\cdots S_{i_1},
\where \hw=\pi_r s_{i_l}\cdots s_{i_1}\in \hW,
\label{Phiprod}
\end{align} 
is a  well defined homomorphism  from $\hW.$
The same holds  for $G_{\hw}.$

As to $\Phi_i$, they
 satisfy the relations 
for $\{T_i\}$ (i.e. the  homogeneous Coxeter relations 
and those with $\pi_r$). So the decomposition in 
(\ref{Phiprod}) has
to be reduced.

It suffices to check these properties
in the polynomial representation, where they are 
obvious since
$S_i=s_i$ for $0\le i\le n.$ A direct proof is based on 
the following 
property of $\{\Phi\}:$ 
\begin{align}
&\Phi_{\hw} X_b\ =\ X_{\hw(b)}\Phi_{\hw},\ \hw\in \hW.
\label{Phix}
\end{align}
Here the elements $\Phi_{\hw}=
\Phi_{\pi_r}\Phi_{s_{i_l}}\cdots
\Phi_{s_{i_1}}$ can be introduced
for any choice of the reduced decomposition. Relation
(\ref{Phix}) fixes them uniquely up to left
(or right) multiplications by functions of $X.$
Since $\Phi_{\hw}-T_{\hw}$ is a combination
of $T_{\hw'}$ for $l(\hw')<l(\hw),$ one gets that  
$\Phi_{\hw}$
does not depend on the particular choice of the reduced 
decomposition of $\hw$,
and $\Phi$ has the desired multiplicative property.
  
We will also use that $\ \Phi_i,\ \phi_i\ $ are 
self-adjoint with 
respect to 
the anti-\-involution (\ref{star}) and therefore
\begin{align}
& \Phi_{\hw}^\star = \Phi_{\hw^{-1}},\ 
S_{\hw}^\star\ =\ G_{\hw^{-1}},  \ \hw\in \hW.
\label{Phistar}
\end{align}
This gives that $G$ and $S$ are $\star$-unitary up to 
a functional
coefficient of proportionality. Explicitly:
\begin{align}
 & (G_{\hw})^\star\,G_{\hw}\ =
\ \bigl((S_{\hw})^\star\, S_{\hw}\bigr)^{-1}\ =\ 
\label{gstarg} 
\\ 
&\prod_{\al,j} \Bigl(
\frac{
t_\al^{1/2}-q_\al^j t_\al^{-1/2} X_\al}{
t_\al^{-1/2}-q_\al^j t_\al^{1/2} X_\al 
}
\Bigr),\ [\al,\nu_\al j]\in \la(\hw).
\notag \end{align}

To define the {\it $Y$-intertwiners} we apply
the involution $\vep$ to $\Phi_{\hw}$ and to $S,G$.
The $Y$-intertwiners satisfy the same
$\star$-relations (\ref{Phistar}).
The formulas can be easily calculated using
(\ref{vphto}). In the case of $GL_n,$ one gets 
the intertwiners from [KnS].
For $i>0,$ the  $X_{\al_i}$ in $\Phi_{s_i}$ are replaced 
by  $Y_{\al_i}^{-1}$ without touching $T$ and $t.$
We use that $\vep$ sends nonaffine $T_i$ to $T_i^{-1},$
$t^{1/2}$ to $t^{-1/2},$ and transpose $X$ and $Y.$

When applying the $Y$-intertwiners to the spherical
polynomials (or any eigenfunctions of the $Y$-operators) 
it is more convenient to use $\tau_+$ instead of $\vep.$ 
The following proposition is from [C1].
We set
\begin{align}
&\Phi_i^b\ =\ \Phi_i(q^{b_{\#}})\ =\
T_i + (t_i^{1/2}-t_i^{-1/2})(X_{\al_i}
(q^{b_{\#}})-1)^{-1},
\label{Phijb}
\\ 
&G_i^b\ =\ (\Phi_i\phi_i^{-1})(q^{b_{\#}})\ =\
\frac{ T_i + (t_i^{1/2}-t_i^{-1/2})(X_{\al_i}
(q^{b_{\#}})-1)^{-1}
}{
t_i^{1/2} + (t_i^{1/2} -t_i^{-1/2})(X_{\al_i}
(q^{b_{\#}})-1)^{-1} }.
\label{Gjb}
\end{align}

\begin{proposition}
Given $c\in P,\ 0\le i\le n$ such that 
$(\al_i, c+d)> 0,$ 
\begin{align}
&E_{b}q^{-(b,b)/2} \ =\ t_i^{1/2}\tau_+(\Phi_i^c) (E_c)
 q^{-(c,c)/2}  
\for b= s_i\llb c\rrb
\label{Phieb}
\end{align}
and the automorphism $\tau_+$ from (\ref{tau}). 
If  $(\al_i, c+d)=0$
then
\begin{align}
&\tau_+(T_i) (e_c) \ =\ t_i^{1/2} e_c, \ 0\le i\le n, 
\label{Tjeco}
\end{align}
which  results in the relations $s_i(E_c)=E_c$ as
$i>0$.  For $b=\pi_r\llb c\rrb,$ where
the indices $r$ are from $O',$
\begin{align}
&q^{-(b,b)/2+(c,c)/2}E_b\ =\ \tau_+(\pi_r)(E_c)\ =\
X_{\om_r}q^{-(\om_r,\om_r)/2}\pi_r(E_c).
\label{pireb}
\end{align}
\label{PHIEB}
\end{proposition}
\qed

We can reformulate the proposition using the spherical 
polynomials.
We will also replace $P$ by the lattice 
$Q\subset B\subset P.$
All above considerations are compatible with the 
reduction
$P\to B.$
Let 
\begin{align}
&\hE_{b}\ =\ \tau_+(\pi_r G_{i_l}^{c_l}\ldots 
G_{i_1}^{c_1})(1),
\where  
\label{ehatb}
 \\ 
&c_1=0, c_2=s_{i_1}\llb c_1\rrb,\ldots,
c_l=s_{i_l}\llb c_{l-1}\rrb, 
\for \pi_b= \pi_r s_{i_l}\ldots s_{i_1}.
\notag \end{align}
These polynomials 
do not depend on the particular choice of
the decomposition of $\pi_b$ (not necessarily
reduced), and are proportional to $E_b$ for all
$b\in B$:
\begin{align}
E_{b}q^{-(b,b)/2} \ =&\ \prod_{1\le p\le l}
\bigl(t_{i_p}^{1/2}\phi_{i_p}(q^{c_p})\bigr)\ 
\hE_b\notag\\ 
=&\ \prod_{[\al,j]\in \la'(\pi_b)}
\Bigl(
\frac{ 
1- q_\al^{j}t_\al X_\al(q^{\rho_k})
 }{
1- q_\al^{j} X_\al(q^{\rho_k})
}
\Bigr)\ \hE_b.
\label{ebebhat}
\end{align}

Combining (\ref{ebebhat}) and (\ref{ebebs}), 
we conclude that
\begin{align}
&\e_{b}\ =\ q^{(-\rho_k+b_-,b_-)/2} \hE_b.
\label{ebebeq}
\end{align}
See [C1] for a proof which is not based on the explicit formulas.

As an application we get that a spherical polynomial 
$\e_b$
is well defined for $q,t\in \C^*$ if 
\begin{align}
&\prod_{[\al,j]\in \la'(\pi_b)}
\bigl(
1- q_\al^{j}t_\al x_a(t^\rho)\bigr)\ \neq\ 0.
\label{ebebnze}
\end{align}

\vskip 0.2cm
\section {Fourier transform on polynomials}
\setcounter{equation}{0}
 
We will begin with the norm-formula for the spherical
polynomials. See [C1] Theorem 5.6, and [M3,C2,M4,C4].
Further $B$ will be any lattice between $Q$ and $P,$ 
$\{b_i, 1\le i\le n\}$ a basis of $B,$ 
$\HH^{\flat}\ $ the 
corresponding intermediate subalgebra of $\HH\ .$
We also set $\hW^{\flat}=B\cdot W\subset \hW,\ $
$\Q_{q,t}[X_b]=\Q_{q,t}[X_b, b\in B],$
and replace $m$ by the least $\tilde{m}\in \N$ such that 
$\tilde{m}(B,B)\subset \Z$ in the definition of the 
$\Q_{q,t}.$

\begin{proposition} 
For $b,c\in B$ and the Kronecker delta $\de_{bc}$,
\begin{align}
 &\lan \e_b,\e_c\ran_\circ\ =\ \lan 
\hE_b,\hE_c\ran_\circ\ =\ 
\de_{bc}\mu^{-1}(q^{b_{\#}})\mu(q^{-\rho_k})=\notag\\  
&\de_{bc}\prod_{[\al,j]\in \la'(\pi_b)}
\Bigl(
\frac{
t_\al^{1/2}-q_\al^jt_\al^{-1/2} X_\al(q^{\rho_k})}{
t_\al^{-1/2}-q_\al^jt_\al^{1/2} X_\al(q^{\rho_k}) 
}
\Bigr). 
\label{normhats}
\end{align}
\label{NORMSHAT}
\end{proposition}
{\it Proof.} Using (\ref{ehatb}) and that $\tau_+$ is
$\star$-unitary, 
$$
\lan \hE_b,\hE_b\ran_\circ = 
\lan \tau_+(G_{\pi_b})(1)\, ,\, 
\tau_+(G_{\pi_b})(1) \ran_\circ =
\lan G_{\pi_b}^\star G_{\pi_b}(1)\, ,\, 1\ran_\circ.
$$
Hence we can apply (\ref{gstarg}) substituting 
$X_\al\mapsto
X_\al(q^{-\rho_k}).$
\qed

%\vfil
\vskip 0.3cm
{\it Discretization.}
This proposition will be interpreted as calculation of the
Fourier transform from the polynomial representation
to the functional representation of $\HH^\flat .$
The latter depends on $n$ independent parameters 
denoted by the formal exponentials
$q^{\xi_i}$ for $\{1\le i\le n\}.$ 

We set
\begin{align}
&X_a(q^{\xi})\equal \prod_{i=1}^n q^{l_i\xi_i},\where 
a=\sum_{i=1}^n l_i b_i\in B,\notag\\ 
&X_a(bw)\ =\ X_a(q^{b+ w(\xi)})\ =\
q^{(a,b)}X_{w^{-1}(a)}(q^{\xi}),\notag\\  
&\hu (g)(bw)\ =\ g(\hu^{-1}bw),\ b\in B, w\in W.
\label{deltag}
\end{align}
Note that $b+w(\xi)=(bw)\llb \xi\rrb.$

These formulas naturally determines
the {\it discretization homomorphism}
$g(X)\mapsto \de(g)(\hw)$ (depending on $\xi$) 
which maps 
the space of rational functions
$g\in\Q_{q,t}(X)$ or more general functions of $X$ to
$$\FF[\xi]\equal \hbox{Funct}(\hW^{\flat},\Q_\xi)
$$ 
formed by
$\Q_\xi$-valued functions on $\hW^{\flat}.$ 

Here and further
$\Q_\xi\equal\Q_{q,t}(q^{\xi_1},\ldots,q^{\xi_n})$.
We  omit $\de$ when it is clear that $\FF[\xi]$
is considered.

This homomorphism can be naturally extended to
the operator algebra 
$\AA\equal \oplus_{\hu\in \hW^{\flat}}\Q_{q,t}(X)\hu.$
For instance, the discretizations $\delta(\hH)$ of 
the operators
$\hH$ for $H\in \HH^{\flat}\ $ are well-defined
and we get the action 
of $\HH^{\flat}\ $ in $\FF_\xi.$
We will call the latter 
the {\it functional representation}. 
It contains the image of
the polynomial representation 
as an $\HH^{\flat}\ $-submodule. 

Explicitly, $\de(X_a)$ is  
(the multiplication by) the discretization of $X_a,$
$\de(\pi_r)=\pi_r,$ for  $\pi_r\in \Pi^{\flat}$, and 
\begin{align}
{\de}(T_i(g))&(\hw))\ =\ 
\frac{t_i^{1/2}X_{\al_i}(q^{w(\xi)})
q^{(\al_i,b)} - t_i^{-1/2}}{
{X_{\al_i}(q^{w(\xi)})q^{(\al_i,b)}- 1}  }\
g(s_i \hw)\notag\\ 
-&\frac{t_i^{1/2}-t_i^{-1/2}}{
{X_{\al_i}(q^{w(\xi)})q^{(\al_i,b)}- 1}  }\
g(\hw) \for 0\le i\le n,\ \hw=bw.
\label{Tfunct}
\end{align}
The map $\de$ is a $\HH^{\flat}\ $-homomorphism.

\vskip 0.2cm
The inner product corresponding to the anti-involution $\star$ of
$\HH^{\flat}\ $ is given by the discretization of the function 
$\mu_\bullet(X)=\mu(X) /\mu(q^\xi):$
\begin{align}
&\mu_\bullet(\hw)\ = \
\prod_{\al,j}
\Bigl(
\frac{
t_\al^{-1/2}-q_\al^jt_\al^{1/2} X_\al(q^{\xi})}{
t_\al^{1/2}-q_\al^jt_\al^{-1/2} X_\al(q^{\xi})
}
\Bigr),\ [\al,\nu_\al j]\in \la(\hw).
\label{muval}
\end{align}
Actually we will use only  
$\de(\mu_\bullet),$  so by $\mu_\bullet$, we mean
the discretization of $\mu_\bullet.$

We put formally $\xi^*=\xi.$ To be more exact,
$(q^{\xi_i})^*\equal q^{-\xi_i}$ for $1\le i\le n.$
Then $\mu_\bullet(\hw)^*=\mu_\bullet(\hw).$ 
The counterpart of $\lan\,\ran$ is the sum
$$\lan f\ran_\xi\ \equal \ 
\sum_{\hw\in \hW^{\flat}}f(\hw).$$

The corresponding pairing is as follows:
\begin{align}
&\langle f,g\rangle_\bullet\ = \ \lan 
fg^*\mu_\bullet\ran_\xi\ =
\sum_{\hw\in \hW^{\flat}}
\mu_\bullet(\hw) f(\hw)\ g(\hw)^*\ =\ 
\langle g,f\rangle_\bullet^*.
\label{innerdel}  
\end{align}
Here  $f,g$ are 
from the $\HH^{\flat}\ $-submodule of finitely supported
functions  $\F[\xi]\subset \FF[\xi].$

Given $H\in \HH^{\flat}\ ,$
\begin{align}
&
\langle H(f),g\rangle_\bullet\  =\  
\langle f, H^\star(g)\rangle_\bullet, \ \
\langle H(f)\,g\mu_\bullet \rangle_\xi  =\  
\langle f\,\eta(H)^\star(g)\rangle_\xi.
\label{conjxi}  
\end{align}
To justify (\ref{conjxi}), one can follow the 
case of $\lan\, \ran$ or simply use the homomorphism 
$\de.$ 
See Proposition \ref{YONE}. Note that the operators $T_i,X_b$
are unitary for the former pairing
and self-adjoint with respect to the latter.

The {\it characteristic functions} $\chi_{\hw}$ and 
{\it delta-functions}
$\de_{\hw}$ in $ \F[\xi] (\hw\in \hW^{\flat})$ are 
defined from the 
relations 
$$
\chi_{\hw}(\hu)\ =\ \de_{\hw,\hu},\ 
\de_{\hw}(\hu)\ =\ \mu_\bullet(\hw)^{-1}\chi_{\hw} 
$$ 
for the
Kronecker delta. We have usual formulas:
$$ X_b(\chi_{\hw})\equal \de(X_b)(\chi_{\hw})\ =\
X_b(\hw)\chi_{\hw},\ b\in B, \hw\in \hW^{\flat}.    
$$
 
When considering concrete special (non-generic) 
$\xi\in \C^n,$
one need to check that all above formulas are well defined
(the denominators in (\ref{Tfunct}), (\ref{muval}) are nonzero). 

Note that $\mu_\bullet(\hw)=1$ if
$q^{b+w(\xi)}=q^\xi,$ i.e. $q^{(b+w(\xi),a)}=q^{(\xi,a)}$ for all
$a\in B,$ provided that $\mu_\bullet$ is well defined. However it 
is not true,
generally speaking, if (\ref{muval}) is used as a definition
of $\mu_\bullet(\hw),$ without any reference to $\mu_\bullet(X).$

The main example discussed in the paper is the following 
specialization: $\xi=-\rho_k.$ More precisely,
$q^{(\xi,b)}\mapsto q^{-(\rho_k,b)}$ for $b\in B.$ The 
$\HH^{\flat}\ $-module  $\F(-\rho_k)$ becomes reducible. 
Namely,
\begin{align} 
\I_\#\equal \oplus_{\hw\not\in \pi_B}\Q_{q,t}\chi_{\hw}\subset
\F(-\rho_k) \and \F_\#\equal \F(-\rho_k)/\I_\#
\label{forho}
\end{align} 
are $\HH^{\flat}\ $-modules. We denote 
$\lan\ ,\ \ran_\bullet$ by  $\lan\ ,\ \ran_\#$
in this case. Respectively, $\lan\ \ran_\#=$
$\lan\ \ran_{-\rho_k}.$
Since $\I_\#$ is the radical of this pairing it is
a $\HH^{\flat}\ $-submodule.

Here $\F_\#$ can be identified with the space
of finitely supported $\Q_{q,t}$-valued functions
on $\pi_B= \{\pi_b, b\in B\}.$ It is irreducible for generic $q,t.$
The $\HH^{\flat}\ $-module
$\FF_\#= \hbox{Funct}(\pi_B,\Q_{q,t})$
is its natural completion. 

Note that $X_a(\pi_b)=X_a(q^{b_{\#}}).$
Therefore 
$\pi_B$ can be naturally identified 
with the set 
$B_{\#}=\{b_{\#}= b-u_b^{-1}(\rho_k)\}$
for generic $k.$

The $\HH^{\flat}\ $-invariance of $\I_\#$
(see [C1] and [C4]) can be seen directly 
from (\ref{aljb}) and
the following general (any $\xi$) relations:
\begin{align}
T_i(\chi_{\hw})\ =\ 
\frac{t_i^{1/2}X_{\al_i}^{-1}(q^{w(\xi)})
q^{-(\al_i,b)} - t^{-1/2}}{
{X_{\al_i}^{-1}(q^{w(\xi)})q^{-(\al_i,b)} - 1}  }\
&\chi_{s_i\hw}\notag\\ 
-\frac{t_i^{1/2}-t_i^{-1/2}}{
{X_{\al_i}(q^{w(\xi)})q^{(\al_i,b)} - 1}  }\
&\chi_{\hw} \for 0\le i\le n,\notag\\ 
\pi_r(\chi_{\hw})\ =\ \chi_{\pi_r\hw} \for \pi_r\in 
\Pi^{\flat},\ 
\hw=bw\in &\hW^{\flat}.
\label{Tfchar}
\end{align}
 
For further references, the delta-counterparts of these 
formulas are:
\begin{align}
T_i(\de_{\hw})\ =\ 
\frac{t_i^{1/2}X_{\al_i}(q^{w(\xi)})q^{(\al_i,b)} - 
t^{-1/2}}{
{X_{\al_i}(q^{w(\xi)})q^{(\al_i,b)} - 1}  }\
&\de_{s_i\hw}\notag\\ 
-\frac{t_i^{1/2}-t_i^{-1/2}}{
{X_{\al_i}(q^{w(\xi)})q^{(\al_i,b)} - 1}  }\
&\de_{\hw} \for 0\le i\le n.
\label{Tfdel}
\end{align}
The action of $\pi_r$ remains the same as for $\{\chi\}$
because $\mu_\bullet(\pi_r\hw)=\mu_\bullet(\hw).$ 

\vskip 0.3cm
{\bf Basic transforms.}
We introduce the {\it Fourier transform} $\vph_\circ$ and the   
{\it skew Fourier transform} $\psi_\circ$ 
as follows:
\begin{align}
&\vph_\circ(f)(\pi_b)\ =\ \lan f(X)\e_b(X)\mu_\circ\ran,
\ b\in B,\notag\\  
&\psi_\circ(f)(\pi_b)\ =\ \lan f(X)^*\e_b(X)\mu_\circ\ran=
\lan \e_b,f \ran_\circ.
\label{fourier}
\end{align}
They are isomorphisms from the space of polynomials 
$\Q_{q,t}[X_b]\ni f$ onto 
$\F_\#.$ The $\vph_\circ$ is $\Q_{q,t}$-linear, the 
$\psi_\circ$ conjugates
$q,t$ (sends $q\mapsto q^{-1},\ t\mapsto t^{-1}$).
The other two transforms will be defined later in this
section. 

We will denote
the characteristic and
delta-functions $\chi_{\hw},\de_{\hw}\in \F_\# $ for 
$\hw=\pi_b$ by $\chi^\pi_b,\de^\pi_b$ respectively. 

\begin{theorem}
The Fourier transform $\vph_\circ$ is unitary, i.e 
sends $\lan \ ,\ \ran_\circ$ to $\lan \ ,\ \ran_\#,$
and induces the  automorphism $\si$
on $\HH^{\flat}\ :$ $\vph_\circ(H(p))\ =\ \de(\si(H))
(\vph_\circ(p))$
for $p(X)\in \Q_{q,t}[X_b]$, $H\in \HH^{\flat}\ .$ 
Respectively, 
$\psi_\circ$
is unitary too and induces the involution $\vep.$ 
Explicitly,
\begin{align}
&\psi_\circ: \sum g_b\e_b(X)\mapsto 
\sum g_b^{*}\de^{\pi}_b\, \in\, \F_\#,\
\for g_b\in \Q_{q,t},
\label{isofour}
\end{align}
which is equivalent to the relations
\begin{align}
&\e_b(X)\ =\ \vep(G_{\pi_b})(1)\for b\in B,
\label{isvep}
\end{align}
and results in the formula
\begin{align}
&\vph_\circ: \sum g_b\e_b(X)\mapsto 
\sum g_b \prod_{\nu} t_\nu ^{l_\nu(w_0)/2}\, 
T_{w_0}^{-1}(\de^{\pi}_{\varsigma(b)})\notag\\ 
&\for g_b\in \Q_{q,t},\ \varsigma(b)=-w_0(b).
\label{isovph}
\end{align}
\label{ISOVD}
\end{theorem}
{\it Proof.} Let us begin with 
$\psi$ (the subscript ``$\circ$'' will
be dropped until  ``noncompact'' Fourier transforms appear). 
We use the definition of
$\e_b,$ the duality relations (\ref{ebdual}), and of 
course the
fact that the polynomial representation is 
$\star$-unitary
(see (\ref{star}),(\ref{staru})). 
The duality relations are needed in the form of the
{\it Pieri rules} from [C4], Theorem 5.4. Namely, we use that
given $a\in B,$
$X_a^{-1}\e_b(X)$ is a transformation of $Y_a(\de_b^\pi)$ when 
the delta-functions $\de_c^\pi$ are replaced by $\e_c$
(see (\ref{Tfdel})).
This readily gives that $\psi$
induces $\vep$ on the double Hecke algebra. Since $\vep$ is
unitary and both the polynomial representation and    
$\F_\#$ are irreducible (for generic $q,t$) we get that $\psi$
is unitary up to a constant. The constant is $1$
because the  image of $\e_1=1$ is
$\de^\pi_0=\de_{\hbox{\small id}}=\chi^\pi_0.$ 

Thus formula (\ref{isofour}) is a reformulation of 
(\ref{normhats}). Another way to check  it is based on the 
technique of intertwiners. See [C1], Corollary 5.2,
Proposition \ref{PHIEB} above, and formula (\ref{ebebeq}).

Turning to $\vph,$ we will use the involution
$\eta=\vep\si$ sending $\ T_i\mapsto T_i^{-1}$
$(0\le i\le n),$ $X_b\mapsto X_b^{-1},\ $ and
$\pi_r\mapsto \pi_r.$ It is an automorphism of 
$\HH^{\flat}\ ,$
conjugating $q,t.$

Recall that $\eta$ extends the conjugation
$\ast$ in the polynomial representation.
Therefore the transformation of
$\HH^{\flat}\ $ corresponding to $\vph$
is the composition $\vep \eta=\si.$ 
This gives that $\vph$ is unitary.

Concerning (\ref{isvep}), we apply $\psi$ to the relations
\begin{align}
&\de^\pi_b\ =\ G_{\pi_b}\,(\de^\pi_0)\,\for b\in B.
\label{isvde}
\end{align}

To get (\ref{isovph}), we can combine (\ref{isofour}) 
with 
(\ref{ebeconj}) which states that
\begin{align}
&\e_b^*\ =\ \prod_{\nu} t_\nu ^{-l_\nu(w_0)/2}\, 
T_{w_0}(\e_{\varsigma(b)}).
\label{ebeconjj}
\end{align}

Let us give another proof of this formula which
doesn't require the exact value of the 
coefficient of proportionality in 
$\e_b^*\ =\ C_b T_{w_0}(\e_{\varsigma(b)})$ and automatically
gives its value.

What we  have proved can be represented as follows:
\begin{align}
&\lan \e_b\,\e_c\mu_\circ\ran\ =\
C_b^{-1}\lan (T_{w_0}(\e_b))^*\e_c\mu_\circ\ran\notag\\ 
&=\ C_b^{-1}T_{w_0}^{-1}(\de^\pi_{\varsigma(b)})(\pi_c).
\label{ebecmu}
\end{align}
Here the left-hand side is $b\leftrightarrow c\,$- symmetric.

On the other hand,
\begin{align}
&T_{w_0}^{-1}(\de^\pi_{b})(\pi_c)\ =\ 
\lan T_{w_0}^{-1}(\de^{\pi}_b)
\de^{\pi}_c\mu_\bullet\ran_{\#},
\label{tdepibc}
\end{align}
where the last expression is $b\leftrightarrow c\,$- symmetric,
thanks to formula (\ref{conjxi}).

Therefore the multiplier $C_b$ doesn't depend on $b$ and must
be exactly as in (\ref{ebeconjj}). 
\qed

\begin{corollary}
Given $b,c\in P,$
\begin{align}
&\lan \e_b\,\e_c\mu_\circ\ran\ =\
\prod_{\nu} t_\nu ^{l_\nu(w_0)/2}\,
T_{w_0}^{-1}(\de^\pi_{\varsigma(b)})(\pi_c).
\label{ebecm}
\end{align}
The pairing $\lan\!\lan f\, ,\, g\ran\!\ran\equal
\prod_{\nu}t_\nu ^{-l_\nu(w_0)/2}
\lan f\, ,\, T_{w_0}(g^\varsigma)\mu_\circ\ran$ is symmetric on
$\Q_{q,t}[X].$  Here $g^\varsigma(X)=g(w_0(X)^{-1}).$
The corresponding anti-involution is 
the composition 
$\,\eta\cdot\star\cdot T_{w_0}\cdot\varsigma=$
$\si^2 \cdot \eta\cdot\star,$
where by $T_{w_0}$ we mean the corresponding conjugation.
It sends
\begin{align}
& T_i\mapsto T_i\ (1\le i\le n),\ \, Y_b\mapsto Y_b,
\ X_b\mapsto T_{w_0}^{-1}X_{\varsigma(b)}T_{w_0}\ 
(b\in P).
\label{ebecinv}
\end{align}
On has:
$\lan\!\lan \e_b\, ,\, \e_c\ran\!\ran$
$=\de_{bc}\mu^{-1}(q^{b_\#})\mu(q^{-\rho_k}).$
Assuming that $0<q<1$ and imposing the
inequalities $1+(\rho_k,\al^\vee)+k_\al>0$
for all $\al\in R_+,$ this paring is positive and therefore
the polynomial representation is irreducible. For instance,
if $k_{\sht}=k=k_{\lng},$ than the latter inequalities 
mean that
$k>-1/h$ for the Coxeter number $h=(\rho,\vth)+1.$
\label{EBECM}
\end{corollary}
\qed

\vskip 0.2cm
We will introduce the other two Fourier transforms
by conjugating (\ref{fourier}):
\begin{align}
&\bar{\vph}_\circ(f)(\pi_b)\ \equal\ \
\lan f(X)\e_b(X)^*\mu_\circ\ran\ =
\ \lan f\, ,\,\e_b\ran_\circ,\notag\\  
&\bar{\psi}_\circ(f)(\pi_b)\ \equal\ \ 
\lan f(X)^*\e_b(X)^*\mu_\circ\ran,\ b\in B.
\label{fourico}
\end{align}
Let us "conjugate" Theorem \ref{ISOVD}.

\begin{theorem}
The Fourier transforms $\bar{\vph},\,\bar{\psi_\circ}$ 
are unitary, i.e 
transform $\lan \ ,\ \ran_\circ$ to $\lan \ ,\ \ran_\#.$
They induce the automorphisms $\si^{-1}$ and $\eta\si$
on $\HH^{\flat}\ $
respectively and send
\begin{align}
&\bar{\vph}_\circ: \sum g_b\e_b(X)\mapsto 
\sum g_b\de^{\pi}_b \for b\in B, g_b\in \Q_{q,t},\notag\\  
&\bar{\psi}_\circ: \sum g_b\e_b(X)\mapsto 
\sum g_b^* \prod_{\nu} t_\nu ^{-l_\nu(w_0)/2}
T_{w_0}(\de^{\pi}_b). 
\label{isofourc}
\end{align}
\label{ISOVDC}
\end{theorem}
\qed

\vskip 0.2cm
{\it Gauss integrals.}
We introduce the {\it Gaussians} $\ga^{\pm 1}$ as   
 $W$-invariant solutions of the following 
system of difference equations:
\begin{align}
&\om_j(\ga)\ =\  
 q^{(\om_i,\om_i)/2} X_i^{-1}\ga, \ 
\om_i(\ga^{-1})\  =\  q^{-(\om_i,\om_i)/2} X_i\ga^{-1} 
\label{gausseq}
\end{align}
for $1\le i\le n.$
In the current setup,
we need the Laurent series
\begin{align}
&\tga^{-1}\equal\sum_{b\in B} q^{(b,b)/2}X_b.
\label{gauser}
\end{align}
Analytically, the values of $\tga$ on any $q^z\ (z\in \C)$ 
are well defined if
$|q|<1$. However we will mainly operate in formal series
in this section.

The multiplication by $\tga^{-1}$ preserves 
the space of Laurent series
with coefficients from $\Q[[q^{1/(2\tilde{m})}]]$. Recall
that $\tilde{m}\in \N$ is the smallest positive integer
such that
$\tilde{m}(B,B)\subset \Z.$ 
For instance, we can multiply $\tga^{-1}$ by 
the $q$-expansions of 
$\e_b(X)$ 
and $\mu_\circ,$
which is a Laurent series with coefficients in 
$\Q[t][[q]].$
The $q$-expansions of the
coefficients of $\e_b$ (which are rational functions in
terms of $q,t$) belong to $\Q[t^{\pm 1}][[q]].$

\begin{theorem}
Given $b,c\in B$ and the corresponding spherical
polynomials $\e_b, \e_c$,
\begin{align}
\langle \e_b \e_c \tga^{-1}\mu_\circ
\rangle\ & =\ 
q^{(b_\#,b_\#)/2+(c_\#,c_\#)/2 -(\rho_k,\rho_k)} 
\e_c(q^{b_\#})\langle \tga^{-1}\mu_\circ\rangle,
\label{epep}
\\ 
\langle
\e_b \e_c^* \tga^{-1}\mu_\circ
\rangle\ & =\ 
q^{(b_\#,b_\#)/2+(c_\#,c_\#)/2 -(\rho_k,\rho_k)} 
\e_c^*(q^{b_\#})\langle \tga^{-1}\mu_\circ\rangle,
\label{epeps}
\\ 
\langle \e_b^* \e_c^* \tga^{-1}\mu_\circ
\rangle\ & =\ 
q^{(b_\#,b_\#)/2+(c_\#,c_\#)/2 -(\rho_k,\rho_k)}\notag\\  
&\times\prod_{\nu}t_\nu ^{-l_\nu(w_0)/2}
T_{w_0}(\e_c^*)(q^{b_\#})\langle \tga^{-1}\mu_\circ\rangle,
\label{epepss}
\end{align}
where the coefficients of
$\mu_\circ,\, \e_b,\, \e_c,$ and $\e_c^*$ are expanded 
in terms of 
positive powers of $q.$
Here the coefficient of proportionality can be 
calculated explicitly:
\begin{align}
&\langle \tga^{-1}\mu_\circ\rangle\ =\ 
\prod_{\al\in R_+}\prod_{ j=1}^{\infty}\Bigl(\frac{ 
1-t_\al^{-1} q_\al^{(\rho_k,\al^\vee)+j}}{
1-      q_\al^{(\rho_k,\al^\vee)+j} }\Bigr).
\label{mehtamu}
\end{align}
\label{EPEP}
\end{theorem}
{\it The proof} is based on the following fact:
$\hH\tga^{-1}=\tga^{-1}\tau_+(\hH)$ in the polynomial
representation extended by the Gaussian, where 
$H\in \HH\ $ 
and the mapping $H\mapsto \hH$ is 
from Theorem \ref{FAITH}, Section 2. Indeed,
the conjugation  by $\tga$
corresponds to $\tau_+$ on the standard generators 
$X_b(b\in B),$ $T_i(0\le i\le n),$
$\pi_r(r\in O').$ To be more exact,
$\hH\tga^{-1}$ coincides with $\tga^{-1}(\tau_+(H))\hat{}$
for all $H\in \HH^\flat\ .$
It can be readily deduced from the $W$-invariance of 
$\tga$
and (\ref{gausseq}).  

The same holds in the functional representations 
$\FF[\xi]$
if we take
\begin{align}
&\ga^{\pm 1}(bw)\equal\ q^{\pm(b+w(\xi),b+w(\xi))/2},
 \ b\in B,\ 
w\in W.  
\label{gaussf}
\end{align}
Here we need to extend the field of constants $\Q_{\xi}$
by $q^{\pm(\xi,\xi)/2}.$ We can easily avoid this by  
removing $\pm(\xi,\xi)/2$ from the exponent
of (\ref{gaussf}), since we need the
Gaussian only up to proportionality. However we prefer  
to stick to the standard 
$\ga(q^z)=q^{(z,z)/2}$ in this section. 

Involving Theorem \ref{ISOVD} we conclude that
the map 
$$
\vph_\ga: f(X)\mapsto f'(\pi_b)=
\ga^{-1}(b_\#)\lan f \e_b\tga^{-1}\mu_\circ\ran,
$$
induces the involution $\tau_+^{-1}\si\tau_+^{-1}=
\tau_-^{-1}$
on $\HH^{\flat}.$ Note that
$\ga(b_\#)$ is nothing else but
$\ga(\pi_b)$ evaluated at $\xi=-\rho_k.$
The map $\vph_\ga$ acts from $\Q_{q,t}[X_b]$
to the $\HH^{\flat}\ $-module $\FF_\#,$ 
where the field of constants is 
extended by  $q^{-(\xi,\xi)/2}.$ 

The automorphism $\tau_-$ fixes the $Y$-operators. Hence
the image of $f=\e_c$ is an eigenfunction of the discretizations
$\de(Y_a)$ of the $Y_a (a\in B)$ corresponding to the same set of
eigenvalues as for $\e_c.$ Let us prove that 
$\vph_\ga(\e_c)$ has to be proportional
to the discretization $\de(\e_c)$ of $\e_c.$

In the first place,
we may assume that $c=0$ employing the $Y$-intertwiners
(see (\ref{ehatb})).
The images of the $Y$-intertwiners with
respect to $\tau_-^{-1}$ can be exactly calculated 
but we do not need explicit formulas here.
It is sufficient to know that they are invertible operators
acting in $\FF_\#.$ The function $g=\vph_\ga(1)$ has additional
symmetries: $t_i^{-1/2}\de(T_i)(g)=$ $g=\de(\pi_r)(g).$ 
We have used that $\tau_-$ fixes $T_i$ and $\pi_r$ for all $i,r$
(see (\ref{taumin})). This readily leads to the
$\hW^{\flat}$-invariance of $g,$ which means that it 
has to be constant.  

Setting $\vph_\ga(\e_c)=
h_c\ga(\pi_c)q^{-(\rho_k,\rho_k)}\de(\e_c)$ 
for $h_c\in \Q_{q,t},$
we need to check that $h_c=1.$ It is true for $h_0.$ 
However
$h_c=h_b$ for all $c,b\in B$
because both the left-hand and the right-hand 
sides of (\ref{epep}) are 
$b\leftrightarrow c$ symmetric. The first symmetry is obvious, 
the 
second is the duality relation
(\ref{ebdual}). 

\vskip 0.2cm
The second formula can be easily deduced from the first
thanks to (\ref{ebeconjj}): 
\begin{align}
&\e_c^*\ =\ \prod_{\nu} t_\nu ^{-l_\nu(w_0)/2}\, 
T_{w_0}(\e_{\varsigma(c)}),\ \varsigma(c)=-w_0(c).
\notag \end{align}
Indeed, $T_{w_0}$ commutes with $\ga$ and its
Fourier transform is the discretization of
$\si(T_{w_0})=T_{w_0}.$ Applying this argument one
more time we get (\ref{epepss}). The (conjugation
of the) latter formula will be important to construct
the inverse Fourier transform in the next section.
A direct proof of (\ref{epeps}), without the conjugation
formula, is also not difficult.
Here it is.

First, $\e_c^*$ are eigenfunctions of the operators
$\eta(Y_a),$ where $\eta = \vep\si$ (see above).
The images of these operators under $\vph_\ga$ are
$\tau_-^{-1}(\eta(Y_a).$ One has: 
$$
\tau_-^{-1} \eta=
\tau_-^{-1}\vep\si=\vep\tau_+^{-1}\si=
\vep\tau_-^{-1}\tau_+=
\vep\si\tau_-=\eta\tau_-.
$$
Therefore  $\vph_\ga$ fixes $\eta(Y_a)$ and 
$\vph_\ga(\e_c^*)$
is proportional to $\de(\e_c^*)$ for any $a,c\in B.$

Second, we introduce $\hat{h}$ from $\vph_\ga(\e_c^*)=
\hat{h}_c\ga(\pi_c)q^{-(\rho_k,\rho_k)}\de(\e_c^*).$
However now is not obvious that the left-hand side of 
(\ref{epeps})
is symmetric ($a\leftrightarrow b$) as it was for 
(\ref{epep}).

We need to involve
$$
\psi_\ga: f(X)\mapsto 
\ga^{-1}(b)\lan f^* \e_b\tga^{-1}\mu_\circ\ran.
$$
It induces the automorphism $\ \tilde{\vep}\ $ 
$ =\ \tau_+^{-1}\vep\tau_+.$ Indeed,
$\lan f^* \e_b\tga^{-1}\mu_\circ\ran\ =$ 
$\lan (f\tga)^* \e_b\mu_\circ\ran.$

Since
$$
\eta(Y_a)=\vep\tau_-^{-1}\tau_+\tau_-^{-1}(Y_a)=
\tau_+^{-1}\vep\tau_+(Y_a)=\tilde{\vep}(Y_a),
$$
$\psi_\ga(\e_c)$ is proportional
to $\de(\e_c^*).$ 
Therefore letting  
$$\psi_\ga(\e_c^*)=
\tilde{h}_c\ga(\pi_c)q^{-(\rho_k,\rho_k)}\de(\e_c^*),
$$
we conclude that $\hat{h}_c=\tilde{h}_b$ for all
$b,c\in B.$ Hence both functions  $\hat{h}, \tilde{h}$
are constants and, moreover, equal $1$ thanks to the 
normalization of $\e_b.$

The explicit formula for $\lan \tga^{-1}\mu_\circ\ran$ 
was calculated in [C5] using the shift operators and
the analytic continuation.
\qed

We remark that the product $\mu\tga^{-1}$  generalizes
the (radial) Gaussian measure  in the theory of 
Lie groups and symmetric spaces. To be more exact, 
the ``noncompact'' case, which
will be considered next, is such a generalization.
The above setup can be called ``compact'', because 
taking the constant term corresponds to the 
integration with respect to the imaginary period. 

We would like to mention that the appearance of the 
nonsymmetric
polynomials has no known counterparts in the classical 
representation
theory even in the so-called group case when $k=1,$ i.e. 
$t=q.$ In this case, the symmetric Macdonald polynomials 
become the characters of the compact Lie groups,
and one can expect the nonsymmetric polynomials to be
somehow connected with these representations. However
it does not happen. There are certain relations of
the (degenerate) nonsymmetric polynomials with the
Demazure character formulas (the Kac-Moody case, 
essentially basic representations). However it looks
more accidental than conceptual. It merely reflects
the fact that both constructions are based on the
Demazure operations.

In the group case,
the symmetrizations of (\ref{epep}) and (\ref{epeps}) 
can be readily deduced from the Weyl character formula.
It seems that these formulas  were not used in 
the harmonic analysis.
Maybe because they cannot be extended to the Harish-Chandra
zonal spherical functions on general symmetric spaces.

\vskip 0.2cm
\section {Jackson integrals}
\setcounter{equation}{0}

The formulas from the previous sections can be generalized
for Jackson integrals taken instead of the constant term
functional (corresponding to the imaginary integration 
over the period). It is a variant of the classical noncompact case.
Another variant is a straightforward analytic integration in the 
real direction, which will not be discussed in the paper.
We keep the same notation, expend all functions 
in terms of non-negative powers of $q,$ and consider
$q^{(\xi,b_i)}$ as independent parameters. One may also treat
$\xi$ as complex vector in a general position 
assuming that $|q|<1.$

The Jackson integral of $f(bw)\in \FF[\xi]$ is
$$ 
\lan f\ran_\xi\ \equal\ \sum f(bw), \where b\in B,\ 
w\in W.
$$ 

Recall that  $X_a(q^z)=q^{(a,z)},\ X_a(bw)=q^{(a,w(\xi)+b)},\ 
\ga(q^z)=q^{(z,z)/2},$ and $\ga(bw)=\ga(q^{w(\xi)+b}).$
Thus 
$$|W|^{-1}\langle \ga\rangle_\xi= 
\sum_{a\in B} q^{(\xi+a,\xi+a)/2}=
\tga^{-1}(q^\xi)q^{(\xi,\xi)/2}$$
for $\tga$ from (\ref{gauser}). 

The Jackson inner product already appeared in 
(\ref{innerdel}). It is
$\langle f,g\rangle_\bullet\ = 
\lan fg^*\mu_\bullet\ran_\xi$ 
for $\mu_\bullet$ from (\ref{muval}). To be more precise, here we 
need the discretization $\de(\mu_\bullet).$

The {\it Fourier-Jackson transform} $\vph_\bullet,$  the   
{\it skew Fourier-Jackson transform} $\psi_\bullet,$
and their bar-counterparts
are as follows:

\begin{align}
&\vph_\bullet(f)(\pi_b)\ =\ 
\lan f(\hw)\e_b(\hw)\mu_\bullet\ran_\xi,
\ b\in B,\notag\\  
&\psi_\bullet(f)(\pi_b)\ =\ 
\lan f(\hw)^*\e_b(\hw)\mu_\bullet\ran_\xi=
\lan \e_b,f\ran_\bullet,
\label{fouriers}
\\ 
&\bar{\vph}_\bullet(f)(\pi_b)\ =\ 
\lan f(\hw)\e_b^*(\hw)\mu_\bullet(\hw)\ran_\xi,\ 
b\in B,\notag\\ 
&\bar{\psi}_\bullet(f)(\pi_b)\ =
\ \lan f(\hw)^*\e_b^*(\hw)\mu_\bullet(\hw)\ran_\xi
=\lan f,\e_b\ran_\bullet.
\label{fourierj}
\end{align}

These transforms  act from subspaces of $\FF[\xi]$ to proper
completions of the
space $\F_\#,$ provided the convergence. 
The involution $\ast$ is  the conjugation of
the values of functions  $f\in \FF[\xi].$ It 
is well defined because the values are rational functions
in terms of $q,t$ (and their certain fractional powers).
One has:
$(f^*)(\hw)= (\de(f^*))(\hw)= (\de(f)(\hw))^*= 
(f(\hw))^*.$
For instance, $\chi_{\hw}^*=\chi_{\hw},$ and 
$\de_{\hw}^*=\de_{\hw}$ 
since $\mu_\bullet(\hw)^*=\mu_\bullet(\hw).$ 
Thus in (\ref{fourier}), (\ref{fouriers})
we can replace $\e_b^*(\hw)$ by 
$\e_b(\hw)^*.$ 

It is straightforward to check that the corresponding
automorphisms of 
$\HH^{\flat}\ $ are the same as in
Theorems \ref{ISOVD} and \ref{ISOVDC}. We simply replace  
$\lan\ \ran$ by  $\lan\ \ran_\xi:$ 
\begin{align}
&\vph \leftrightarrow \si,\ \bar{\vph} \leftrightarrow 
\si^{-1},
 \psi \leftrightarrow \vep,\ \bar{\psi} \leftrightarrow \eta\si.
\label{fouraut}
\end{align}

These automorphisms commute with the anti-involution $\star,$
so all transforms are $\star$- unitary up to a 
coefficient of
proportionality provided the convergence and
the irreducibility of the corresponding $\HH^\flat$-modules.
Here $\xi$ is generic. For the spherical specialization we
can be more precise.

\begin{theorem}
i) For $\xi=-\rho_k,$ the Fourier transforms from 
(\ref{fourier}), (\ref{fouriers})
act from the space of delta-function $\F_\#$ to the discretization
of the space $\Q_{q,t}[X_b]$ upon the restriction to the set
$\pi_B.$ They send
the form $\lan\, ,\, \ran_\#$ to $\lan\, ,\, \ran_\circ,$
and satisfy the relations 
\begin{align}
&\bar{\vph}_\circ\cdot \vph_\bullet=\hbox{id}=
\vph_\circ\cdot\bar{\vph}_\bullet,\notag\\ 
&\bar{\vph}_\bullet\cdot \vph_\circ=\hbox{id}=
\vph_\bullet\cdot\bar{\vph}_\circ,\notag\\ 
&\psi_\circ\cdot\psi_\bullet=\hbox{id}=
\bar{\psi}_\circ\cdot\bar{\psi}_\bullet,\notag\\ 
&\psi_\bullet\cdot\psi_\circ=\hbox{id}=
\bar{\psi}_\bullet\cdot\bar{\psi}_\circ.
\label{fourinver}
\end{align}
ii) Setting $\widehat{f}=\vph_\circ(f)\in \F_\#$
for $f(X)\in \Q_{q,t}[X_b],$ the inverse transform reads 
\begin{align}
&f(X)\ =\ \prod_{\nu} t_\nu ^{-l_\nu(w_0)/2}
\lan\, \widehat{f}\ T^\bullet_{w_0}(\e_{\varsigma(b)})
\mu_\bullet\ran_\#,
\label{fpibh}
\end{align}
where $T^\bullet$ acts on $\e_{b}$ via formula 
(\ref{Tfdel})
for $\de_{b}^\pi=\de_{\pi_b}:$
\begin{align}
&T^\bullet_i(\e_b) = 
\frac{ t_i^{1/2}q^{(\al_i,b_\#)} - t^{-1/2}}{
  q^{(\al_i,b_\#)} - 1  }\e_{s_i(b)}
-\frac{t_i^{1/2}-t_i^{-1/2}}{
  q^{(\al_i,b_\#)} - 1  }\e_b,
\label{Tbull}
\end{align}
extended to $W$ naturally. Similarly,
we have the Plancherel formula
\begin{align}
&\lan f\, g \mu_\circ\ran =
t_\nu ^{-l_\nu(w_0)/2}
\lan \widehat{f}\, T_{w_0}(\widehat{g}^\varsigma)
\mu_\bullet\ran_\#\notag\\ 
&\for f,g\in \Q_{q,t}[X_b],\ h(\hw)^\varsigma=
h(\varsigma(\hw)).
\label{Planch}
\end{align}
\label{INVER}
\end{theorem}
{\it Proof} of i) is based on (\ref{fouraut}) and the
irreducibility of the polynomial representation and
its delta-counterpart $\F_\#.$ The same argument gives
the unitarity of the Fourier transforms under 
consideration
up to proportionality. 
Thanks to the normalization of $\mu_\circ$ and 
$\mu_\bullet,$
$$\lan \de_0 ,\de_0 \ran_\#\ = \ 1\ =\  
\lan 1 , 1 \ran_\circ.
$$
So the coefficient of proportionality is $1.$

The explicit inversion formula from ii) is straightforward:
\begin{align}
&\vph_\circ(\prod_{\nu} t_\nu ^{-l_\nu(w_0)/2}
\lan \widehat{f}\  T^\bullet_{w_0}
(\e_{\varsigma(b)})\mu_\bullet\ran_\#)\notag\\ 
&=\prod_{\nu} t_\nu ^{-l_\nu(w_0)/2}
\lan \widehat{f}\  T^\bullet_{w_0}
\phi_\circ(\e_{\varsigma(b)})\mu_\bullet\ran_\#\notag\\ 
&=\lan \widehat{f}\  T^\bullet_{w_0}
(\prod_{\nu} t_\nu ^{-l_\nu(w_0)/2}
\phi_\circ(\e_{\varsigma(b)}))\mu_\bullet\ran_\#\notag\\ 
&=\lan \widehat{f}\  T_{w_0}
T^{-1}_{w_0}(\de_{b})\mu_\bullet\ran_\#=\widehat{f}
(\pi_b).
\label{phihatf}
\end{align}
Here we use (\ref{isovph}). Let us check
the Plancherel formula. The anti-involution corresponding
to the symmetric form $\lan f\, g\mu_\circ\ran$ 
is 
$$\diamond\equal\eta\cdot\star=\star\cdot \eta.
$$ 
The one for the form 
$\lan \widehat{f}\, T_{w_0}(\widehat{g}^\varsigma)
\mu_\bullet\ran_\#$
can be readily calculated as follows:
\begin{align}
&\lan\, \widehat{f}\ T_{w_0}
(\vph_\circ(H(g)){}^\varsigma)
\mu_\bullet\ran_\# 
=\lan\, \widehat{f}\ T_{w_0}\si(H^\varsigma)
(\widehat{g}^\varsigma)
\mu_\bullet\ran_\# \notag\\ 
=&\lan\, (T_{w_0}\si\varsigma(H)T_{w_0}^{-1})^\diamond
(\widehat{f})\ T_{w_0}(\widehat{g}^\varsigma)
\mu_\bullet\ran_\#\notag\\ 
=&\lan\, (\si^{-1}(H))^\diamond
(\widehat{f})\ T_{w_0}(\widehat{g}^\varsigma)
\mu_\bullet\ran_\#\notag\\ 
=&\lan\, \si(H^\diamond)
(\widehat{f})\ T_{w_0}(\widehat{g}^\varsigma)
\mu_\bullet\ran_\#=
\lan\,(H^\diamond(f))\, \widehat{}\,\,\, 
 T_{w_0}(\widehat{g}^\varsigma)
\mu_\bullet\ran_\#.
\notag \end{align}
Here we use that the Fourier transform $\vph$ induces 
$\si$
on $\HH\ ,$ the relations (\ref{etay}):
$$T_{w_0}\si\varsigma(H)T_{w_0}^{-1}=\si^{-1}(H),$$
and (see the last line) formula 
$$\si^\diamond=\diamond\cdot\si\cdot\diamond=
\eta\si\eta=
\vep\si\vep=\si^{-1}.
$$ 
\qed

%\vskip 0.2cm
Part ii) of the Theorem generalizes the main theorem
about $p$-adic spherical transform due to Macdonald 
(symmetric case) and Matsumoto 
(see [Ma] and resent [O4,O5]). The
classical $p$-adic spherical transform acts 
from $\F_\#$ to the
polynomials. So we need to reverse ii). 
The $\e^*_b$ generalize the $p$-adic
(nonsymmetric) spherical functions due to Matsumoto.
The limit $q\to \infty$ of the
following Corollary is exactly the theory of spherical
Fourier transform. The case of generic $\xi$ (in place of
$-\rho_k$) is presumably connected with the general 
(non-spherical)
Fourier transform on affine Hecke algebras [KL1] and is
expected to be directly related to [HO2]. 

Note that the
latter paper is devoted to affine Hecke algebras with
arbitrary labels. Only equal labels are considered in [KL1].
There are Lusztig's papers towards nonequal labels. However
the general labels  unfit the ``geometric'' methods.

%\vskip 0.2cm
\begin{corollary}
Setting 
$\widehat{f}=\bar{\vph}_\bullet(f)\in \Q_{q,t}[X_b]$
for $f\in \F_\#,$ its inversion and the
Plancherel formula are as follows:  
\begin{align}
&f\ =\ \prod_{\nu} t_\nu ^{l_\nu(w_0)/2}
\lan\, \widehat{f}\  T^{-1}_{w_0}
(\e^*_{\varsigma(b)})\mu_\circ\ran,
\label{fpibhi}\\
&\lan\, f\ g \mu_\bullet\ran_\#\ =\
t_\nu ^{l_\nu(w_0)/2}
\lan\, \widehat{f}\ T_{w_0}^{-1}(\widehat{g}^\varsigma)
\mu_\circ\ran\notag\\ 
&\for f,g\in \Q_{q,t}[X_b],\ 
h(X_b)^\varsigma=h(X^{-1}_{w_0(b)}).
\label{Planchi}
\end{align}
\label{INVERI}
\end{corollary}

%\vskip 0.2cm
{\it Gauss-Jackson integrals.} In a sense, these integrals
are the missing part of the Harish-Chandra theory
of spherical functions and the theory of the $p$-adic
spherical transform. The formulas for the integrals of 
the Gaussians with respect
to the Harish-Chandra (zonal) transform exist in the group case only
($k=1$). The Gaussians cannot be    
added to the $p$-adic theory as well.

Let us transfer Theorem \ref{EPEP} from the compact case,
with $\lan \ran$ as the integration, to the Jackson case.
The proof remains essentially unchanged. In fact, the
formulas below are $\ast$-conjugations of (\ref{epep}) with
a minor reservation about (\ref{mehtamu}).
See also Theorem 7.1 from [C5]. 

\begin{ theorem}
Given $b,c\in P$ and the corresponding spherical polynomials 
$\e_b,\e_c$,
\begin{align}
\ \ \langle
\e_b\, \e_c^*\,\ga{\mu_\bullet}
\rangle_\xi\ &=\
 q^{-(b_\#,b_\#)/2-(c_\#,c_\#)/2 +(\rho_k,\rho_k)} 
\e_c(q^{b_\#})\langle \ga{\mu_\bullet}\rangle_\xi,
\label{eejack} \\ 
\ \ \langle \e_b^*\, \e_c^*\,\ga \mu_\bullet
\rangle_\xi\ &=\
 q^{-(b_\#,b_\#)/2-(c_\#,c_\#)/2 +(\rho_k,\rho_k)} 
\e_c^*(q^{b_\#})\langle \ga{\mu_\bullet}\rangle_\xi,
\label{eejacks}
\end{align}
\begin{align}
\ \ \langle \e_b\, \e_c \ga\mu_\bullet
\rangle_\xi\ &=\ 
q^{-(b_\#,b_\#)/2-(c_\#,c_\#)/2 +(\rho_k,\rho_k)}\notag\\  
&\times\prod_{\nu}t_\nu ^{l_\nu(w_0)/2}
T^{-1}_{w_0}(\e_c)(q^{b_\#})\langle 
\ga\mu_\bullet\rangle_\xi,
\label{eejac}
\\ 
\langle \ga{\mu_\bullet}\rangle_\xi= 
&\mu(q^\xi,t^{-1})\, |W|^{-1}\langle \ga\rangle_\xi\ 
 \prod_{\al\in R_+}\prod_{ j=0}^{\infty}
\Bigl(\frac{ 
1- t_\al^{-1}q_\al^{-(\rho_k,\al^\vee)+j}}{
1- q_\al^{-(\rho_k,\al^\vee)+j} }\Bigr). 
\label{memujack}
\end{align}
Here by $\mu(q^\xi,t^{-1})$ we mean the right-hand side of
(\ref{mu}) evaluated at $X=q^\xi$
where all $t_\al$ are replaced by $t_\al^{-1}.$
\label{EEJACK}
\end{ theorem}
\qed

In these formulas, $t$ is generic or complex provided that
the discretization $\de(\mu_\bullet)=$ 
$\mu_\bullet(bw)$ is well defined and
$\e_{b}, \e_c$ exist. Note that the right-hand side of 
(\ref{memujack}) has to 
be replaced by the limit when $k_\nu\in \Z_+$.

We come to the inversion and Plancherel formulas.
The pairings $\lan\, f \, g\mu_\circ \ran$ and 
$\lan\, f \, g\mu_\bullet\ran_\xi$ can be naturally extended
to the $\HH^\flat $-modules $\Q_{q,t}[X_b]\tga^{-1}$ and
the discretization-image 
$\de(\Q_{q,t}[X_b]\ga)$ of $\Q_{q,t}[X_b]\ga$ 
in $\FF[\xi].$ We will use them for the Plancherel formulas.

We set $c_\xi=\lan\ga\mu_\bullet\ran_\xi.$
The Fourier transforms 
$$\vph^c\ =\ c^{-1}_\xi\vph_\bullet \and 
\bar{\vph}^c_\bullet\ =\
c^{-1}_\xi\bar{\vph}_\bullet
$$
transfer 
$$\de(Q_{q,t}[X_b]\ga)\ \to \
\Q_{q,t}[X_b]\tga^{-1}.
$$
Respectively, $\psi^c_\bullet, \bar{\psi}^c_\bullet$
preserve the first module.

In the compact case, setting  $c=\lan\ga\mu_\circ\ran$,
$$\vph^c_\circ=c^{-1}\vph_\circ,\ 
\bar{\vph}^c_\circ=c^{-1}\bar{\vph}_\circ,\ 
\psi^c_\circ=c^{-1}\psi_\circ, \
\bar{\psi}^c_\circ=c^{-1}\bar{\psi}_\circ,
$$ 
the first two transforms act from the second module to the first
and the last two preserve the second module.

They satisfy the counterparts of the inversion formulas
(\ref{fourinver}):
\begin{align}
&\bar{\vph}^c_\circ\cdot \vph^c_\bullet=\hbox{id}=
\vph^c_\circ\cdot\bar{\vph}^c_\bullet,\notag\\ 
&\bar{\vph}^c_\bullet\cdot \vph^c_\circ=\hbox{id}=
\vph^c_\bullet\cdot\bar{\vph}^c_\circ,\notag\\ 
&\psi^c_\circ\cdot\psi^c_\circ=\hbox{id}=
\bar{\psi}^c_\circ\cdot\bar{\psi}^c_\circ,\notag\\ 
&\psi^c_\bullet\cdot\psi^c_\bullet=\hbox{id}=
\bar{\psi}^c_\bullet\cdot\bar{\psi}^c_\bullet.
\label{fourinverj}
\end{align}
It is straightforward to reformulate part ii) of Theorem
\ref{INVER} for $\vph^c$ and the pairings
$\lan\, f \, g\mu_\circ \ran\ $ on 
$\Q_{q,t}[X_b]\tga^{-1}$  
and  
$\lan\, f \, g\mu_\bullet\ran_\xi\ $ on 
$\de(\Q_{q,t}[X_b]\ga).$ 

Theorem \ref{EEJACK}
and the above facts remain valid when $\xi=-\rho_k.$
In this case, 
\begin{align}
&\mu\ = \mu^{(k)}=\prod_{\al \in R_+}
\prod_{i=0}^\infty \frac{(1-X_\al q_\al^{i}) 
(1-X_\al^{-1}q_\al^{i+1})
}{
(1-X_\al t_\al q_\al^{i}) 
(1-X_\al^{-1}t_\al^{}q_\al^{i+1})},\
\notag \end{align}
\begin{align}
&\mu(q^{-\rho_k},t^{-1})=\prod_{\al \in R_+}
\prod_{i=0}^\infty \frac{
(1-q_\al^{-(\rho_k,\al^\vee)+i}) 
(1-q_\al^{ (\rho_k,\al^\vee)+i+1})
}{
(1- t_\al^{-1} q_\al^{-(\rho_k,\al^\vee)+i}) 
(1- t_\al^{-1} q_\al^{ (\rho_k,\al^\vee)+i+1})},
\label{muu}
\\ 
&\mu_\bullet(bw)= \mu_\bullet(q^{-w(\rho_k)+b}) =
\prod_{[\al,j]\in \la'(bw)}
\Bigl(
\frac{
t_\al^{-1/2}-t_\al^{1/2} q_\al^{(\al^\vee,\rho_k)+j}}{
t_\al^{1/2}-t_\al^{-1/2} q_\al^{(\al^\vee,\rho_k)+j}
}
\Bigr),
\label{muvalr}
\end{align}
for $\la'(bw)=\{[-\al,j]\ |\ [\al,\nu_\al j]\in \la(bw)\}$,
$\la(bw)=R_+^a\cap (bw)^{-1}(R_-^a).$
The function $\mu_\bullet(bw)$ is always well defined
and nonzero only as $\pi_b=b u_b^{-1}$
(for generic $q,t$).
See (\ref{jbset}) and (\ref{muval}).
In this case, all $\al$ in the product (\ref{muvalr}) are
positive.

Formula (\ref{memujack}) now reads as follows:
\begin{align}
& \langle \ga\mu_\bullet\rangle_\#\ =\ 
|W|^{-1}\langle \ga\rangle_{\#} 
\prod_{\al\in R_+}\prod_{ j=1}^{\infty}\Bigl(\frac{ 
1- q_\al^{(\rho_k,\al^\vee)+j}}{
1-t_\al^{-1}q_\al^{(\rho_k,\al^\vee)+j} }\Bigr). 
\label{hatmu}  
\end{align}
Recall that $\lan\ \ran_\xi$ for $\xi=-\rho_k$ is denoted by
$\lan\ \ran_\#.$
The inversion formulas for $\vph$ get simpler:
\begin{align}
&\bar{\vph}_\circ\cdot \vph_\bullet\ =\
|W|^{-1}\lan \ga\ran_\#\,\hbox{id}\ =\
\vph_\circ\cdot\bar{\vph}_\bullet,\notag\\ 
&\bar{\vph}_\bullet\cdot \vph_\circ\ =\
|W|^{-1}\lan \ga\ran_\#\,\hbox{id}\ =\
\vph_\bullet\cdot\bar{\vph}_\circ.
\label{invervj} 
\end{align}

\vskip 0.2cm
{\bf Macdonald's $\eta$-identities.}
As a by-product, we can represent $\lan \ga\ran_\#$
as a theta-like product times a certain finite sum, 
provided that the 
left-hand side of (\ref{hatmu}) contains 
finitely many nonzero terms.
It happens when a certain $\Z_+$-linear combination of 
$k_{\sht},k_{\lng}$
is from $-\N.$
The main example (which will be used later) is as follows. 
Let $q$ be generic such that $|q|<1.$

We call a root $\al\in R_+$ {\it extreme}
if $h_\al(k)\equal (\rho_k,\al^\vee)+k_\al$ does not coincide 
with any $(\rho_k,\be^\vee)$ for positive roots $\be,$
and {\it strongly extreme} if also $(\al,\om_i)>0$ for all 
$1\le i\le n.$  
Note that $h_\vth(k)$ is $kh$ 
for the {\it Coxeter number}
$h=(\rho,\vth)+1$ in the case of coinciding $k.$  

Here $k_{\sht},k_{\lng}$ are treated as independent parameters.
We omit the subscript of $k$ in the simply-laced case:
$\rho_k=k\rho,\ $ $ h_\al(k)= ((\rho,\al)+1)k.$
In this case, there is only one extreme root, namely, 
$\vth.$
Let us list the extreme roots for the other root systems
(the notation is from [B]). All of them are short:

$B_n$) all short roots $\ep_i\ (1\le i\le n)$ where
$h_{\ep_i}(k)=2k_{\sht}+2(n-i)k_{\lng};$

$C_n$) $\vth=\ep_1+\ep_2$ with $h_{\vth}(k)=
2(n-1)k_{\sht}+2k_{\lng}$ and $\al=\ep_1-\ep_n$ with 
$h_\al(k)=$
$n\, k_{\sht};$

$F_4$) 
$0011, 0121, 1121, \vth=1232$ with $h$ equal to 
$3k_{\sht}$, $4k_{\sht}+2k_{\lng}$ 
 $4k_{\sht}+4k_{\lng}$, $6k_{\sht}+6k_{\lng}$ respectively;

$G_2$) $\al_1\ (h_{\al_1}(k)=2k_{\sht}),\ 
\vth=2\al_1+\al_2\
(h_\vth(k)=3(k_{\sht}+k_{\lng})).$

Thus $\vth$ is a unique strongly extreme root
for all root systems except for  $F_4.$ 
The root $1121$ for $F_4$ is strongly extreme too. 

The definition can be modified
by imposing one of the following conditions
a) $k_{\lng}=k_{\sht}\ $, b) $k_{\sht}=0,\ $ 
c) $k_{\lng}=0.$ Under 
either constraint,
$\vth$ becomes a unique strongly extreme root. The below theorem
holds in these cases.

\begin{theorem}
i) Let us assume that $q$ is generic, $|q|<1$, and
\begin{align}
&(\rho_k,\al^\vee)\ \not\in\ \Z \setminus\{0\}
\hbox{\ for\ all\ } \al\in R_+.
\label{munonz}
\end{align}
These conditions result in 
$$(\rho_k\,,\,\al^\vee)-k_\al\ \not\in\ \Z \setminus\{0\}
$$
for $\al\in R_+$ and therefore are sufficient for 
the right-hand side of (\ref{hatmu}) to exist and be nonzero
as well as for
the sum in the definition of $\langle \ga\mu_\bullet \rangle_\#.$ 
The latter sum is finite if and only if
$h_\al(k)\in -\N$ for a strongly extreme root $\al$.

ii) If $h_\vth(k)=-1$ and (\ref{munonz}) is satisfied (which 
automatically holds as $k_{\sht}=k_{\lng}$),
the measure $\mu_\bullet(bw)$ takes exactly 
$|\Pi^{\flat}|$
nonzero values. All of them equal  $1$ and  
\begin{align}
& \langle \ga \rangle_\#\ =\ \sum_{b\in B}
q^{(b+\rho_k,b+\rho_k)/2}\ 
=\
\Bigl(\sum_{\pi_r\in \Pi^{\flat}}
q^{(\om_r+\rho_k,\om_r+\rho_k)/2}\Bigr)\times \notag\\ 
& \prod_{\al\in R_+}\prod_{ j=1}^{\infty}\Bigl(\frac{ 
 1-t_\al^{-1}q_\al^{(\rho_k,\al^\vee)+j} }{
1- q_\al^{(\rho_k,\al^\vee)+j} }\Bigr), \where \om_0=0. 
\label{gakneg}  
\end{align}
In the exceptional case 
$\{R=F_4,\, \al=1121,\, h_\al(k)=-1\},$
the number of nonzero values of $\mu_\bullet(bw)$ 
is greater than $|\Pi|=1.$ 
\label{GANEG}
\end{theorem}

{\it Proof.} 
Concerning the first claim, 
$(\rho_k\,,\,\al^\vee)-k_\al$ is zero for simple $\al=\al_i.$
If $\al$ is not simple, then there exist a simple root $\al_i$
of the same length as $\al>0$ such that 
$\al^\vee-\al_i^\vee$ is a positive coroot $\be^\vee.$
Then $(\rho_k,\be^\vee)=(\rho_k\,,\,\al^\vee)-k_\al$ 
cannot be from $\Z\setminus \{0\}$ due to (\ref{munonz}). 

There must be at least one root $\al$ satisfying
$h_\al(k)\in -\N$ to make the Jackson sum
$
\langle \ga\mu_\bullet \rangle_\#\ =\ \sum_{b\in B} 
q^{(b_\#,b_\#)/2} \mu_\bullet(\pi_b)
$  
finite. 
Assumption (\ref{munonz}) gives that at least 
one such $\al$ has to be strongly extreme.
By the way,  
if there are two such roots $\al\neq\be$ 
(this may happen for $F_4$ only)
then $h_\be(k) -h_\al(k) \in \Z,$ which contradicts 
(\ref{munonz}).

Hence the above sum is finite. 
Indeed, recall that $\mu_\bullet(\hbox{id})\ =\ 1$ and 
\begin{align}
&\mu_\bullet(\pi_b)\ =\ \prod
\Bigl(
\frac{
t_\al^{-1/2}-t_\al^{1/2} q_\al^{(\al^\vee,\rho_k)+j}}
{t_\al^{1/2}-t_\al^{-1/2} q_\al^{(\al^\vee,\rho_k)+j}
}
\Bigr) \for \al\in R_+,
\label{muhatpi}
\\ 
-( b_-, \al^\vee )&>j> 0 
\iif  u_b^{-1}(\al)\in R_-,\ 
-( b_-, \al^\vee )\ge j > 0 \hbox{\ otherwise}. 
\notag \end{align}
See (\ref{muvalr}) and (\ref{jbset}). 

Let us assume now 
that $\al=\vth,\ h_{\vth}(k)=-1,$ and $b\neq 0.$ 
The case of $F_4, 1121$ is
left to the reader.
Then $\mu_\bullet(\pi_b)$ can be nonzero only for 
$b$ from the $W$-orbits
of the minuscule weights $b_-=-\om_r.$
Really, otherwise $-(b_-,\vth)\ge 2,\ $ 
$[\vth,1]\in \la'(\pi_b),\ $ 
and the product (\ref{muhatpi}) is zero.
Moreover, $u_b^{-1}$ sends the roots 
$\{\al^\vee_i,i\neq r,i>0\}$
to $R^\vee_+$ since it is minimal such that 
$u_{b}^{-1}(\om_r)=b_-.$
Also $u_b^{-1}(-\vth)\in R_+^\vee$ 
to make the product (\ref{muhatpi})
nonzero.

However $\{-\vth,\al^\vee_i,i\neq r,i>0\}$ form a 
basis of $R^\vee.$
Therefore these conditions determine $u_{b}^{-1}$ 
uniquely and 
it has to coincide with $u_{b_+}^{-1}.$ I.e. 
$b=b_+=w_0(-b_-)=\om_{r^*},$  
the $\la'$-set of  $\pi_{b_+}=\pi_{r^*}$ is empty, and 
$\mu_\bullet(\pi_{r^*})=1.$ See  (\ref{pi}).

Thus it suffices to evaluate the Gaussian:
$$
((\om_r)_\#,(\om_r)_\#)=
(\om_r-u_r^{-1}(\rho_k)\,,\,\om_r-u_r^{-1}(\rho_k))=
(\om_{r^*}+\rho_k\, ,\,\om_{r^*}+\rho_k).
$$
\qed

Another method of proving (\ref{gakneg})
is based on the following symmetrization
of (\ref{hatmu}) due to  [C5], formula (1.11): 
\begin{align}
& \sum_{b_-\in B_-} q^{(b_--\rho_k\, ,\, b_--\rho_k)/2}
\De_\bullet(b_-)\ =\ 
\Bigl(\sum_{b\in B}  q^{(b+\rho_k\, ,\, b+\rho_k)/2}  \Bigr)
\times
\label{dehatpii}\\  
&
\prod_{\al\in R_+}\prod_{ j=1}^{\infty}
\Bigl(\frac{ 
1- q_\al^{(\rho_k,\al^\vee)+j}}{
1-t_\al^{-1}q_\al^{(\rho_k,\al^\vee)+j} }\Bigr),\
\ \De_\bullet(b_-)\ =\notag\\ 
&\prod
\Bigl(
\frac{
(t_\al^{-1/2}-t_\al^{1/2} 
q_\al^{(\al^\vee,\, \rho_k)+j-1})
(1-q_\al^{(\al^\vee,\, \rho_k-b_-)})}
{(t_\al^{1/2}-t_\al^{-1/2} 
q_\al^{(\al^\vee,\, \rho_k)+j})
(1-q_\al^{(\al^\vee,\, \rho_k)})
}
\Bigr).
\label{dehatpi}
\end{align}
The latter product is over the set from (\ref{muhatpi})
for $b=b_-,$ i.e. the relations 
$-( b_-, \al^\vee )\ge j> 0$ for all $\al\in R_+$ must hold.
In (\ref{dehatpi}), only
$b_-=-\om_r$ have nonzero $\De_\bullet(b_-).$
So it is sufficient  to check that
the latter are $1$ (which is true). This approach 
seems to be more convenient for computing explicit 
formulas when $h_\vth(k)=-m,\ m>1.$ They result in
generalizations of (\ref{gakneg} ).

We mention 
that when  $h_{\al}(k)=-m$ for an extreme but not strongly
extreme root $\al,$ 
then the sum for  $\lan \ga\mu_\bullet\ran_\#$
and, equivalently,  the  sum in (\ref{dehatpii}) 
are not finite.
However the corresponding variant of the theorem allows 
to reduce the summation dramatically.

Also note that the parameter $k_{\sht}$ 
is free in the case of $B,C,G,F,$ 
provided that it is in a general position.
For instance,
$\vth=\ep_1$ for $B_n$ and there is  only one constraint
$h_{\vth}(k)=2k_{\sht}+2(n-1)k_{\lng}=-1.$

The formula (\ref{gakneg}) is one of modifications
of the Macdonald identities [M6] closely related to
the Kac-Moody algebras (cf. [K],Ch.12). 
In the last section, we will connect it with one-dimensional 
representations of $\HH^{\flat}.$

\vskip 0.2cm
\section {Semisimple representations} 
\setcounter{equation}{0}

In this section we begin with representations $\HH^\flat$
with the cyclic vector and find out when they are
semisimple and pseudo-unitary. 
The latter means the existence of a nondegenerate hermitian
$\star$-invariant form. We add ``pseudo'' because
the involution $\ast$ acts on constants via
the complex conjugation only for $|q|=1$
and real $k.$ Also the form is not supposed to be positive. 
By definite, we mean a form with nonzero squares of the
eigenfunctions with respect to either $X$ or $Y$-operators.
So it is a substitute for the positivity.  
We will give necessary and sufficient conditions 
for induced modules to possess nonzero semisimple and
pseudo-unitary quotients.

It is worth mentioning that the main theorem generalizes
the classification of semisimple representations of affine 
Hecke algebra of type $A$ for generic $t.$
See e.g. Section 3 from [C9], references therein, and [Na]. 
In the $A$-case, the combinatorial part was clarified
in full. It is directly related to
the classical theory of Young's bases.

There are also partial results for classical root systems
in the author's works and  recent papers by A.Ram 
(he also considered some special systems, 
for instance, $G_2$). 
Still there is no complete answer. Generally, one can use
the classification of all irreducible representations 
from [KL1]. However checking the semisimplicity 
is far from immediate.

Actually the importance of semisimple representations 
is somewhat doubtful from the $p$-adic viewpoint, 
in spite of interesting combinatorial applications.
The reason is that the Bernstein-Zelevinsky operators
are not normal with respect to the natural unitary
structure which comes from the $p$-adic theory.
 
The theory of double affine Hecke algebras, especially  
as $|q|=1$, does require such representations.
In reasonably interesting representations,
either $X$-operators or $Y$-operators are normal.
The case of $A_{n-1}$ will be considered below.

\vskip 0.2cm
The notation is from of the previous sections. 
All $\HH^b$-modules (the lattice $B$ is fixed)
will be defined over the field $\Q_{q,t}$ 
or over its extension $\Q_\xi=$ $\Q_{q,t}(q^{\xi_1},
\ldots,q^{\xi_n}),$
where $q^{\xi_i}=q^{(b_i,\xi)}$ for a basis $\{b_i\}$ of $B.$
Here $q,t,\xi$ can be generic as above
or arbitrary complex. In the latter case, 
$\Q_\xi$ is a proper subfield of $\C,$
provided the existence of the involution $\ast$ on  
$\Q_{q,t}\subset \Q_\xi$ taking 
$$
q^{1/\tilde{2m}}\mapsto q^{-1/\tilde{2m}},\
t_\nu^{1/2}\mapsto t_\nu^{-1/2},\ q^{\xi_i}\mapsto q^{-\xi_i}.
$$
This involution is extended to the anti-involution 
(\ref{star})
of $\HH^\flat.$ We need several general definitions.

\vskip 0.2cm
{\it Eigenvectors.}
A vector $v$ satisfying 
\begin{align}
&X_a(v)\  =\ q^{(a,\xi)}\,v \hbox{ \ for \ all\ } 
a\in B.
\label{xcyc}
\end{align}
is called an 
$X$-eigenvector of weight $\xi.$ We set
\begin{align}
&V^s_X(\xi)\equal 
\{v\in V\ |\ (X_a-q^{(a,\xi)})^s(v)=0 \for a \in B\},
\notag\\ 
&V_X(\xi)\ =\ V^1_X(\xi),\ \, V_X^\infty(\xi)=
\cup_{s>0}\, V^s_X(\xi).
\label{vxxi}
\end{align}
Taking $Y_a$ instead of $X_a$ in (\ref{xcyc}), we define 
$Y$-eigenvectors of weight $\xi$ and introduce 
$V_Y^s(\xi).$

The action of the group $\hW$ on the weights is affine:
$\hw(\xi)=\hw\llb \xi\rrb.$ See (\ref{afaction}).
The $\hW^\flat$-{\it stabilizer} of $\xi$ is
\begin{align}
&\hW^\flat[\xi]\equal \{\hw\in \hW^\flat,\ 
q^{(\hw\llb\xi\rrb-\xi\, 
,\, a)}=
1\for \hbox{all\ } a\in B\}.
\label{hwflatxi}
\end{align}

From now on  weights $\xi$ and $\xi'$ will be identified if
$q^{(a,\xi'-\xi)}=1$ for all $a\in B.$ We put 
$q^\xi=q^{\xi'}$ if it
is true. 

The {\it category} $\mathcal{O}_X$ is formed by the 
modules $V$ such that
$V=\oplus_{\xi}V^{\infty}_X(\xi)$ and the latter spaces are
finite-dimensional, where the summation is over different $q^\xi.$
This is supposed to hold
for a proper extension of  the field of constants.
In the definition of $\mathcal{O}_Y,$ we substitute 
$Y$ for $X.$

Using the decomposition (\ref{hatdec}) from
the PBW-type Theorem \ref{FAITH}, one can checks that an 
$\HH^\flat-$
modules from $\mathcal{O}_X\cap \mathcal{O}_Y$ with finitely many
generators are finite-dimensional. Indeed, the space 
$\tU=\Q_{q,t}[Y_b] U$ is finite-dimensional for any 
finite-dimensional subspace $U$ of such
a module. If $U$ is preserved by the affine Hecke subalgebra 
$\h_X^\flat$
generated by $\Q_{q,t}[X_b] U$ and 
$T_i,\, 1\le i\le n,$ then $\tU$
is a finite sum of $\HH^\flat u$ for proper $u\in U.$ 
Therefore it
is an $\HH^\flat-$ submodule of $V.$ Since there are finitely
many generators of $V$ we get the required.

In this reasoning, it is 
not necessary to assume that the spaces $V_X^\infty(\xi)$ and
$V_Y^\infty(\xi)$ are finite-dimensional. If we know that
$V$ is finitely generated then it is sufficient to
use that  
\begin{align}
&\oplus_{\xi}V^{\infty}_X(\xi)\ =\ 
V=\oplus_{\xi}V^{\infty}_Y(\xi).\notag
\end{align}

\vskip 0.2cm
{\it Intertwiners.}
If the $X$-weight of $v$ is $\xi,$ then
$\Phi_{\hw} v$ is an $X$-eigenvector of weight
$\hw\llb \xi \rrb$ for the $X$-intertwiners $\Phi$ from
(\ref{Phi}). Here we can
take $S_{\hw}$ or $G_{\hw}$ instead of $\Phi_{\hw}$ (ibid.), 
because they are proportional to $\Phi_{\hw}.$ 
The denominators of $\Phi,S$ or $G$ have to be 
nonzero upon the evaluating at  $q^\xi$ if the latter
is not generic.
 
Given a reduced decomposition $\hw=\pi_r s_{i_l}\ldots s_{i_1},$
\begin{align}
&\Phi_{\hw} v= \pi_r\Phi_{i_l}(q^{\xi\{l-1\}})\cdots
\Phi_{i_1}(q^{\xi\{1\}})(v)\notag\\ 
&\for \xi\{0\}=\xi,\
\xi\{j\}=s_{i_p}\llb\xi\{p-1\}\rrb,\notag\\ 
&\Phi_i(q^\xi) \ =\ 
T_i + (t_i^{1/2}-t_i^{-1/2})(q^{(\xi+d\, ,\, 
\al_i)}-1)^{-1}, 
\ i\ge 0.
\label{Phiindu}
\end{align}
Cf. (\ref{Phijb}) and (\ref{ehatb}). 

In the $Y$-case, the $\Phi_i$ becomes $\vep(\Phi_i).$ 
For instance,
$q^\xi$ is replaced by $q^{-\xi}$ in the latter formula
as $j>0$ without touching $T,t.$ Constructing 
$\xi\{p\}$ in terms of the initial $Y$-weight $\xi$ 
remains unchanged.
There is another variant of the $Y$-intertwiners, 
technically more convenient.
It has been used in Proposition (\ref{PHIEB}). 
We will use it 
in the main theorem below. Namely, we 
can take  $\tau_+\bigl(\Phi_{\hw}(q^{-\xi})\bigr)$ 
instead of  
$\vep(\Phi_{\hw})(q^{\xi}).$ These operators
intertwine the $Y$-eigenvectors as well. 
Given $v$ of $Y$-weight $\xi,$
$v'=$ $\tau_+\bigl(\Phi_{\hw}(q^{-\xi})\bigr)v$ is a 
$Y$-eigenvectors of
weight $\hw\llb \xi\rrb.$ Notice the opposite sign of 
$\xi$ here.
Recall that $\tau_+$ can be interpreted as
the conjugation by the Gaussian fixing 
$X,$ $T_i(1\le i\le n),$ and the anti-involution $\star.$ 
Explicitly,
\begin{align}
&v'\ =\ \tau\bigl(\pi_r\Phi_{i_l}(q^{-\xi\{l-1\}})\cdots
\Phi_{i_1}(q^{-\xi\{1\}})\bigr)(v),
\label{Phiindy}
\end{align}
where the representation (and the notation) from 
(\ref{Phiindu})
is used.

\vskip 0.2cm
{\it Induced and semisimple modules.}
A $\HH^\flat$-module $V$ over $\Q_\xi$
is called $X$-{\it cyclic} of weight $\xi$
if it is generated by an element $v\in V$ 
satisfying (\ref{xcyc}).

The {\it induced} module generated by  $v$ with  (\ref{xcyc})
regarded as the defining relations 
is denoted by $\i_X[\xi].$ 
Its weights constitute the orbit  $\hW^\flat(\xi).$
Respectively, the weights of any cyclic modules of weight $\xi$
belong to this orbit.
The map $H\mapsto Hv$ identifies $\i_X[\xi]$
with the affine Hecke algebra $\h_Y^{\flat}\subset \HH^\flat.$
Note that the $\xi$-eigenspace $\i(\xi)$ of $\i=\i_X[\xi]$
becomes a subalgebra of $\h_Y^{\flat}$ under this identification. 

Similarly, the $Y$-cyclic modules 
are introduced for $Y_b$ instead of  $X_b:$ 
$$\i_Y[\xi]\cong
\h_X^{\flat}\ =\ \lan X_b,\, T_i\ \mid\ b\in B,i>0 
\ran.
$$  

By $X$-{\it semisimple}, we mean a $\HH^\flat$-module $V$  
coinciding with $\oplus_\xi V_X(\xi)$
over an algebraic closure of the field of constants. 
Here $\xi$ constitute the 
set of different $X$-weights of $V,$ 
denoted by $\hbox{Spec}_X (V)$
and called the $X$-{\it spectrum}.
If $V$ is a cyclic module of weight $\xi$
then the spectrum is defined over the field  $\Q_\xi.$
Indeed the spectrum of $\i_X[\xi]$ is the
orbit  $W^\flat\llb \xi\rrb.$ 
The definition of a $Y$-semisimple module is  for $Y$ taken
instead of $X.$

The functional representation $\F[\xi]$
introduced in Section 4
is $X$-cyclic for generic $\xi.$ 
It is isomorphic to $\i_X[\xi_o]$ for
the weight $\xi_o\in \hW\llb \xi\rrb,$ satisfying
the inequalities 
\begin{align}
& t_\al^{-1} X_{\tal} (q^\xi_o)  \neq 1\for \tal\in \tR_+,
\label{pluset}
\end{align}
provided that $\hW^\flat[\xi]=\{1\}.$ It is proved  in [C1] 
(after Corollary 6.5). 

The polynomial representation, which will be denoted
by $\v,$ is  $Y$-cyclic for generic $q,t.$
Any $E_b$ can be taken as its cyclic vector.
One has $\hbox{Spec}_Y (\v)=\{b_\# \mid b\in B\}.$

\vskip 0.2cm
{\it Pseudo-unitary structure.}
A quadratic form $(\ ,\ )$
is called respectively
$\ast$-bilinear and {\it pseudo-hermitian} if
$$
(ru,v)\ =\ r(u,v)\ =\ (u,r^*v) \for  r\in \Q_{q,t} \and
(u,v)\ =\ (v,u)^*.
$$  
We will always consider nondegenerate forms unless
stated otherwise.
  
A $\HH^\flat$-module $V$ is  {\it pseudo-unitary} if
it is equipped with a $\star$-{\it invariant} pseudo-hermitian form:
$$
(H(u),v)\ =\ (u,H^\star(v)) \hbox{ \ for \ all\ } u,v\in V,\ H\in 
\HH^\flat.
$$
We call it {\it definite} if $(v,v)\neq 0$ for all $v,$ and 
{\it $X$-definite} if 
$(v,v)\neq 0$ for all $X$-eigenvectors $v\in V$ assuming that 
the field of constants contain all eigenvalues.
By $X$-unitary, we mean 
$X$-semisimple $V$ with an $X$-definite form.
Strictly speaking, they should be  called
$X$-pseudo-unitary, but we skip ``pseudo'' in 
the presence of $X.$
The same definitions will be used for $Y.$

We do not suppose  $X$-definite or $Y$-definite forms to
be positive (or negative) hermitian forms. 
Note that $X$-unitary $V$ from the category 
$\mathcal{O}_X$ is 
$X$-semisimple. Vice versa, pseudo-unitary $X$-semisimple $V$ 
with the simple $X$-spectrum is $X$-unitary. 

When $\Q_\xi\subset \C,$ the involution $\ast$ is the restriction 
of the complex conjugation, and  $(u,u)>0$ for all 
$0\neq u\in V,$ 
then we call the form positive unitary
without adding ``pseudo'', $X,$ or $Y.$

Any irreducible $\HH^\flat$-quotients $\i_X[\xi]\to V$ can be defined
over $\C$ assuming that $|q|=1, k_\al\in \R, \xi\in \R^n$ in
the following sense. We claim
that there is a continuous deformation of
any triple $q,k,\xi$ to such a triple, which 
preserves the structure of
$\i_X[\xi],$ including quotients, submodules, semisimplicity, and 
the pseudo-hermitian invariant form. This simply means 
that the 
variety of all triples $\{q,k,\xi\}$ has points 
satisfying the above
reality conditions.
Thus $\ast$ can be assumed to be a restriction of 
the complex conjugation to the field of constants upon
such deformation.

\vskip 0.2cm
{\it Semisimple range.}
Given a weight $\xi,$ its {\it plus-range} $\Up_+[\xi]$
is the set of all elements
$\hw\in \hW^\flat$ satisfying the following condition
\begin{align}
& \tal=[\al,\nu_\al j]\in \la(\hw)\ \Rightarrow \ 
t_\al X_{\tal} (q^\xi)=q^{(\al\, ,\, \xi)+
\nu_\al(k_\al+j)}  \neq 1.
\label{upset}
\end{align}

We put $\hu\mapsto \hw$ if\ a) either $\hw=s_i\hu$
for $0\le i\le n$ provided that $l(\hw)= l(\hu)+1$ or\
b) $\hw=\pi_r\hu$ for  $\pi_r\in \Pi^\flat.$
The {\it boundary} $\ddot{\Up}_+[\xi]$
 of  $\Up_+[\xi]$ is by definition
the set of all $\hw\not\in \Up_+[\xi]$ such that $\hu\mapsto \hw$
for an element $\hu\in \Up_+[\xi].$ 
It may happen in case a) only. The condition
$l(\hw)= l(\hu)+1$ will be fulfilled automatically once  
$\hw=s_i\hu.$

The {\it minus-range}  $\Up_-(\xi)$ is defined for 
$t_\al^{-1}$
in place of  $t_\al.$ 
The {\it semisimple range}  
$\Up_*[\xi]$ is the intersection
$\Up_+[\xi]\cap \Up_-[\xi].$
The {\it zero-range}  $\Up_0[\xi]$ is introduced for
$t_\al^{0}=1,$ i.e. without $t_\al$ in  (\ref{upset}).
The boundaries are defined respectively.
Notice that  
\begin{align}
&\Up_-[\xi]= \Up_+[-\xi],\ \Up_0[\xi] = \Up_0[-\xi],
\notag\\ 
&\ddot{\Up}_*[\xi]\subset \ddot{\Up}_+[\xi]
\cup \ddot{\Up}_-[\xi].
\label{rangepm}
\end{align}

We will use the compatibility of the semisimple and
zero ranges with the right multiplications.
Namely, 
\begin{align}
&\Up_*[\hw\llb \xi\rrb]\ =\ \Up_*[\xi]\,\hw^{-1}
\hbox{\ for \ every \ } \hw\in \Up_{*}[\xi].
\label{rangem}
\end{align}
It readily results  from (\ref{tlaw}).
The same holds for $\Up_0\ $ (but, generally speaking, 
not for 
$\Up_\pm$). 

\vskip 0.2cm
{\it Semisimple stabilizer.}
The following {\it semisimple stabilizer} of $\xi$ 
\begin{align}
&\hW_*^\flat[\xi]\equal \{\hw\in \hW^\flat[\xi]\, \mid
\, \hw\in \Up_*[\xi]\}\ =\ \hW_*^\flat[-\xi].
\label{unstab}
\end{align}
will play an important role in the main theorem of this section.
It is a subgroup of $\hW^\flat[\xi].$ Indeed, if 
$\hu,\hw\in  
\hW_*^\flat[\xi],$ then 
the set $\la(\hu)\cup\hu^{-1}(\la(\hw))$
satisfies (\ref{upset}) 
and the corresponding relation for $\Up_-[\xi]$
because $\hu$ does not change $q^\xi.$
The set $\Up_-[\xi]$
contains $\la(\hw\hu)$ even if $l(\hw\hu)<l(\hw)+l(\hu).$
Apply (\ref{tlaw}) to the product of the reduced decompositions
of $\hw$ and $\hu,$ which may be nonreduced, to see it.

Moreover, this very argument gives that
\begin{align}
&\Up_*[\xi] \hw = \Up_*[\xi] 
\hbox{\ for\ every\ }
\hw\in \hW_*^\flat[\xi] \and
\label{unstabi}
\\ 
&\hW_*^\flat[\hw\llb \xi\rrb]= 
\hw\, \hW_*^\flat[\hw\xi]\, \hw^{-1} \hbox{\ for\ all\ }
\hw\in \Up_*[\xi].
\label{stabc}
\end{align}
The last property is directly connected with (\ref{rangem}).

Similarly, one can introduce the groups 
$\hW_\pm^\flat[\xi]$
and  $\hW_0^\flat[\xi].$ They preserve the corresponding
$\Up[\xi]$ acting via the right multiplication.

\vskip 0.2cm
{\it Sharp case.}
When considering spherical and self-dual
representations, the {\it sharp range} $\Up_+[-\rho_k]$ 
will be used, which is simply $\Up_+[\xi]$ as 
$\xi=-\rho_k,$
and other $\Up$ for such $\xi.$

Recall that by $c\leftrightarrow\!\rightarrow b,$  
we mean that 
$b\in B$ can be obtained from $c\in B$ by
a chain of transformations i),ii),iii) (the simple arrows)
from Proposition \ref{BSTIL}. Here $c\leftrightarrow$ indicates that 
transformations of type iii) may be  involved. The latter correspond 
to
the applications of $\pi_r\in \Pi^\flat$ and are invertible. 
Otherwise
we put  $c\rightarrow\!\rightarrow b.$ 

This usage of arrows is actually a particular case of the 
$\hu\mapsto \hw$ upon the restriction to  
$\Up_+(-\rho_k).$  
Indeed, the
latter is the set of all $\hw=\pi_b$
satisfying (\ref{upset}). If $\pi_c\in \Up_+(-\rho_k)$ and
$c\rightarrow b$ corresponds to $\pi_b=s_i\pi_c,$
then  $\pi_b\in \Up_+(-\rho_k)$ if and only if
\begin{align}
&t_\al q^{(\al\, ,\, c_- -\rho_k)}\neq 1 \for i>0,\ 
\al=u_c(\al_i),\notag\\ 
&t_\al q^{1+(\al\, ,\, 
c_- -\rho_k)}\neq 1 \for i=0,\ \al=u_c(-\vth).
\label{upsharp}
\end{align}
Using $c_\#=c-u_c^{-1}(\rho_k)=u_c^{-1}(-\rho_k+c_-)$
and $([\al, j]\, ,\, d)=j,$
\begin{align}
&t_\al q^{(\al_i\, ,\, c_\#+d)}\neq 1 \for i\ge 0. 
\label{upsha}
\end{align}
If $c\leftrightarrow b,$ i.e. $\pi_b=\pi_r\pi_c$ for 
$\pi_r\in \Pi^\flat,$
then  always $\pi_b\in \Up_+(-\rho_k).$
\vskip 0.2cm

{\bf Main Theorem.}
We turn  to a general
description of semisimple and pseudo-unitary
representations. It will be reduced to a certain combinatorial
problem which can be managed in several cases. We permanently
identify the cyclic generator 
$v\in \i_X[\xi]$ with $1\in \h_Y^{\flat}$, which can be
uniquely extended to
a $ \h_Y^{\flat}$-isomorphism 
$\i_X[\xi]\cong\h_Y^{\flat}.$

\begin{theorem}
i) Assuming that
\begin{align}
&\Up_*[\xi]\subset\Up_0[\xi] \and
\j_\xi\equal\sum \h_Y \Phi_{\hw}(q^{\xi})\, \neq\, \h_Y, 
\label{uniset}   
\end{align}
where the summation is over $\hw\in \ddot{\Up}_*[\xi]$,
the quotient  
$U=U_X^\xi\equal \h_Y^{\flat}/ \j_\xi$
is an $X$-semisimple 
$\HH^\flat$-module with a basis of $X$-eigenvectors
$\{\Phi_{\hw}(q^{\xi})\}$ for $ \hw\in \Up_*[\xi].$
Its $X$-spectrum is $\{\hw(\xi)\}$ for such $\hw.$ 

If $V$ is
an irreducible $X$-semisimple $\HH^\flat$-module 
with a cyclic
vector of weight $\xi$ then (\ref{uniset})
holds and $V$ is a quotient of $U,$ provided that
$t_{\sht}\neq 1\neq t_{\lng}.$

ii) The elements  $\Phi_{\hw}v',\ $ $S_{\hw}v',\ $ 
$G_{\hw}v'$
are well defined for all $v'$ from the $\xi$-eigenspace
$U(\xi)$ of $U$ and  $\hw\in \Up_*[\xi].$
They induce isomorphisms $U(\xi)\cong U(\hw(\xi)).$

The group $\hW_*^\flat[\xi]$ acts in $U(\xi)$ by automorphisms
via $\hw\mapsto S_{\hw}.$ A quotient $V$ of $U$ is an
irreducible $\HH^\flat$-module if
and only if  $V(\xi)$ is an irreducible  
$\hW_*^\flat[\xi]$-module.
Moreover, provided (\ref{uniset}),
an arbitrary irreducible  representation $V'$ of this 
group can be uniquely extended to 
an irreducible $X$-semisimple  $\HH^\flat$-quotient $V$ of  
$\i_X[\xi]$ 
such that $V'=V(\xi).$
Here $G$ can be used instead of $S.$ 

The above statements hold for  
 $\i_Y[\xi]$ and its semisimple quotients
when $\Phi_{\hw}(q^{\xi})$
are replaced by $\tau_+\bigl(\Phi_{\hw}(q^{-\xi})\bigr)$ everywhere.

iii) Under the same constraint (\ref{uniset}), 
the function 
\begin{align}
&\mu_\bullet(\hw)\equal
\prod_{[\al,\nu_\al j]\in \la(\hw)}
\Bigl(
\frac{
t_\al^{-1/2}-q_\al^jt_\al^{1/2} X_\al(q^{\xi})}{
t_\al^{1/2}-q_\al^jt_\al^{-1/2} X_\al(q^{\xi})
}
\Bigr)
\label{muuni}
\end{align}
has neither poles nor zeros at $\hw\in \Up_*[\xi].$ Moreover,
\begin{align}
&\mu_\bullet(\hu)\mu_\bullet(\hw)
\ =\ \mu_\bullet(\hu\hw)  \hbox{\ whenever\ }
\hw\in \hW_*^\flat[\xi]
\label{mumono}
\end{align}
for $\hu,\hw\in\Up_*[\xi],$ 

A $\HH^\flat$-quotient $V$ 
of $U$ is pseudo-unitary if and only if
the eigenspace $V(\xi)$ is a pseudo-unitary
$\hW_*^\flat[\xi]$-module. The latter means that
the form $(\,,\,)$ is nondegenerate on $V(\xi)$ and
\begin{align}
&(S^\xi_{\hw})^\star\, S^\xi_{\hw}=\mu_\bullet(\hw)
\for \hw\in \hW_*^\flat[\xi],\ 
S^\xi_{\hw}\equal  S_{\hw}(q^\xi).
\label{muext}
\end{align}

iv) A pseudo-hermitian form on  $e',e''\in V(\xi),$
can be uniquely extended to 
an pseudo-hermitian invariant form on $V:$
\begin{align}
&(S_{\hu}e'\, ,\, S_{\hw}e'')\ =\ 
\de'_{\hu,\hw}\,\mu_\bullet(\hu)\,
(e'\, ,\, S_{\hu^{-1}\hw}e''),\where 
\label{muformx}
\\ 
\hu,\hv\in \Up_*[\xi]&,\ 
\de'_{\hu,\hw}=1 \for \hu^{-1}\hw\in \hW_*^\flat[\xi] \and
=0 \hbox{\ otherwise}.
\notag \end{align}
If $(e',e')\neq 0$ as $V(\xi)\ni e'\neq 0$ then $V$ is 
$X$-unitary. 

An arbitrary irreducible pseudo-unitary ($X$-unitary) 
$\HH^\flat$-module with a cyclic vector
of weight $\xi$ can be obtained by this construction as 
$t_{\sht}\neq 1\neq t_{\lng}.$

In the $Y$-case,  
$$\tG^\xi_{\hw}\equal
\tau_+\bigl(G_{\hw}(q^{-\xi})\bigr)$$
is taken instead of 
$S_{\hw}(q^\xi);$ these operators send $V(\xi)$ to 
$V(\hw\llb \xi\rrb).$
The extension from a pseudo-hermitian form on $V(\xi)$ is 
\begin{align}
&( \tG^\xi_{\hu} e'\, ,\, \tG^\xi_{\hw} e'')\ =\ 
\de'_{\hu,\hw}\,\mu_\bullet(\hu)^{-1}\,
(e'\, ,\,  \tG^\xi_{\hu^{-1}\hw} e'').
\label{muformy}
\end{align}
\label{PUNIT}
\end{theorem}

{\it Proof.}
We set $e_{\hw}= \Phi_{\hw}(1)$ for $\hw\in \hW^\flat$ 
provided that   $\Phi_{\hw}(1)$ is well defined. It
is an $X$-eigenvector of weight $\hw(\xi).$
Recall that we permanently identify $v$ and $1.$
Then $\h_Y e_{\hw}$ is a $\HH^\flat$-submodule because
$\HH^\flat= \h_Y\cdot \Q_{q,t}[X_b].$ 

Assuming (\ref{uniset}), $\j_\xi$ and $U$ 
are  $\HH^\flat$-modules. Moreover,
all $e_{\hw}$ are invertible in $\h_Y^\flat$ and 
their images linearly generate $U$ as 
$\hw\in \Up_*[\xi].$   
Indeed, if $\Phi_i e_{\hw}$ loses the invertibility 
then $s_i{\hw}$ leaves
$\Up_*[\xi],$  $l(s_i\hw)=l(\hw)+1,$ and
$s_i\hw\in \ddot{\Up}_* [\xi].$ By (\ref{uniset}), 
$\ddot{\Up}_* [\xi]\subset \Up_0[\xi].$ We conclude that 
the element $T_i(e_{\hw})$ is
proportional to $e_{\hw}$ in this case. 

This means that each $T_i(e_{\hw})$ is a linear 
combination of $e_{\hu}$
for any $n\ge i\ge 0,$ where  $\hw,\hu\in \Up_* [\xi].$
Of course it holds for $\pi_r\in \Pi^\flat$ too.
Thus $\{e_{\hw'}\}$ linearly generate $U.$
Moreover, given $\hw\in \Up_*[\xi],$ 
the map $v'\mapsto  \Phi_{\hw}v'$
induces an isomorphism $U(\xi)\cong U(\hw(\xi)).$

Summarizing, (\ref{uniset}) implies that
$\{e_{\hw\hu}\}$ form a basis of $U$ when   \\ 
 {(a)} $\hw$ are representatives of all classes
$\Up_*[\xi]/\hW^\flat[\xi],$  \\ 
 {(b)} $\hu\in\hW^\flat[\xi]$  and $\{e_{\hu}\}$ form a basis
of $U(\xi).$

Let us discuss the usage of $S$ and $G.$ 
The elements  $S_{\hw}(1)\in \h_Y^\flat$ are well 
defined and for
$\hw\in \Up_+[\xi].$ The sum 
$\j_\xi^+\equal\sum \h_Y S_{\hw}\ ,$
$\hw\in \ddot{\Up}_+[\xi],$ is a  
$\HH^\flat$-submodule of $\i_X[\xi].$
Finally,  $S_{\hw}(1)$ are invertible
for $\hw\in  \Up_*[\xi],$  which results 
directly from the 
definition of
the latter set. 
The same holds for  $\Up_-[\xi]$
with $S$ being replaced by $G.$ Turning to $\i_Y[\xi],$
$G$ and $S$ serve respectively $\Up_+[\xi]$ and 
$\Up_-[\xi].$  

\vfil
\vskip 0.2cm
Now let us check that the existence of irreducible 
$X$-semisimple quotients
$V\neq \{0\}$ of $\i_X[\xi]$ implies (\ref{uniset}),
and all such $V$ are quotients of $U.$ 
We suppose that $t_{\sht}\neq 1\neq t_{\lng}.$
Actually we will impose  
$t_{\sht}\neq \pm 1\neq t_{\lng}$ 
in the argument below, leaving the consideration of
$t=-1$ to the reader. This allows separating 
$T_i-t_i^{1/2}$
and $T_i+t_i^{-1}.$ When $t_i=-1$ the $T_i$ is not semisimple, 
but the theorem still holds in this case.

The images $e_{\hw}$ of the vectors 
$\Phi_{\hw}(1)$ in $V$ are
nonzero as $\hw\in \Up_*[\xi].$ 
If  $\Up_0[\xi]\not\supset\Up_*[\xi],$  we can
pick $e'=e_{\hw}$ such that
$S_i(e')=e'\ $ and  $X_{\al_i}(e')=e'.$  This gives
that $X_{\al_i}$ is
not semisimple in the two-dimensional module 
generated by $e'$
and $e''=T_i(e').$ The dimension is two because 
the characters
of $\h_X^i\equal$ $\lan T_i,X_{\al_i} \ran$ 
send $X_{\al_i}\mapsto t_i^{\pm}\neq 1.$ So it cannot be one.

Now let
$\hu$ belong to $\Up_*[\xi]\ $ and  
$\hw=s_i\hu\in \ddot{\Up}_+[\xi].$
The case  $\hw=s_i\hu\in \ddot{\Up}_-[\xi]$
is completely analogous. We will skip it.
Then $l(\hw)=l(\hu)+1.$ Both  $e'=\Phi_{\hu}(e)$ and 
$e''=\Phi_i(e')=$ $(T_i-t_i^{1/2})(e')$ are $X$-eigenvectors.
Here by $e,$ we mean the image of $1$ in $V.$
Let us prove that $e''=0.$

Supposing that $e''\neq 0,$ it 
cannot be proportional to $e'.$ 
Indeed, this may only happen if $(T_i+t_i^{-1/2})(e')=0,$ 
which is impossible since $\hu\in \Up_*[\xi]$
$\subset \Up_-[\xi].$

We note that if $V$ is $X$-unitary, then restricting
the $X$-definite form to the linear span of $e',e''$
readily leads to a contradiction with the relation  
$T_i^\star=T_i^{-1}.$ This argument doesn't require
the irreducibility of $V,$ but involves the $X$-unitary
structure. 

If $e''\neq 0$ then it generates $V$ as a 
$\HH^\flat$-module 
and as a $\h_Y^\flat$-module thanks to the irreducibility. 
This means that 
\begin{align}
&\h_Y^\flat(T_i-t_i^{1/2})\Phi_{\hw}(1) + J_V=\h_Y^\flat,
\where V=\h_Y^\flat/J_V,
\notag \end{align}
for a proper ideal $J_V\subset \h_Y^\flat$ which is  a 
$\HH^\flat$-submodule of $\i_X[\xi].$
Hence, $J_V$ contains 
$ \h_Y^\flat(T_i+t_i^{-1/2})\Phi_{\hw}(1)$
and $e''$ has to be  proportional $e'.$ This is 
impossible because
$X_{\al_i}(e'')=t_i^{-1}$ whereas 
$X_{\al_i}(e')\neq t_i^{\pm 1}.$ 

Summarizing, we almost completed parts i),ii) from
the theorem. The remaining
claims are  \\ 
 {a)} the equivalence
of the irreducibility of the $\HH^\flat$-module $V$
and the $\hW_*^\flat[\xi]$-module $V(\xi),$  \\ 
 {b)} the possibility of taking an arbitrary irreducible 
representation of $\hW_*^\flat[\xi]$ as $V(\xi).$

Applying the intertwining operators,
both claims readily follow from the statements which have been  
checked.

Let us turn  to iii),iv).
We take $e',e''\in V(\xi)$ and denote  
$e'_{\hw}=S_{\hw}e'$ and 
 $e''_{\hw}=S_{\hw}e'',$ where $\hw\in \Up_*[\xi],$ for 
a certain pseudo-unitary irreducible quotient 
$V$ of $\i_X^\flat.$
They are nonzero $X$-eigenvectors. 
If the eigenvalues of $e'_{\hu}$
and $e''_{\hw}$ are different, i.e. 
$\hu^{-1}\hw\not\in \hW_*^\flat[\xi],$
then these vectors are orthogonal with respect to the 
invariant form of $V.$

Applying (\ref{gstarg}),
\begin{align}
&(e'_{\hw}\, ,\, e''_{\hw})\ =\ 
(S_{\hw}e'\, ,\, S_{\hw}e'')\ =
\ ( S_{\hw}^\star\, S_{\hw}e'\, ,\, e'')\ =
\label{eemu} 
\\ 
&(e',e'')\,\prod_{\al,j} \Bigl(
\frac{
t_\al^{-1/2}-q_\al^j t_\al^{1/2} X_\al(q^\xi)}{
t_\al^{1/2}-q_\al^j t_\al^{-1/2} X_\al(q^\xi) 
}
\Bigr),\ [\al,\nu_\al j]\in \la(\hw).
\notag \end{align}
The latter product is nothing else but $\mu_\bullet(\hw)$ from 
(\ref{muuni}).
Thus the extension of $(\ ,\ )$ from $V(\xi)$ to $V$ is unique
if it exists. Explicitly:
\begin{align}
&(S_{\hu}e'\, ,\, S_{\hw}e'') = 
\de'_{\hu,\hw}\,\mu_\bullet(\hu)\,
(e'\, ,\, S_{\hu^{-1}\hw}e''),
\where e',e''\in V(\xi),
\label{muforx}
\\ 
&\hu,\hv\in \Up_*[\xi],\ \
 \de'_{\hu,\hw}=1 \for \hu^{-1}\hw\in \hW_*^\flat[\xi] \and
0 \hbox{\ otherwise}.
\notag \end{align}

To check that this formula  extends correctly an
arbitrary given definite form on  $V(\xi),$ it suffices 
to compare 
$$
(S_{\hw}e''\, ,\, S_{\hu}e')^* \and 
(S_{\hu\hv}(S_{\hv}^{-1}e')\, ,\, S_{\hw}e'')
$$
and to show that they coincide whenever $\hv\in  
\hW_*^\flat[\xi].$
Using (\ref{muforx}), we need the relations
$$
\mu_\bullet(\hw^{-1}\hu)\mu_\bullet(\hw)=
\mu_\bullet(\hu),\
\mu_\bullet(\hu\hv)\mu_\bullet(\hv^{-1})=\mu_\bullet(\hu),
$$
which are nothing else but (\ref{mumono}). Claim iv) is verified.

Relation (\ref{mumono}) obviously holds when 
$l(\hu\hw)=
l(\hu)+l(\hw),$ because the $\mu_\bullet$-product for
$\hu\hw$ is taken over
$\la(\hu\hw),$ which is a union of $\la(\hw)$ and 
the set $\hw^{-1}(\la(\hu)).$ The $\mu_\bullet$-product
over the latter set 
coincides with the product over  $\la(\hu)).$   

Thanks to (\ref{tlaw}), we can omit the constraint
$l(\hu\hw)=l(\hu)+l(\hw)$ in this argument. Indeed,
the extra positive affine roots 
$\tal=[\al,\nu_\al j]\in \tla^+(\hu\hw)$
appear in $\tla(\hu\hw)$ together  with 
$-\tla\in -\tla^+(\hu\hw).$ Here we take the product
of the reduced decompositions of $\hu$ and $\hw,$ which
may be nonreduced, to construct $\tla(\hu\hw).$ However  
the product of the  $\mu_\bullet$-factors corresponding 
to $\tal=[\al,\nu_\al j]$ and $-\tal$ is 
\begin{align}
&\Bigl(
\frac{
t_\al^{-1/2}-q_\al^j t_\al^{1/2} X_\al(q^{\xi})}{
t_\al^{1/2}-q_\al^j t_\al^{-1/2} X_\al(q^{\xi})
}
\Bigr)
\Bigl(
\frac{
t_\al^{-1/2}-q_\al^{-j} t_\al^{1/2} 
X_\al^{-1}(q^{\xi})}{
t_\al^{1/2}-q_\al^{-j} t_\al^{-1/2} 
X_\al^{-1}(q^{\xi})
}
\Bigr)
\ =\ 1.
\notag \end{align}
This concludes iii) and the proof of the theorem.
\qed

Thanks to Proposition \ref{PIOMG}, there is a reasonably
simple and explicit combinatorial reformulation of the first
part of the condition (\ref{uniset}). 
However it doesn't guarantee
the second part.

\begin{proposition} The relation
$\Up_*[\xi]\subset\Up_0[\xi]$
is equivalent to the following property of $\xi.$
Let $\ddot{\Up}$
$=\{\tbe=[\be,\nu_\be j]\}$ be the set
of all affine positive roots such that
$t_\be^{\pm 1} X_{\tbe}(q^\xi)=1.$ Then an arbitrary 
positive affine $\tga=[\ga,\nu_\ga i]$ such that
$X_{\tga}(q^\xi)=1$ is a linear combination of some 
$\tbe\in \ddot{\Up}$
with positive rational (not always integral) coefficients. 
Using the closure $\bar{\CC_a}$ of the affine Weil chamber 
$\CC_a$ from Section 1, 
\begin{align}
&\cup_{\hw}\hw\llb -\bar{\CC_a} \rrb
= \{z\in \R^n,\ (z,\be)+j<0\}
 \for \hw\in \Up_*[\xi],\ \tbe\in \ddot{\Up}.
\label{semibase} 
\end{align}
\label{UNISETG}
\qed
\end{proposition}

\vskip 0.2cm
{\it Generic cyclic representations.}
We will start the discussion of the theorem with the 
irreducible $\i_X[\xi].$ Following [C1], Theorem 6.1,
let us assume that $\hW_*^\flat[\xi]$
is conjugated in $\hW$ to
the span of a proper subset of $\{s_0,s_1,\ldots,s_n\},$
maybe $\emptyset$ but smaller than the whole set of 
affine simple reflections. It is possible only if $q$ 
is not a root of unity, because otherwise  $\hW_*[\xi]$ 
is infinite. By the way, this constraint results in the 
implication 
\begin{align}
&\Up_*[\xi]\subset\Up_0[\xi] \Rightarrow
\j_\xi=\sum \h_Y \Phi_{\hw}\, \neq\, \h_Y, \where
\hw\in \ddot{\Up}_*[\xi].
\notag \end{align}
We do not prove/use this claim in the paper. 

Apart from 
roots of unity,
$\i_X[\xi]$ is irreducible if and only if 
$\Up_+[\xi]=\hW^\flat=\Up_-[\xi].$ Thus its semisimplicity
is equivalent to the triviality of the stabilizer 
$\hW_*^\flat[\xi].$  

The theorem can be applied to a general description
of the finite-dimen\-signal semi\-simple representations.
We continue to suppose that $q$ is generic.

\begin{proposition}
i) Under condition (\ref{uniset}), the induced module 
$\i_X[\xi]$
has a nonzero finite-dimensional $X$-semi\-simple 
quotient if
and only if there exists a set of roots 
$T\subset R$
such that for any $b\in B$ there exist $\be\in T$
with $(b,\be)>0,$ i.e. $\be$ do not belong to any 
halfplane
in $\R^n,$  
and for every $\be\in T,$
\begin{align}
(\xi,\be^\vee)&+k_\be \in -1-\Z_+
\hbox{\ or\ } (\xi,\be^\vee)-k_\be \in -1-\Z_+ \for \be<0,\notag
\\ 
&(\xi,\be^\vee)+k_\be \in -\Z_+
\hbox{\ or\ } (\xi,\be^\vee)-k_\be \in -\Z_+ \for \be>0.
\label{semifi}
\end{align}
\label{SEMIFI}
\end{proposition}
{\it Proof}. Let us
assume that for all $\al\in R$ such that
$(b,\al)> 0,$ either
$(\xi,\al^\vee)\pm k_\al \not\in -1-\Z_+$
for both signs of $k_\al$ as $\al<0$ or the same inequalities
hold with $-\Z_+$ as $\al>0.$ 
Then $\la(lb)$ for $l\in \Z_+$
consists of $\tal=[\al,\nu_\be j]$ for such $\al$ with integers
$j$ satisfying inequalities $0< j\le l(b,\al^\vee)$ as 
$\al<0$
and with $0\le j$ otherwise. We take $b$ from 
$\Up_*[\xi],$ so
all these $\tal$  satisfy (\ref{upset}) and its 
counterpart
with $t_\al^{-1}$ instead of $t_\al.$
Since $l$ is arbitrary positive, we get that $\Up_*[\xi]$ is
infinite if the relations (\ref{semifi}) don't hold. 

Let us check that this condition guarantees that
$\Up_*[\xi]$ is finite. 

In the first place, any given infinite
sequence of pairwise distinct $b^i\in B$
is a (finite) union of subsequences $\{\tb^{i}\}$ 
such that there exist $\be_\pm\in T$ satisfying 
$\pm(\tb^{i},\be_\pm)>0.$
Taking one such infinite subsequence, we may assume that
there exist $\be_\pm\in T$ such that $\pm(b^i,\be_\pm)>0.$
If the latter scalar products are bounded, we
switch to a new sequence $b^i$ formed by proper 
differences of the
old $\{b^i\}$ ensuring that all $(b^i,\be_\pm)=0.$ 
Then we find
the next pair ${\be_\pm}'$ such that 
$\pm(b^i,{\be_\pm}')>0,$ form the next
sequence of differences, and so on until all $T$ is exausted.

We conclude that any sequence 
$\{b^i\}$ is a union of subsequencies
$\{\tb^i\}$ such that $|(\tb^i,\be)|\to \infty$ as 
$i\to \infty$
for certain $\be\in T.$  
Following the proof
of Proposition \ref{LAPOSI} (see (\ref{lawb}),
(\ref{lawbb})),
given $w\in W,$
the sets $\la(w\tb^i)$ contain the roots
$[\be,j\nu_\be]$  where $\max\{j\}\to \infty$ as 
$i\to \infty.$ 

Note that the proof gets more transparent via
the geometric interpretation from
Proposition \ref{PIOMG}. 
\qed

The generalized Macdonald identities correspond
to the set 
$T=\{\al_1,\ldots,\al_n,-\vth\}$ taken together
with $\xi=-\rho_k.$
The sign of $k_\be$ in the relations (\ref{semifi}) from the Proposition
is  plus for $\vth$ and minus otherwise. 
See the end of Section 8.

The constraint (\ref{uniset}) will be fulfilled, 
for instance, if 
$(\xi,\al^\vee)\not\in \Z$ for all $\al\in R.$
The simplest example is as follows.
Letting $k_{\lng}=k=k_{\sht},$  $\xi=-k\rho,$  
we take $k=-m/h$ for the
Coxeter number $h$  provided that $(m,h)=1.$ 
This representation is spherical. 
In the $A_n$ case, one can follow [C9] to describe
all irreducible representations under consideration
(for generic $q$).

\vskip 0.2cm
{\bf $GL_n$ and other applications.}
Let us give a description of the $X$-semisimple and general
representations for the double Hecke algebra 
associated to $GL_n$ :
\begin{align}
&\HH_n\equal\langle X_1^{\pm 1},\cdots,X_n^{\pm 1}, 
\pi, 
T_1,\cdots,T_{n-1}\rangle,
\where
\label{HHgln}
\end{align}

(i)\ \,$\pi X_i=X_{i+1}\pi$ ($i=1,\ldots n-1$) and 
$\pi^n X_i=q^{-1} X_i \pi^n$ ($i=1,\ldots,n$),

(ii)\,\, $\pi T_i=T_{i+1}\pi$ ( $i=1,\ldots n-2$) and 
$\pi^n T_i=T_i \pi^n$ ($i=1,\ldots,n-1$),

(iii)$X_iX_j=X_jX_i$ ($1\le i,j\le n$), 
$\ T_iX_iT_i=X_{i+1},\ i<n,$
 
(iv)\,\, $T_i T_{i+1} T_i=T_{i+1} T_{i} T_{i+1},\ $ 
$T_i T_j=T_j T_i$ as $i<n-1\ge j,\ |i-j|>1,$

(v)\ \,$(T_i-t^{1/2})(T_i+t^{-1/2})=0$ for 
$T_1,\ldots,T_{n-1},\
\ t=q^k.$ 

Switching here from $\C[X_i^{\pm 1}]$ to the
subalgebra $\C[X_iX_j^{-1}]$ of Laurent polynomials
of degree zero, we come to the double Hecke 
algebra $\HH_n^o$ of type $A_{n-1},$ i.e. that for $SL_n.$
To be more exact, we must also set $\pi^n=1,$ which is possible 
because $\pi^n$ becomes central thanks to (i,ii), and in this
definition the root lattice $Q$ is used for 
$X$ in place of the 
whole weight lattice $P.$

Later on we will assume that $\pi$ is invertible
and identify the representations under the automorphisms
$\pi\mapsto c\pi$ of $\HH_n$ for $c\in\C^*.$ 

The following construction is a "cylindrical" counterpart of that
from Section 3 from [C9] (see also [Na]) and generalizes 
the corresponding
affine classification due to Bernstein and Zelevinsky.
It is a direct application of the Main Theorem.
We will publish the details elsewhere. 

Let $\De=\{\de^1,\cdots,\de^m\}$ be a set of 
skew Young diagrams without empty rows of the total 
degree $n$\,
(i.e. with $n$ boxes in their union),\ 
$C=\{c_1,\cdots,c_m\}$ a set of complex
numbers such that $q^{c_a - c_b}\not\in t^{\Z}q^{\Z}$
as $a\neq b.$
 
We number the boxes
of these diagrams in the {\it inverse order}, i.e.
from the last box of the last row of
$\de^m$ through the first
box of the first row of $\de^1,$ and set
$$\xi_l\ =\ c_{p_l}+k(i_l-j_l),\ l=1,\ldots, n,$$ 
where $\de^{p_l}$ contains box $l$ 
($p_1=m,\ldots, p_n=1$), and
$i_l,j_l$ are the row and column numbers of the 
$l$-th box in
the corresponding diagram. 
We say that $\xi$ is associated with the {\it pair} 
$\{\De,C\}.$ 

Concerning skew diagrams, if box $l$ belongs to the 
$i$-th row of 
the $p$-th diagram (so $i_l=s, p_l=p$) then
$m_i< j_l\le n_i$ for the "endpoints"
$m_i\le n_i$ of the rows 
(if they coincide the row is empty) 
satisfying the inequalities 
\begin{align}
&m_i\ge m_{i+1},\ n_i\ge n_{i+1},\ 1\le i < 
(\hbox{\ number \ of
\ rows \ in \ } \de_p).
\label{mini}
\end{align}

\begin{theorem}
Let $q^i t^j\neq 1$ for any $(i,j)\in \Z^2\setminus 
(0,0).$   
Then condition (\ref{uniset}) is satisfied
if $\xi$ is associated with a $\De$-set of skew diagrams and
the corresponding $C$-set of complex numbers. The quotient
$U=U_X^\xi$ from the Main Theorem is irreducible and 
$X$-unitary for such $\xi.$ Its $X$-spectrum is simple.
An arbitrary irreducible $X$-semisimple $\HH_n$-module can
be constructed in this way for a proper pair $\{\De,C\}.$
Isomorphic modules correspond to the pairs which can be obtained
from each other by a permutation of the skew diagrams in
$\De$ together with the corresponding numbers in $C,$ 
and by 
adding arbitrary integers to the $c$-numbers. 
\label{HHGLNG}
\end{theorem}
{\it Proof.} See [C1].
\qed

Note that the weight $\xi$ from the theorem
can be used to construct 
an irreducible representation of the affine Hecke algebra
$$\h_n\equal \langle X_1^{\pm 1},\ldots, X_n^{\pm 1},
T_1,\ldots,T_{n-1}\rangle.$$
Cf. [C1] and [C9]. Inducing this module from $\h_n$ to $\HH_n,$ we 
get $U_X^\xi.$ The same holds for the
subalgebra $\h_n^o$ of $\h_n$ with $\C[X_iX_j^{-1}]$
instead of $\C[X_i^{\pm 1}]$ and for the natural restriction
of $\xi$ to  $\C[X_iX_j^{-1}].$

\vskip 0.2cm
{\it Periodic diagrams.}
Assuming that $q$ is not a root of unity, 
let $t^r=q^s$ for integral $r,s$ where $s>0$ and 
$(r,s)=1.$
Using $t=q^k,$ we get $k=s/r.$ 
Let $r\in \Z_+.$ An infinite skew diagram $\de$ without
empty rows is called
$r$-{\it periodic} of degree $n$ if 
it is invariant with respect to translation by a vector
$w=(v, v-r)$ in the $(i,j)$-plane for a positive integer 
$v$ such that the subdiagram of $\de$ with the rows 
$1\le i \le v$
contains $n$ boxes. 

Here we need to fix the first row in the diagram,
although the isomorphism classes of the resulting 
representations will not depend on such choice. The above 
subdiagram is called the {\it basic subdiagram}. It is 
of course skew.

Any skew subdiagram which is the fundamental
domian for the action of $\Z w$ on $\de$ (so it must contain
$n$ boxes) is called a {\it fundamental subdiagram}. 
It may have empty rows.
An example is the basic subdiagram. However there are infinitely 
many non-basic ones
unless $\de$ is an infinite column.  

Similarly, $\de$ is called $r$-periodic for negative 
$r$ if it is
invariant with respect to $w=(v+r,v)$ and has a fundamental 
subdiagram with $n$ boxes. The basic diagram is formed by 
the first $v$ consecutive columns in this case. 

\vfil
\vskip 0.2cm
{\it Partitions.}
A periodic skew diagram is naturally
a portion of the $\{i,j\}$-plane
between its {\it upper-left boundary} 
(the boundaries of a box are
its four sides) and its {\it lower-right boundary}. 
The boundaries are by definition continuous
curves made from segments where either $i\in \Z$ 
or $j\in \Z.$ 
They satisfy the skew condition. We will call them
{\it skew paths}.

By an $s$-{\it partition} of a skew diagram $\de$
of degree $n$ 
we mean its representation as a disjoint union of $s$ 
skew subdiagrams $\{\de_0,\cdots,\de_{s-1}\}$ which are 
the portions of $\de$ between consecutive paths 
from a certain system of skew paths without crossings.
The latter means that different paths from such a system may
have coinciding subpaths but do not
cross each other. We number the paths upwards.
We will call such partitions {\it increasing}.

Note that some of the subdiagrams in the partition
can be empty and they may have empty rows.

The subdiagrams in a given partition 
are {\it partially ordered} naturally. Namely,
$\de'>\de''$ if at least one box of $\de'$ is in the
same column and higher than a box from $\de''$ or
is in the same row and to the left of a box from $\de''.$

Using this definition,
increasing skew partitions can be introduced without 
using the paths as follows.
First, it is not allowed that $\de_a> \de_b$
and at the same time $\de_a < \de_b$ for different $a,b.$
Second, $a>b$ if $\de_a > \de_b.$ 

For instance, the representation of a skew diagram
$\de$ of degree $n$ as a union 
of all its boxes counted from the
last box in the last row through the first box in the 
first
row (i.e. with respect to the inverse order) is
an increasing $n$-partition. However there are different
increasing orderings of these boxes unless $\de$ is a column
or a row. 

We will need the partitions of the fundamental subdiagrams.
They appear naturally when decomposing semisimple
representations of $\HH_n$ under the
action of $\h_n.$

\vfil
\vskip 0.2cm
{\it New pairs.}
Extending these definitions to $\De,$ we call it 
$r$-periodic of 
degree $n$ if all $\de^p$ are $r$-periodic and the total
degree of $\{\de^p\}$ is $n.$ 
The weights $\xi$ associated with a {\it new pair}
$\{\De,C\}$ are those for {\it any} 
choice of fundamental
subdiagrams of $\{\tde^p\}$ calculated using the formulas above.

Here we need {\it recalculating}  the $c$-numbers as follows.
If $c$ is assigned
to periodic $\de$ then the corresponding $c$-number of a 
fundamental subdiagram $\tde$
will be $c+(i-j)k$ for the coordinates $(i,j)$ 
(inside $\de$)
of the first box in 
the first nonempty row of the subdiagram $\tde$. 
Upon such a recalculation,
we call the set of fundamental subdiagrams a
{\it fundamental subpair}.

Here we need to choose the first line,
i.e. the one with $i=1.$ It is was also necessary in the
definition of the basic subdiagram.
Changing the position of the first row
results in recalculating the corresponding $c$-number.
If the first box (with $i=1,j=1$)
remains in the same principal diagonal $i=j+$const, 
then there will be no $c$-change. 

The constraint for the $c$-numbers of the new (periodic) 
pairs remains the same as for the old pairs.
Thanks to the condition $t^r=q^s,$   
it reads:
$$c_a-c_b\ \not\in\ (1/r)\Z \hbox{\ \ for\ all\ } 
1\le a,b\le m.
$$ 

An $s$-{\it partition} of $\{\De,C\}$ is a collection 
$\{\tde_a^p,\ 0\le a \le s-1,\ 1\le p\le m\}$
of increasing skew partitions $\tde^p = \cup_a \tde_a^p$
of some fundamental subdiagrams $\tde^p$ of $\de^p.$ 
Here, first, the $c$ numbers of $\tde_a^p$ are 
recalculated as the ones for $\de^p.$ Second,
the integer $a$ is added to the 
resulting $c$-number for every $\tde_a^p.$  

For the representation of a skew diagram
$\de$ of degree $n$ as a union 
of all its boxes counted from the
last box  through the first one, considered above,
if $c$ is assigned to $\de,$
then the $c$-number of the $a$-th
box, denoted $\de_a,$ 
is $c+(i_a-j_a)k+a,$ where $0\le a\le n-1.$ 

\vskip 0.2cm
{\it Equivalence.}
New pairs $\{\De,C\}$ and $\{\De',C'\}$ are called
equivalent if they can be obtained from each other by a permutation
of the components and the corresponding $c$-numbers
combined with adding arbitrary integers to the $c$-numbers. 
So only $c_i\,\hbox{mod}\, \Z$ matter in the (new) equivalence 
classes. Also we may shift a diagram $\de$ 
in the $(i,j)$-plane
by $(i_o,j_o)$ replacing at the same time 
$c$ by $c-k(i_o-j_o).$ 

Given a new pair,
the fundamental subpairs and their $s$-partitions 
are called equivalent if they coincide combinatorially
inside $\De.$ One checks that two partitions are 
equivalent
if and only if the corresponding pairs (the diagrams and
their $c$-numbers) coincide up to a permutation of the 
components together with the $c$-numbers. 

Given a subpair or its partition, 
treated as a set of skew diagrams with the $c$-numbers 
assigned by the above construction, 
one naturally constructs its weight and
then defines the corresponding
irreducible representations of $\h_n$ and/or $\h_n^o.$ 

\vskip 0.2cm
{\it $SL$-Equivalence.}
The subpairs 
will be called $SL-${\it equi\-va\-lent} if the
corresponding weights can be obtained from 
each other by adding a common constant.
Respectively,
their partitions are $SL$-equivalent if they
may be obtained from each other by 
proper permuting the $a$-components.

The latter permutation may be necessary when 
we add a constant which is an integer.
Indeed, we add $a$ to the $c$-number (and
the weight) of $\tde_a^p$ always
assuming that $0\le a\le s-1.$ So we may need the 
reduction modulo $s$ here, 
which leads to a certain geometric restructuring
of the diagram due to the relation $s=kr.$
After such restructuring, a permutation
of the components $\tde_a^p$ of the 
resulting $s$-partition may be needed, 
because we consider only increasing partitions. 

The corresponding irreducible representations 
of $\h_n^o$ remain in the same isomorphism 
classes under the $SL$-equivalence.

%\vfil
\vskip 0.2cm
\begin{theorem}
i) Given $r,s,n,$ irreducible $X$-semisimple representations 
of $\HH_n$ up to multiplication 
of $\pi$ by a nonzero constant 
are in one-to-one correspondence with the 
equivalence classes of new pairs $\{\De,C\}$ consisting 
of $r$-periodic sets $\De$ of degree $n$ and the 
$c$-numbers. The relations 
(\ref{uniset}) hold for the weight $\xi$ associated with 
$\{\De,C\},$ and the corresponding representation is 
$U=U_X^\xi$ from the Main Theorem. It is $X$-unitary and
its $X$-spectrum is simple. 
\vfil
ii) Let $s=1.$ Upon restriction to the affine Hecke algebra 
$\h_n,$ the above module $U=U_X^\xi$ associated to
$\{\De,C\}$ is a direct sum of the 
$\h_n$-modules associated to
the fundamental subpairs, which are, we recall, the
fundamental skew subdiagrams of $\{\de^p, 1\le p\le m\}$ 
supplied with the corresponding (recalculated) $c$-numbers. 
Each of these representations appears
exactly once in the decomposition.
\vfil
iii) For arbitrary $s>0,$ the module $U$ is isomorphic 
to a
direct sum of the $\h_n$-modules associated with the 
increasing $s$-partitions of the fundamental subdiagrams
$\{\tde^p\},$ where, recall, the integer $a$ is 
added to the $c$-number of the partition component 
$\tde_a^p.$  We treat such partitions as  sets of 
$ms$ skew diagrams equipped 
with the $c$-numbers and define the $\h_n$-modules
correspondingly. The diagrams $\tde_a^p$
can be empty for some $a$; 
they do not contribute to the decomposition.

iv) The previous claim holds for the algebra $\HH_n^o$ 
and the irreducible module $U=U_X^\xi$ defined for 
this algebra in terms of the restriction of
the above weight $\xi$ to $\C[X_iX_j^{-1}]$ (this $U$ 
is not a restriction of the above $U$  to $\HH_n^o$). 
The following change is necessary because 
the translations $c_i\mapsto c_i+c$ for  $1\le i\le m,$
and $c\in \C$ doesn't change the equivalence classes. 
Namely, the summation in the $\h_n^o$-decomposition 
of $U$ has to be reduced to 
$SL$-nonequivalent fundamental subpairs and  
nonequivalent $s$-partitions.  

v) The $\HH_n^o$-module $U_X^\xi$ is   
finite-dimensional
if and only if  $\De=\{\de\}$ and $\de$ is either
the infinite column as $r>0$ or the infinite row as 
$r<0.$
In either case, $r=\pm n$ and $s$ is an arbitrary natural number 
relatively prime to $n.$ When $s=1$ the dimension of this 
representation is $1.$ Generally, it equals $s^{n-1}.$
\label{HHGLN}
\qed 
\end{theorem}

Part v) readily follows from ii) because the number
of nonequivalent fundamental subdiagrams of a periodic
diagram is infinite unless it is a column (infinite) 
or a row. Let us describe the partitions
in this case and apply it to the dimension formula.  

\vskip 0.2cm
An $s$-partition of the $n$-column or $n$-row is 
a decomposition 
$$n=n_0+\cdots+n_{s-1} \hbox{\ for\ } n_a\in \Z_+.
$$
The numbers $n_a$ are identified with the consecutive 
segments, maybe empty, of the column (or the row), counted 
from the last box through the first one. Recall that
the $c$-number of the $a$-th segment (and the 
corresponding weight) are increased by $a,$ by the
definition. The dimension of the corresponding irreducible
representation of $\h_n^o$ is $n!/(n_0!\cdots n_{s-1}!).$

Since we may add a common integer modulo $s$ to the $c$-numbers
without changing the $SL$-equivalence class, there will
be exactly $s$ different decompositions of $n$ in every 
$SL$-equivalence class of the $s$-partitions. 
It readily gives the dimension formula $s^{n-1}$ from v).
A direct calculation (without the $\h_n^o$-decomposition)
is not difficult as well. See Theorem \ref{GANEG}, formula
(\ref{semibase}) and 
Theorem \ref{FFNEGAT} below, where the "row-column"
representations are considered in detail for arbitrary
root systems.

\vskip 0.2cm
Recently 
Berest, Etingof, and Ginzburg checked that the relation 
$k=\pm s/n,$ provided that $s\in \Z_+$ is not divisible by $n,$
is necessary for the existence of finite-dimensional 
irreducible representations of $\HH^o_n$
in the rational limit.  They also found that 
if such a representation exists then its dimension
is divisible by  $s^{n-1}$ and its $S_n$-character is
divisible by
the one for $\C[Q/sQ]$ for the root lattice 
$$
Q\ =\ 
\oplus_{i=1}^{n-1} \Z(e_i-e_{i+1})\ \subset\ 
\Z^n=\oplus_i \Z e_i.
$$

See also [CO] (the complete classification of irreducible
representations as $n=2$) and a paper by Dez\'el\'ee 
with some
results about finite-dimensional representations in the
rational limit. 

When restricted to the nonaffine
Hecke subalgebra $\H_n=\langle T_i,\ 0<i<n\ran,$
the module from v) becomes isomorphic to the 
$t$-counterpart
of $\C[Q/sQ].$ Here we use that the above $\h_n^o$-modules
associated with the decompositions
$n=n_0+\cdots+n_{s-1}$ are the corresponding parabolically
induced $\H_n$-modules. One may
use (\ref{semibase}) and
Theorem \ref{FFNEGAT} below to generalize this statement
to arbitrary root systems.

Thus the standard deformation argument, combined with
what was proven in the paper by 
Berest, Etingof, Ginzburg, readily shows 
that our ``column-row'' representations are 
well defined and
remain irreducible upon the rational degeneration, and, moreover,
that their character formula holds 
exactly, i.e. without an unknown multiplier. 

Actually here one
can employ a more general method, based on the
existence of an isomorphism of $\HH\ $ and its
rational degeneration $\HH''$ upon proper completion.
For instance, it holds for $q=e^h, t=q^k$ as $k\neq 0$
over the formal series in terms of $\sqrt{h},$ which was 
conjectured by Etingof. 
It can be proved by combining the Lusztig-type 
isomorphism discussed in the Appendix of [C1]
and the formula
connecting the trigonometric and rational {\it differential}
Dunkl operators. The first connects $\HH\ $ and its trigonometric
degeneration $\HH'$, the second leads to an
identification of $\HH'$ and $\HH''.$
It is not just a formal isomorphism.
It establishes an explicit isomorphism between finite-dimensional
representations of $\HH\ $ as $q,t$ are not roots of 
unity and
$\HH''$ for nonzero $k.$

\vfill
The representations from v) are unique
finite-dimensional representations among all irreducible 
representations, not only among the
$X$-semisimple ones. To see this we will generalize
the Bernstein-Zelevinsky classification.

We would like to mention recent work [Va] devoted to
the $K$-theoretic classification of irreducible representations
of $\HH\ $ 
for arbitrary roots systems (with coinciding $k$)
subject to the constraint $t^r=q^s$ from the Theorem.
It is similar to [KL1],
however involves more sophisticated geometric methods and 
combinatorial problems. In the case of $GL_n,$
the Vasserot classification  is expected to be equivalent
to Theorem \ref{HHGLNGEN}. It also leads to a 
description of 
finite-dimensional representations including the dimension
formula for arbitrary root systems. See Theorem 
\ref{FFNEGAT}
below.

\vskip 0.2cm
{\it General representations of $\HH_n$.}
The classification will be given 
in terms of periodic diagrams but with the following relaxation
of the conditions from (\ref{mini}):

\begin{align}
&n_i\ge n_{i+1} \hbox{\ or\ } \{n_{i}=n_{i+1}-1 
\hbox{ \ and\ }
m_i\ge m_{i+1}-1\}.
\label{minig}
\end{align}

We call such diagrams {\it generalized}.
Respectively, we define the {\it generalized pairs}
and then construct the corresponding 
weights using the same formulas as above. 
For generic $q,t,$ an arbitrary
irreducible representation of $\HH_n$ is induced from
an irreducible representation of $\h_n.$ This can be
proved following [C1]. Therefore all of them are
infinite-dimensional for such $q,t.$ 

\begin{theorem}\label{HHGLNGEN}
Imposing the relation $t^r=q^s$ for $r>0, s>0, (r,s)=1,$  
all $X$-cyclic irreducible $\HH_n$-modules up to
isomorphism and multiplication $\pi$ by a constant are in
one-to-one correspondence with the generalized new
pairs $\{\De,C\}$ up to the equivalence which is 
defined in the same way as above.
\qed 
\end{theorem}
 
These modules are quotients of the
induced modules $I=I_{\De,C}$ corresponding to 
the weights associated with $\{\De,C\}.$
The construction of such weights
involves the basic fundamental subdiagrams 
formed by consecutive rows of $\{\tde^p\}.$ The change
of such subdiagrams can be compensated by a proper change
of the corresponding $c$-numbers.
  
The induced modules $I$ are isomorphic to the affine Hecke
algebra 
$$\h_n^y\equal\lan \pi,\ T_i,\ 1\le i\le n-1\ran.
$$
Let $I'=I/(\sum\h_n^y(1-s_i))$ summed over the
simple (nonaffine) reflections $s_i$ corresponding to adjacent
boxes in the rows of the basic fundamental
subdiagrams of $\De.$ It is an $\HH_n$-module and 
the module from
the theorem is a unique irreducible nonzero 
quotient of $I'.$
Cf. [C11] and references therein.

Here $r>0.$
Note that we can skip the case of $r<0$ because 
it follows 
formally from the positive case 
thanks to the automorphism of 
$\HH_n$ sending 
$$X_i\mapsto X_i,\ \pi\mapsto \pi,\ q\mapsto q,\ \,  
T_i\mapsto -T_i,\ t\mapsto t^{-1}. $$

The $w$-periodicity and
other definitions remain unchanged. However
the following reservation 
about the fundamental subdiagrams and their 
$s$-partitions
is needed. 

We assume that they are intersections of
a given periodic diagram with sufficiently big {\it skew}
diagrams. Thus they can be of generalized type,
but their upper-left and lower-right boundaries 
in a given generalized $r$-periodic 
diagram are intersections with some skew paths. 

Similarly, the partitions are portions of fundamental
subdiagrams between consecutive skew paths. The $s$-partitions
of fundamental subpairs of $\{\De,C\}$
describe all {\it submodules} of the corresponding irreducible
representations upon the restriction to $\h_n.$ 
Here we continue considering the increasing
partitions: $\de>\de'$ if it holds for some
disjoint skew diagrams containing them. 

Concerning the finite-dimensional $\HH_n^o$-modules, it is 
clear that the irreducible $X$-cyclic representations
of $\HH_n^o$ have infinitely many non-isomorphic 
$\h_n^o$-submodules
unless in the case of the infinite column (recall that 
$r>0$).
More directly, it is not difficult to construct an 
infinite chain of $X$-eigenvectors with pairwise 
distinct eigenvalues if $\De$ is not a column. 
This construction essentially holds for 
arbitrary root systems however the 
finite-\-dimen\-sional representations are not described yet. 
There is another approach suggested by Vasserot.

\vfil 
\vskip 0.2cm
{\it Roots of unity.}
We go back to arbitrary root systems.
Now let $q$ be a primitive root of unity of order $N.$
Then $\hW_*^\flat[\xi]$ always contains 
\begin{align}
&A(N)\equal (N\cdot A)\cap B, \where
A\equal \{a\in P^\vee \mid (a,B)\subset \Z\},
\label{aforb}
\\ 
&\hbox{\ so\ } A=Q^\vee \for B=P, 
\hbox{\ and \ always\ } A\supset Q^\vee.
\notag \end{align}
Recall that $Q\subset Q^\vee$ and  $P\subset P^\vee,$
therefore  sometimes  $A\not\subset B$ and
$A(N)\neq N\cdot A.$

Here and further $q'=q^{1/(2\tilde{m})}$ is assumed to be 
a primitive root of order $2\tilde{m}N$ unless otherwise stated.
Recall
that $\tilde{m}\in \N$ is the smallest positive integer 
such that $\tilde{m}(B,B)\subset \Z.$ 
If $q'$ is a primitive root of order $m'N$ for  
$m'\,\mid\, 2\tilde{m},$
then we shall replace $A$ by $m'(2\tilde{m})^{-1}A$ in the 
definition of $A(N).$

Let us describe semisimple irreducible 
finite-dimensional quotients $V$ of $\i_X[\xi]$
assuming that $\hW_*^\flat[\xi]=A(N),$ i.e. the
stabilizer of $\xi$ is the smallest possible, 
and $\Up_+[\xi]=\hW^\flat=\Up_-[\xi].$
The main constraint (\ref{uniset})
becomes $\Up_0[\xi]=\hW^\flat$
in this case.

Theorem \ref{PUNIT} states that 
the spectrum of $\{X_b\}$ in any $V$ has to be
simple, the dimension is always
$|W|\cdot|B/A(N)|,\ $ and $\{V\}$ are  
in one-to-one correspondence
with one-dimensional characters $\varrho(a)$ of the group
$A(N),$ which determine the action via $S$ of $A(N)$ on
$V(\xi).$ 
To be more exact, the $V$-quotients up to isomorphisms
are in one-to-one correspondence
with the pairs 
$\{\,\varrho,\ \hbox{orbit \ }q^{\hW\llb \xi \rrb}\, \}.$
We denote them 
\begin{align}
&V\equal V[\xi,\varrho],\  
V[\xi,\varrho]\simeq V[\xi',\varrho]\iif
q^{\xi'}=q^{\hw\llb \xi\rrb}
\for \hw\in \hW.
\label{vxivar}
\end{align}
The isomorphisms here are established using
the intertwiners  $\Phi_{\hw}.$

Concerning the unitary structure,
let us calculate $\mu_\bullet$ on $A(N).$
Setting $N_\al=N/(N,\nu_\al),$ the scalar products 
$(a,\al^\vee)=(a,\al)/\nu_\al$
are divisible by $N_\al$ as $a\in A(N),\, \al\in R.$ 
On the other hand, $q_\al=q^{\nu_\al}$ is a 
primitive root of order
$N_\al.$ Using formula (\ref{lambi}),
\begin{align}
\mu_\bullet(a)\equal
&\prod_{\al\in R\mid (a,\al)>0}
\Bigl(
\frac{
t_\al^{-N_\al/2}-t_\al^{N_\al/2} q^{N_\al(\al,\xi)}}{
t_\al^{N_\al/2}-t_\al^{-N_\al/2} q^{N_\al(\al,\xi)}  
}
\Bigr)^{((a,\al^\vee)/N_\al)}\ =\notag\\ 
&\prod_{\al\in R_+}
\Bigl(
\frac{
t_\al^{-N_\al/2}-t_\al^{N_\al/2} q^{N_\al(\al,\xi)}}{
t_\al^{N_\al/2}-t_\al^{-N_\al/2} q^{N_\al(\al,\xi)}  
}
\Bigr)^{((a,\al^\vee)/N_\al)}
\label{muani}
\end{align}
for $a\in A(N).$ By the way, it is clear from this formula
that $\mu_\bullet$ is a homomorphism from $A(N)$ to the 
multiplicative group of the field  
$\Q(t_\al^N,q^{N(\xi,\al)}).$ It was already checked 
in the theorem.  
We come to the following
proposition.

\begin{proposition}\label{GENRE}
Assuming that $q$ is a primitive root of unity of 
degree $N$ and  
$$
\hW_*^\flat[\xi]\ =\ A(N),\ \Up_+[\xi]=\hW^\flat=\Up_-[\xi],
$$ 
$X$-semisimple irreducible quotients of $\i_X[\xi]$ 
exist if and 
only if
\begin{align}
&\Up_0[\xi]\ =\ \hW^\flat \hbox{\ i.e.\ } 
X_{\tal}(q^\xi)\neq 1
\hbox{\ for\ all\ } \tal\in \tR.
\label{genrep}
\end{align}
Such quotients $V[\xi,\varrho]$
are described by the  characters $\varrho$ of $A(N).$
The pseudo-unitary structure exists if and only if
$$
\varrho(a)^*\varrho(a)\ =
\ \mu_\bullet(a) \hbox{\ for\ the\ generators\ of\ }
A(N). 
$$
\end{proposition}

\vskip 0.2cm
{\it Comment on finite stabilizers.} 
When $q$ is a root of unity, 
the case of finite stabilizators
is, in a sense, opposite to the case of representations 
$V[\xi,\varrho].$ Indeed,
the group $\hW_*^\flat[\xi]$ is finite 
if and only if the set
$\Up_*[\xi]$ is finite. The check is simple. 
If the latter set is infinite than it contains
infinitely many elements in the form $a\hw$ for a 
fixed $\hw$
and  $a\in A(N),$ because the quotient 
$\hW/A(N)$ is finite.
Note that $A(N)$ is normal in $\hW,$ so this quotient 
is a group. 
We conclude that $\hW_*^\flat[\hw\llb \xi\rrb]$ is infinite. 
However it is conjugated  to $\hW_*^\flat[\xi]$ due to 
(\ref{stabc}), so the latter has to be infinite as well.

Recall that 
the existence of the pseudo-unitary structure
on $V$ reads:
\begin{align}
&S_{\hw}^\star S_{\hw} = \mu_\bullet(\hw)
\hbox{\ in\ } V(\xi) \for
\hw\in \hW_*^\flat[\xi] \hbox{\ of\ finite\ order.} 
\label{sunif}
\end{align}

We are going to check that here $\mu_\bullet(\hw)$
is $\pm 1$ and, moreover, automatically $1$
in many case. 

\begin{proposition}\label{SECORD}
For an element $\hw$  
in $\hW_*^\flat[\xi]$
of finite order,
$\mu_\bullet(\hw)=\pm 1$
without any assumptions on $q,t,\xi.$ Moreover,
$\mu_\bullet(\hw)=1$ if the order of $\hw$ is
$2M+1,$ i.e. odd, or $2(2M+1).$
\end{proposition}

{\it Proof.} First  we consider the case $\hw^2=1.$ 
Then $\la(\hw)=-\hw(\la(\hw))$ and each 
$X_{\tal}(q^\xi)$ 
appears in the product (\ref{muuni}) together with 
$X^{-1}_{\tal}(q^\xi)$ thanks to
$\hw(q^{\xi})=q^\xi.$ If $X_{\tal}(q^\xi)\neq$
$X^{-1}_{\tal}(q^\xi)$ then the corresponding binomials will
annihilate each other, i.e. their contribution to the product
will be $1.$ However they can conside. There are two 
possiblilities: $X_{\tal}(q^\xi)=\pm 1.$ 
If $X_{\tal}(q^\xi)=-1$ then  
\begin{align}
&
\Bigl(
\frac{
t_\al^{-1/2}-q_\al^jt_\al^{1/2} X_\al(q^{\xi})}{
t_\al^{1/2}-q_\al^jt_\al^{-1/2} X_\al(q^{\xi})
}
\Bigr)\ =\ 1.
\notag \end{align}
The equality $X_{\tal}(q^\xi)=1$ is impossible thanks to
the constraint (\ref{uniset}). So we get the complete
annihilation of all factors in $\mu_\bullet(\hw).$

Now we assume that $\hw$ is of finite order. Then
$\mu_\bullet(\hw)$ has to be a root of unity.
Under the reality conditions
$|q|=1,\ $ $k_\al\in \R,\ $ $ \xi\in \R^n\ ,$
$\mu_\bullet(\hw)$ is real and can be $\pm$ only.
Arbitrary triple $\{q,k,\xi\}$ 
can be continuously deformed to the triple above,
without changing the group  
$\hW_*^\flat[\xi]$ and the value of $\mu_\bullet(\hw).$
Indeed, the structure of $\hW_*^\flat[\xi]$ 
depends on certain multiplicative relations
among  $q_{(\al,\xi)}$ and $q_\al^{k_\al}.$
Therefore we need to deform $\{q,k,\xi\}$ within some
subvariety, which is possible. We get the claim about
$\pm 1.$

It gives that the orders of 
$$
\mu_\bullet(\hw) \and 
\mu_\bullet(\hw^{2M+1})=
\mu_\bullet(\hw)^{2M+1}
$$ 
always coincide. If the order of $\hw$ is odd or 
$2(2M+1),$ then 
$\mu_\bullet(\hw)$ has to be $1$ thanks to the claim about
the elements of the second order. 
\qed

\vskip 0.2cm
\section {Spherical representations} 
\setcounter{equation}{0}

It is a  continuation of the previous section. 
The parameters $q,t,\xi$ are still arbitrary
(nongeneric). We assume that 
$t_{\sht}\neq \pm 1\neq t_{\lng}$ 
unless otherwise stated. Let us continue the 
list of the basic
definitions.

\vskip 0.2cm
{\it Spherical and cospherical representations.}
A $\HH^\flat$-module $V$ is called $X$-{\it spherical} if 
it is generated by
an element $v,$ a {\it spherical vector}, such that
\begin{align}
&T_i(v)\ =\ X_{\al_i}(q^{\rho_k})v=\ 
t_i\, v
\for 1\le i\le n,\notag\\  
&X_a(v)\ =
\ X_{a}^{-1}(q^{\rho_k})v=q^{-(a\, ,\, \rho_k)}v
\for a\in B,
\label{txv}
\end{align}
It is $X$-cospherical if the space of linear form 
$\varpi$ satisfying
\begin{align}
&\varpi(T_i(u))\ =\ t_i\, \varpi(u)
\for 1\le i\le n,\notag\\  
&\varpi(X_a(u))\ =\ q^{-(a\, ,\, \rho_k)}\varpi(u)
\for a\in B,\ u\in V,
\label{txvv}
\end{align}
is nonzero, and at least one of them, called 
a {\it cospherical form}, 
does not vanish identically on each 
$\HH^\flat$- submodule of $V.$

We will mainly call the vectors $v$ satisfying 
(\ref{txv})
$X$-{\it invariants} and the cospherical forms 
$X$-{\it coinvariants}. 

It will be assumed that the spaces of 
invariants and coinvariants
are finite-dimensional. 
It is true for all modules in the category 
 $\mathcal{O}_X.$

Note that $\HH^\flat$-quotient of a spherical representation 
is spherical, and nonzero submodules of a cospherical module are
cospherical.

An $X$-spherical module is $X$-cospherical and vice versa
for pseudo-unitary modules $V$ in the category 
$\mathcal{O}_X.$ 
The bi-form identifies the spaces of invariants
$v$ and coinvariants $\varpi.$ Recall that it is 
nondegenerate $\star$-invariant and $\ast$-hermitian.
One can drop the last condition in this definition.

The cospherical form is constructed as follows. If $v$ is
a given spherical vector then the linear form $(u,v)$ 
is a coinvariant and
has no $\HH^\flat$-submodules in its kernel because the latter have
to be orthogonal to $\HH^\flat v= V.$ 
Therefore $(u,v)$ is a cospherical form. 

To construct an $X$-spherical vector from a 
$X$-cospherical
form, we use that the coinvariants vanish
on $V_X^\infty(\xi)$ 
(these spaces are finite-dimensional) unless
the $X$-character $\xi$ is from the second
line of (\ref{txv}). Indeed, 
they obviously vanish on $V_X(\xi)$ for such $\xi,$ so 
we can go to $V_X^\infty(\xi)/V_X(\xi)$ 
and continue by induction.
Therefore an arbitrary coinvariant 
$\varpi(u)$ can be represented as
$(u,v)$ for an invariant $v.$ 
The latter has to be a spherical vector
for a cospherical form $\varpi.$ Otherwise
the orthogonal complement of $\HH^\flat v\neq V$ would
belong to the kernel of $\varpi,$ which is impossible.
This complement is nonzero since $V\in\mathcal{O}_X.$

Concerning the coincidence of the space of invariants
and coinvariants for $V\in \mathcal{O}_X,$ 
it is not necessary to assume
that the form is $\HH^\flat$-invariant. 
It is sufficient to have a nondegenerate 
$\ast$-bilinear and $\h_X^\flat$-invariant form.

A $Y$-spherical vector, a cospherical form $\varpi,\,$ 
$Y$\--in\-variants,
and $Y$-co\-in\-va\-riants are defined by the relations: 
\begin{align}
&T_i(v) = t_i\, v,\  
\pi_r(v) = v,\for i>0,\ \pi_r\in \Pi^\flat, \notag\\  
&Y_a(v)\ =\ q^{(a\, ,\, \rho_k)}v \for a\in B,
\label{tyv}
\\ 
&\varpi(T_i(u)) = t_i\, \varpi(u),\ 
\varpi(Y_a(u)) = q^{(a\, ,\, \rho_k)}\varpi(u).
\notag \end{align}
Note the different signs of $ (a\, ,\, \rho_k)$
for $X$ and $Y.$

\vskip 0.2cm
A module $V$ with the one-dimensional space of 
invariants $v$ has a unique nonzero
spherical submodule. It is the $\HH^\flat$-span of $v.$
Dualizing, a module $V$ with the one-dimensional space of 
coinvariants
$\varpi$ has a unique nonzero cospherical quotient. 
It is $V$ divided by the sum of all 
$\HH^\flat$-submodules of $V$ inside the kernel of 
$\varpi.$ 
If $V$ is already (co)spherical then it has no proper
(co)spherical submodules or respectively quotients.

The main example is $\i_X[\xi]$ which has a unique  
$Y$-cospherical quotient because
the space of $Y$-coinvariants $\varpi$ from 
(\ref{tyv}) for the
module  $\i_X[\xi]$ is one-dimensional.
Cf. [C1], Lemma 6.2. 

The polynomial representation 
$\v\equal\Q_{q,t}[X_b]=\Q_{q,t}[X_b, b\in B]$ is $Y$-spherical 
and maps onto an arbitrary given $Y$-spherical module.  
The constant term is its $Y$-coinvariant. If $q,t$ are generic
then $\v$ is irreducible and automatically $Y$-cospherical.
When $q$ is a root of unity there are
infinitely many irreducible spherical modules.

In the definition of spherical and cospherical modules, 
one may take
any one-dimensional  character of the affine Hecke subalgebra 
with respect to $X$ or $Y$ in place of  
(\ref{txv}),  (\ref{txvv}),
(\ref{tyv}).
The spherical represenatation for generic $q,t$ remains isomorphic
to the space of polynomilas however with different
formulas for the $\hT$-operators depending on 
the choice of
the character. The functional representations, 
which are defined via the discretization of these operators,
also change. We will stick to the ``plus-character'' 
in this paper, 
leaving a straightforward generalization to the reader.

It is worth mentioning that relations between 
the spherical and 
cospherical representations constructed for different 
one-dimensional characters
play an important role in the theory of $\HH\ $  
as well as for the 
classical representations of affine Hecke algebras. An example is 
a general theory of the shift operators.
We will not consider these questions in the paper.

%\vfil
\vskip 0.2cm
{\it Primitive modules.}
The definition is of general nature but these modules 
are mainly needed for finite-dimensional representations. 
Recall that finite-dimensional $\HH^\flat$-modules exist
only when either $q$ is a root of unity or $q$ is generic 
and $k$ is special rational.  

By $Y$-{\it principal} spherical 
$\HH^\flat$-modules, we mean  $Y$-sphe\-ri\-cal 
$\HH^\flat$-modules
equipped with nondegenerate $\star$-invariant 
$\ast$-bilinear forms.
So these modules are also $Y$-cospherical.
Vice versa, the modules which are spherical and 
at the same time cospherical
have to be principal thanks to Proposition \ref{UNIABS}.
In this definition, we do not assume the pairing to be 
$\ast$-hermitian. We only need it to be nondegenerate.   

\vskip 0.2cm
A $Y$-principal module
is called {\it $Y$-primitive} if the dimension of the space
of $Y$-coinvariants is one. Using the nondegenerate form,
we conclude that the dimension of the space of $Y$-invariants 
has to be one as well. Note that these dimensions always coincide
for $Y$-principal spherical modules in the category $\mathcal{O}_Y,$
as it was discussed above (for $X$ instead of $Y$).

A primitive module has no nontrivial spherical submodules
and cospherical quotients.
The key property  of such  modules, a variant of the 
Schur Lemma,
is that a nonzero $\HH^\flat$-homomorphism
between two primitive modules $V_1\to V_2$ is an isomorphism. 
Indeed, the image of $V_1$ is spherical and therefore 
must coincide with $V_2.$ Hence $V_2$ 
is a cospherical quotient of $V_1$ which is impossible unless
the kernel is zero.
We also have the following proposition.

\vfil
\vskip 0.2cm
\begin{proposition} \label{PRIMLEM}
Let $V$ be a $Y$-principal module from $\mathcal{O}_Y$
embedded  into a finite
direct sum $V_{\hbox{\small sum}}$ 
$=\oplus V_i$ of $Y$-primitive modules such that
its intersections with all $V_i$ are nonzero.  
Then $V$  coincides with  $V_{\hbox{\small sum}}.$ 
\end{proposition}

{\it Proof.}
The $Y$-invariants of $V_{\hbox{\small sum}}$
are linear combinations  
$\sum c_i v_i$ of the spherical vectors 
$v_i\in V_i.$  
The space of $Y$-invariants of $V$ is smaller than
the space of invariants of
$V_{\hbox{\small sum}}$ if $V\neq V_{\hbox{\small sum}}.$

Similarly, $Y$-coinvariants of $V_{\hbox{\small sum}}$ 
are linear combinations   
$\varpi=\sum c_i \varpi_i$ of the cospherical forms  
$\varpi_i$ of $V_i$ natuarlly extended to the sum.
Since $V$ has nonzero intersections with $V_i$
all such $\varpi\neq 0$ remain nonzero when restricted to $V.$
Indeed, $c_i\neq 0$ for certain $i$ and if $\varpi(V)=0$
then the cospherical form $\varpi_i$ of $V_i$
vanishes on $V\cap V_i,$ which is impossible.

Using the coincidence of the
dimensions of the spaces of invariants and coinvariants
for $Y$-principal modules in $\mathcal{O}_Y$ we conclude that  
$V=V_{\hbox{\small sum}}.$
\qed
\vfil

Later we will see that a $Y$-primitive $X$-semisimple module is 
irreducible. Also
a $Y$-principal $X$-cyclic module, i.e. that generated 
by an $X$-eigenvector, is $Y$-primitive and a $Y$-spherical
$X$-cyclic irreducible module which has at least one 
coinvariant is primitive (the next proposition).

We note that a $Y$-spherical module such that
$V=\oplus_{\xi}V^{\infty}_X(\xi)$
has to be finite-dimensional
(in particular, one can take $V\in\mathcal{O}_X$). 
Indeed, $V=\oplus_{\xi}V^{\infty}_Y(\xi)$ because
it is a quotient of the $Y$-induced module with the 
"plus-character" as the weight. Finitely generated modules
which has such decompositions for $X$ and $Y$ together
are finite-dimensional. See the previous section.

Assuming that a finite-dimensional $Y$-spherical module 
$V$ possesses at least 
one $Y$-coinvariant $\varpi\neq 0,$  let 
$\tV$ be the quotient of  $V$ with respect to the intersection 
of all radicals of all $\ast$-bilinear $\star$-invariant 
forms. Since $V$ is finite-dimensional there exists a  
$Y$-coinvariant with the kernel precisely coinciding with the 
intersection of the kernels of all coinvariants. This  
coinvariant makes $\tV$ cospherical and therefore $Y$-principal. 
Obviously it is the universal 
$Y$-principal quotient of $V$. 

It is worth mentioning that
spherical (finite-dimensional)
representations in the classical theory of affine 
Hecke algebra are defined with respect
to its nonaffine Hecke subalgebra and always have isomorphic
spaces of invariants and coinvariants  
provided that the latter algebra is semisimple.
Therefore irreducible spherical modules are primitive in the
above sense, which readily gives that
(affine) primitive modules are nothing else but 
irreducible quotients of the polynomial representation
of the affine Hecke algebra.  

\vskip 0.2cm
 
In the next proposition, the field of constants is 
assumed algebraically closed. 
Note that the definition of
principal and primitive modules does not depend on the particular
choice of the field of constants. They remain principal
or primitive over any extension of   
$\Q_{q,t}$ and vice versa. If they are defined over $\Q_{q,t}$
and become principle/primitive over its extension 
then they
are already principle/primitive over $\Q_{q,t}.$
Recall that the $\ast$-bilinear forms above and below
are not supposed to be $\ast$-hermitian. 
By the radicals, we mean the left radicals of such forms.  

\begin{proposition} \label{PRIMS}
A $Y$-spherical module $V\in \mathcal{O}_X$ is principal
if and only if the dimensions of its $X$-eigenspaces 
$V(\xi)$ are no greater than one and there exists a
$Y$-coinvariant of $V$ which is nonzero at all 
nonzero $X$-eigenvectors. 
The sum  $V_{\hbox{\small sum}}$
$=\oplus V_i$ of $Y$-principal modules from 
$\mathcal{O}_X$
is $Y$-principal if and only if
\begin{align}
&  \hbox{Spec}_X(V_i)\, \cap\, \hbox{Spec}_X(V_j)\ =\ \emptyset
\for i\neq j.
\label{sumpr}
\end{align}
If $Y$-principal $V\in \mathcal{O}_X$ is generated by its 
$X$-eigenvectors then it is a direct sum of primitive modules.  
\end{proposition}

{\it Proof.} In the first place,
given a $Y$-spherical vector 
$v\in V$ (by definition, it generates $V$),
the coinvariants $\varpi$ are in 
one-to-one correspondence
with  $\star$-invariant $\ast$-bilinear forms on $V$. 
Namely, any such form can be represented
as $\varpi'( f_1 f_2^*)$ for $f_1,f_2\in \v$ 
where $\varpi'$ is the pullback of $\varpi$
with respect to the $\HH^\flat$-homomorphism $\v\to$      
$V,$ sending $1$ to $v.$ I.e. $(v_1,v_2)_\varpi=$
$\varpi'( f_1 f_2^*),$ where 
$$v_1=f_1(v)=f_1(X)v,\ v_2=f_2(v)=f_2(X)v.$$
It follows from Proposition \ref{UNIABS}. 

The form may be degenerate and
non-hermitian. Recall that hermitian forms 
(sarisfying $(u,v)=(v,u)^\ast$)
correspond to coinvariants of $\v$  
such that $\varpi'( f^\ast)= \varpi'( f)^\ast.$
We do not use hermitian forms in the proposition.

The radical of $(\ ,\ )_\varpi$ is 
the greatest $\HH^\flat$-submodule of $V$ inside the kernel
of $\varpi.$ Hence this form is nondegenerate if and only
if the coinvariant $\varpi$ is a cospherical vector. 

Let us check that
in the category $\mathcal{O}_X,$ a $Y$-coinvariant  
$\varpi$ is a cospherical form if and only if
\begin{align}
& \varpi(e)\neq 0 \hbox{ \ for\ all\ } 
X-\hbox{eigenvectors\ }
e\in V. 
\label{ndco}
\end{align}

These relations are obviously
sufficient since any submodule of $V$ has at least one 
nonzero $X$-eigenvector.
Let us verify that (\ref{ndco}) is necessary. 
If $ \varpi(e)=0$ for such $e,$ then
$\ \varpi( \h_Y^\flat\, \x(e))\ $ $=\ 0\ $ 
$=\ \varpi( \HH^\flat(e)),$ 
where by $\x$ we mean the algebra generated by $X_b$
for $b\in B.$

Now, if there exists a weight $\xi$ such that 
dim $V(\xi)>1,$ then $V$ cannot be cospherical because
any given  coinvariant vanishes on a proper linear
combination of $e_1$ and $e_2.$ 
Thus the first claim holds as well as relations
(\ref{sumpr}). 

Note that the same argument proves that the dimension
of the space of $Y$-coinvariants of a 
$X$-cyclic module is no greater than one.
In particular, $Y$-principal $X$-cyclic modules are
$Y$-primitive.
Indeed, otherwise $\HH^\flat v\neq V$ for
an arbitrary $X$-eigenvector $v$ of $V$ because
 $\HH^\flat v$ 
belongs to the kernel of a proper coinvarinat of $V$ vanishing at $v,$
which always exists if there are at least two non-proportional
coinvariants. 
      
Coming to the last claim,
for each $X$-eigenvector $e\neq 0$ in $V,$ let $\w_e$ be the 
hyperplane of all coinvariants of
$V$ vanishing at $e$ in the space $\w$ of all coinvariants.
The codimensions of $\w_e$ are  one
because $V$ is principal. The intersection of all $\w_e$
consists of the coinvariants vanishing at all 
$X$-eigenvectors of $V,$ so it is zero
because $V$ is generated by its $X$-eigenvectors. 

Since we have at least two different $\w_e$
we may pick finitely many  coinvariants $\varpi_i$ 
(two are sufficient for the induction step)
which do not belong all to
any particular $\w_e$ and consider the direct sum of
the quotients $V_i$ of $V$ by the corresponding radicals.
The resulting map $V\to\oplus V_i$ is injective.
Indeed, for any $X$-eigenvector $e^o\neq 0$ 
in its kernel, the hyperplane $\w_{e^o}$
would contain all $\varpi_i,$ which is impossible.   

Applying this procedure repeatedly to the resulting 
factors $V_i$ (the quotients of $V$ 
are generated by $X$-eigenvectors as well as $V$)
we can eventually make  all $V_i$ in the sum primitive. 
We may also assume that neither $V_i$ can be removed from
the sum $\oplus V_i,$ i.e. that the
intersections $V\cap V_i$ are all nonzero, and use the
previous proposition.
\qed

\vfil 
\vskip 0.2cm
{\bf Semisimple spherical representations.} 
Generally speaking, primitive modules can be reducible.
However in this paper
we will mainly stick to the semisimple modules.
Let us check that $Y$-primitive $X$-semisimple modules
are irreducible. 

If such $V$ has a 
$\HH^\flat$-submodule $U$, then its orthogonal complement $U^\bot$
with respect to the (unique) $\star$-invariant 
$\ast$-bilinear form has 
zero intersection with $U.$ Here we use that the $X$-spectrum of 
$Y$-principal module is simple. Hence $U\oplus U^\bot=$
$V.$ However this is impossible because the space of 
$Y$-invariants
of $V$ is  one-dimensional and generates $V$ as a 
$\HH^\flat$-module.

Let $V$ be an  irreducible $X$-semisimple $\HH^\flat$-module 
which is also $Y$-co\-sphe\-rical. Then
$V\, =\,$ $\i_X[\xi]/\j$ for  proper $\xi$ and proper 
$\HH^\flat$-submodule
$\j$ such that 
$\j_\xi\subset \j\subset  \i_X[\xi]$
for $\j_\xi\ =\ \sum \h_Y \Phi_{\hw},$ where the summation 
is over $\hw\in \ddot{\Up}_*[\xi].$

It is $Y$-cospherical if and 
only if
$\varpi_+(\j)=0$ for the plus-character of 
$\i_X[\xi]\simeq \h_Y^\flat:$ 
\begin{align}
&\varpi_+(T_i H)=t_i^{1/2}\varpi_+(H),\ 
\varpi_+(\pi_r H)=\varpi_+(H),
\ \varpi_+(1)=1.
\label{vapl}
\end{align}
The condition $\varpi_+(\j)=0$
readily results in $\Up_*[\xi]=\Up_+[\xi]$
since $\varpi_+(\Phi_{\hw})\neq 0$ for 
$\hw\in\ddot{\Up}_-[\xi].$  
We see that 
the condition $\Up_+[\xi]\subset \Up_0[\xi]$
is fulfilled for $\xi,$ and 
for any other $X$-weight $\xi',$ under the assumption
that $V$ is irreducible, $X$-semisimple,
and $Y$-cospherical.

In such $V,$
the eigenspace $V(\xi)$ is one-dimensional
and  $S_{\hw}=1$ in $V(\xi)$ for $\hw\in \hW^\flat[\xi].$
Indeed, $\varpi_+(S_{\hw}(1_\xi))=1$ in $\h_Y^\flat$ due 
to the 
normalization of  $S$ from (\ref{Phi}). Here 
$1_\xi$ is $1$ from
$\h_Y^\flat$ with the action of $\{X\}$ via $\xi.$ 
Therefore the same
holds in any quotients of $\h_Y^\flat$ by submodules belonging
to Ker$\, \varpi_+\, .$ Since $V(\xi)$ is an irreducible 
representation of $\hW^\flat[\xi],$ we get that 
$S_{\hw}\mapsto 1.$ 

Let us demonstrate that  
the conditions 

a) $\Up_+[\xi]\subset \Up_0[\xi],\ $ b) 
$\dim \,V(\xi)=1,$ and   

c) $S_{\hw}=1$ in $V(\xi)$ whenever  
$\hw\in \hW^\flat[\xi]$ 

\noindent 
are sufficient to make $V$ irreducible
and $Y$-cospherical.

Given a weight $\xi,$ first 
we need to check that  
\begin{align}
&\j_\xi=\sum \h_Y \Phi_{\hw}\, \neq\, \h_Y, 
\hbox{\ summed\ over\ } \ddot{\Up}_+[\xi],
\notag \end{align}
assuming that $\Phi_{\hw}$ are well defined. It 
is obvious since $\varpi_+(\j_\xi)=0.$ Moreover, 
the latter
relation gives that 
$U=\i_X[\xi]/\j_\xi$ has a unique cospherical quotient.
It is the quotient by 
the span of all $\HH^\flat$-submodules
in the kernel of $\varpi_+$ on $U.$ 
Let us denote it by $V.$
The restriction of $\varpi_+$ to the eigenspace
$V(\xi)$ is nonzero. Really,
otherwise it is identically zero on the
whole $V.$ The quotient of $V$ by the $\HH^\flat$-span of
the kernel of $\varpi_+$ on $V(\xi)$ is irreducible and 
cospherical. Hence it coincides with $V$ due to the uniqueness
of the cospherical quotients. Moreover,
dim$\,V(\xi)=1,$ and the action
of  $\hW^\flat[\xi]$ is trivial. 
We come to the following theorem.

\begin{theorem} 
i) The induced representation $\i_X[\xi]$ 
possesses a nonzero irreducible $Y$-co\-sphe\-rical 
$X$-semi\-simple quo\-tient $V$ if and only if 
\begin{align}
&\Up_0[\xi]\, \supset\, \Up_+[\xi]\, \subset\, 
\Up_-[\xi].
\label{sphset}
\end{align}
Such a quotient is unique and has one-dimensional 
$V(\xi)$
with the trivial ($S_{\hw}$ act as  $1$) action of 
the group  $\hW^\flat[\xi].$ Its $X$-spectrum is simple.

This $V$ is $X$-unitary
if and only if $\mu_\bullet(\hw)=1$ for all 
$\hw\in \hW_*^\flat[\xi].$
If such $V$ is  finite-dimensional then
it is spherical and, moreover, primitive, i.e.
has a unique  nonzero $Y$-invariant  and a 
unique $Y$-coinvariant up to proportionality.  
Setting $V=V_{\sph}[\xi],$ 
the modules $V_{\sph}[\xi']$ and $V_{\sph}[\xi]$ are isomorphic
if and only if $q^{\xi'}=q^{\hw\llb \xi\rrb}$ for 
$\hw\in \Up_+[\xi].$ 

ii) Finite-dimensional $V_{\sph}[\xi]$ can be 
identified with
\begin{align}
& \FF'[\xi]\equal \hbox{Funct}(\Up'_+[\xi],\Q_\xi)
\where \Up'_+[\xi]\equal\Up_+[\xi]/ \hW_+^\flat[\xi].
\label{funpr}
\end{align}
Here the action of $\HH^\flat$  is introduced via 
formulas
(\ref{Tfchar}) and is well defined on $\FF'[\xi].$
The characteristic functions 
$\{\chi_{\hw},\, \hw\in \Up'_+[\xi]\}$ 
form a basis of $\FF'[\xi]$ and are permuted by 
the intertwiners $S_{\hw}$ which
become $\hw$ in this module.

The pseudo-hermitian form of $V$
is proportional to 
\begin{align}
&\langle f,g\rangle'\equal 
\sum_{\hw\in \Up'_+[\xi]}
\mu_\bullet(\hw) f(\hw)\ g(\hw)^*.
\label{innerpr}  
\end{align}
\label{SPHSS}
\end{theorem}

{\it Proof.} The "cospherical" part of i) has been 
already checked.
The spectrum of $V$ is simple thanks to
Theorem \ref{PUNIT}.
The condition (\ref{muext}),  which is necessary and sufficient
for the existence of an invariant $X$-hermitian
form on $V,$ means that 
$\mu_\bullet(\hw)=1$ for all $\hw\in \hW_*^\flat[\xi].$ 
See also (\ref{sunif}).

Since the space of $Y$-coinvariants of $V_{\sph}[\xi]$ is
one-dimensional, any $\star$-invariant 
$\ast$-bilinear form on $V$ is hermitian
(skew-symmetric) up to proportionality. 
It is because $V$ has the "real" structure by the construction: 
we assume that $\ast$ fixes $q^\xi.$
 
Note that any pseudo-hermitian form  on
$V$  is $X$-definite. The scalar squares of $X$-eigenvectors are
nonzero because the $X$-eigenspaces are one-dimensional.
Indeed, if one of these squares vanishes then the 
radical of the form is a nonzero $\HH^\flat$-submodule 
of $V,$ which contradicts its irreducibility.

If such $V$ is also finite-dimensional 
then it is spherical and principal.
The  space of spherical vectors of $V$ is one-dimensional
because so is the space of $\varpi.$ 
See Proposition \ref{PRIMS}.

By the way, 
we get a surjective $\HH^\flat$-homomorphism $\v\to V$ sending $1$
to the spherical generator $1_{\sph}$  of $V.$
Its kernel has to be an intersection of maximal ideals of
$\v.$ This argument  gives another demonstration
of the simplicity of the $X$-spectrum of $V.$
%\vfil
  
Part ii) is nothing else but a reformulation of i). 
Really,
we may introduce (the action of) $s_i$ and $\hw$ in 
$V_{\sph}[\xi]$ as $S_i$ and $S_{\hw}.$ 
Then the formulas from ii) will
follow from i).
\qed

%\vfil
As an immediate application of the proposition, 
we get that
the irreducible  representations $V[\xi,\varrho]$ 
from Proposition \ref{GENRE}
cannot be, generally speaking, principal (spherical and cospherical). 
They are cospherical if and only if $\varrho=1.$ For such 
$\varrho,$ 
$\mu_\bullet(a)$ have to be $1$ on $A(N)$ to provide the
existence of the $X$-unitary structure. Formulas 
(\ref{muani})
combined with the constraint (\ref{genrep}) indicate that
it may happen only for very special values of the parameters.

%\vfil
\vskip 0.2cm
{\it Generic spherical representations.} 
Let us give a  general description
of spherical modules  as $q$
is a primitive root of unity of order $N\ge 0.$

We will use the lattices 
$A(N)= (N\cdot A)\cap B$ from (\ref{aforb}), where 
\begin{align}
&A= \frac{m'}{ 2\tilde{m}}\, \{a\in P^\vee \mid (a,B)\subset \Z\}, 
\ m'\hbox{\ divides\ } 2\tilde{m} 
\label{aforbb}
\end{align}
Where $\ q'=q^{1/(2\tilde{m})} 
\hbox{\ is \ a\ primitive\ root\ 
of \ order \ } m'N.$

To start with, 
let us assume that $t_\nu$ are generic and establish
that the quotient $\widehat{\v}$ of $\v$ by the radical 
of the standard form of the polynomial representation is 
well defined, $Y$-spherical, $Y$-cospherical,
$Y$-semisimple,  $X$-semisimple, and  
irreducible.

In the first place, all polynomials $\e_b$
are well defined thanks to the intertwining formulas.
See  (\ref{ebebnze}). Then 
$e_b=E_b(q^{-\rho_k})\e_b$ are well defined because
$E_b(q^{-\rho_k})\neq 0$ due to (\ref{ebebs}). However
$E_b\, (b\in B)\, $ form a basis of $\v.$ So do $\e_b,$
and the formulas (\ref{normhats}) for  
$\lan \e_c,\e_c\ran_\circ,$ which were obtained
for generic $q,$ hold in the considered case as well,
i.e. define the form $\lan\, ,\, \ran_\circ$ on $\v.$ 
Its radical $\v_o=$Rad$\, \lan\, ,\,\ran_\circ$
can be readily calculated.

\vfil
\begin{proposition}
i) For generic $t_\nu,$ the $\HH^\flat$-module 
$\widehat{\v}=\v/\v_0$
is well defined and
\begin{align}
&\widehat{\v}\simeq\oplus_b \Q_{q,t}\e_b
\hbox{\ as\ } 
-(\al_i^\vee,b_-)\le N_i\equal N/(N,\nu_i),\, 
1\le i\le n,
\label{vprime}
\\ 
&\{-(\al_i^\vee,b_-)=N_i\}\ \Rightarrow\  
\{u_b^{-1}(\al_i)\in R_-
\hbox{\ i.e.\ } \al_i\in \la(u_b^{-1})\}.
\notag \end{align}

ii) Given $b_-$ from (\ref{vprime}), let $W'$ and $W^\bullet$ be 
subgroups of $W$ generated respectively by $s_i$ such that
$(\al_i^\vee,b_-)=0$ and $-(\al_i^\vee,b_-)=N_i.$
Their longest elements will be denoted by $w'_0$ and 
$w^\bullet_0.$
Then all possible 
$u_b^{-1}$ from (\ref{vprime}) are represented in the form 
$$
w\,w^\bullet_0\for w\in W \hbox {\ such \ that\ }
l(w\, w^\bullet_0\, w'_0)=l(w)+l(w^\bullet_0)+
l(w'_0).
$$
The $Y$-spectrum of $\widehat{\v}$ is simple 
and this module is irreducible.
\label{VPRIME}
\end{proposition}
{\it Proof.} 
Given $b_-\in B_-,$
formula (\ref{lapiom}) provides the
following description of all possible $u_b:$
\begin{align}
&\al\in \la( u^{-1}) \Rightarrow (\al,b_-)\neq 0.
\label{lapiomm} 
\end{align}
Equivalently, $l(u^{-1}w'_0)=l(u^{-1})+l(w'_0)$ for 
the longest element $w'_0$ in the centralizer
$W'_0$ of $b_-$ in $W.$

If $-(\al_i,b_-)=N_i,$ then
$\al_i\in \la(u_b^{-1})$ for such $i,$   
i.e. $u_b^{-1}$ is divisible by $w^\bullet_0$ 
in  the following exact
sense: $l(u_b^{-1})$ $=l(u_b^{-1}\,w^\bullet_0)+
l(w^\bullet_0).$
Here $w^\bullet_0\in W^\bullet_0$ is the longest element.  
This gives the description of $u_b$ from the
proposition. 

Since $k$ is generic, $\e_b$ and $\e_c$ 
have coinciding $q^{b_\#}$ and $q^{c_\#}$ if and only
if $b_-=c_-\mod A(N)$ and $u_b=u_c.$ Recall that 
$A\subset B\subset P.$ So $b_--c_-\in N\cdot P$
in this case. 
Provided (\ref{vprime}) for $b$ and $c$,
we get that  
$b_-=c_-+\sum \pm N_j\,\om_j $ for some
indices $j$ and proper signs. Here 
$(b_-,\al_j^\vee)=N_j$
and $(c_-,\al_j^\vee)=0$ for $+N_j\om_j$
and the other way round for $-N_j\om_j.$

Let us assume that $u_b=u=u_c.$ Transposing 
$b_-$ and $c_-,$
we can find $j$ with the plus-sign 
of $N_j\om_j$ in this sum. 
Then $u^{-1}=u_b$ is divisible by $s_j.$
However $l(u^{-1}s_j)=l(u^{-1})+1$ because $u=u_c$.
This contradiction gives that
$u_b$ and $u_c$ cannot coincide and that
the spectrum Spec${}_Y(\widehat{\v})$
is simple. Since $\widehat{\v}$ is spherical 
(i.e. generated by $\e_0=1$)
and $Y$-pseudo-unitary it has to be irreducible.
\qed

\vfil
Note that $\widehat{\v}$ is  
$X$-semisimple in all known examples.
Moreover its $X$-spectrum  is simple.
It is likely that this is always true. The next 
theorem gives a different (and more general)
construction of spherical representations involving
the central characters. 

\begin{theorem}
i) The center $\z$ of $\HH^\flat$ contains finite sums
\begin{align}
&H=\sum_{a\in A(N)}  F_a\, S_a \hbox{\ such\ that\ } 
F_a\in \Q_{q,t}[X_b],
\label{zroo}
\\ 
& b( F_{a})=F_a \for b\in A(N),
\ F_{w(a)}=F_a \for w\in W, \and\notag\\ 
&F_a\cdot\prod_{[\al,j]} 
(t_{\al}^{1/2}q_\al^j X_{\al}-t_\al^{-1/2})^{-1}\,\in\, 
\Q_{q,t}[X_b],\ [\al,\nu_\al j]\in \la(a).
\notag \end{align}

ii) The algebra $\z_X\in \Q_{q,t}[X_b]$ consisting of 
$\,W\lsmash A(N)$-invariant polynomials 
is central as well as the algebra
$\z_Y$ of $\,W\lsmash A(N)$-invariant polynomials
in terms of $Y.$ 
The algebra $\HH^\flat$ is a free module of dimension 
$|W|\,|\, W^\flat/A(N)\,|^2$ over  $\z_X\cdot\z_Y.$

The elements (\ref{zroo}) for rational 
$F_a\in \Q_{q,t}(X_b)$
constitute the center $\z_{\hbox{\small loc}}$ of the localization  
$\HH_{\hbox{\small loc}}^\flat$ of 
$\HH^\flat$ by nonzero elements of $\z_X.$
The algebra $\HH^\flat_{\hbox{\small loc}}$ is  
of dimension $|\, W^\flat/A(N)\,|^2$ over  
$\z_{\hbox{\small loc}}.$

iii) Given a central $X$-character, a 
homomorphism of algebras $\ze:$ $ \z_X\to \Q_{q,t},$ 
let us denote its kernel by Ker$\,\ze.$ Then 
\begin{align}
&\v_\ze\,\equal\, \v\,/\,(\v\cdot\hbox{Ker}\ \ze)
\label{sroo}
\end{align}
is a $Y$-spherical module of dimension $|\, W^\flat/A(N)\,|.$ 
It is $X$-semisimple if and only if 
$\v\cdot\hbox{Ker}\,\ze$ is an intersection of maximal ideals
of $\v.$ 

iv) The module $\v_\ze$ is  irreducible for generic $\ze.$
If it is irreducible then 
$\v^\circ_\ze=$ Hom$(\v_\ze,\Q_{\xi})$ is isomorphic to  
$V[\xi,1]$  from Proposition \ref{GENRE},
where Ker$\,\ze\in$ Ker$\,\xi.$ The action of $\HH^\flat$ on
$\v^\circ_\ze$ is via the $\Q_{q,t}$-linear
anti-involution $H\mapsto H^\circ$ fixing 
 $T_i\, (0\le i\le n)\ ,$ fixing $X_b\, (b\in B),\ $ and
sending $\pi_r\mapsto \pi_r^{-1}.$ 
\label{SROO}
\end{theorem}
{\it Proof.}
Upon the $X$-localization, the defining property of the
intertwiners $S_{\hw}$ readily gives that
the elements (\ref{zroo}) for rational 
$F_a\in \Q_{q,t}(X_b)$
form the center of the localization  
$\HH^\flat_{\hbox{\small loc}}.$ 
Any element of the latter can be 
represented in this form with the summation over all 
$\hw\in \hW^\flat$
and arbitrary rational coefficients. Conjugating such sums
by $X_b$ and $S_{\hw},$ we
get  (\ref{zroo}). The integrality
condition is sufficient (but not necessary) to go back to 
$\HH^\flat.$

We readily get that $\z_X$ is central. 
The automorphisms of $PGL^c_2(\Z)$  preserve the center of 
$\HH^\flat.$ In particular, $\z_Y$ is central.
The calculation of the rank of this algebra over 
$\z_X\cdot\z_Y$
is a direct corollary of the PBW-theorem (see 
(\ref{hatdec})).

We leave the calculation of the  rank of
$\HH^\flat_{\hbox{\small loc}}$ over 
$\z_{\hbox{\small loc}}$ to the
reader;  it will not be used in the paper.
Note that the calculation of the rank of $\HH^\flat$ over the center
without the $X$-localization is a difficult problem even for
$q=1=t_{\sht}=t_{\lng}.$

Switching to $\v_\ze,$ its dimension is  
$|\, W^\flat/A(N)\,|$
because it is a quotient of $\HH^\flat$
by the left ideal generated by Ker$\,\ze$ and 
Ker$\,\varpi_+\subset \h_Y^\flat$ for 
$\varpi_+$ from (\ref{vapl}). 

The $X$-semisimplicity of  $\v_\ze$ readily implies that
the $X$-spectrum of this module is simple and $\ze$ 
is associated with $\xi.$ The module
 $\v^\circ_\ze=$ Hom$(\v_\ze,\Q_{\xi})$ is 
$Y$-cospherical and
$X$-semisimple with the same spectrum. 
The action of $\HH^\flat$ via $H\mapsto H^\circ$ is well defined
because this anti-involution preserves $\ze.$
This module is irreducible if and
only if  $\v_\ze$ is irreducible. We can apply
Theorem \ref{SPHSS} in this case. It results in
$\v^\circ_\ze=V[\xi,1].$
\qed

\vskip 0.2cm
\section {Gaussian and self-duality} 
\setcounter{equation}{0}

We follow the notation from the previous section.
A $\HH^\flat$- module $V$ is {\it self-dual} if the involution 
$\vep$ from (\ref{vep}) becomes
the conjugation by 
a certain $\ast$-linear involution of $V$ 
upon the restriction to End$(V).$ 
A module is called 
$PGL_2^c(\Z)$-invariant if $\vep, \tau_\pm$ act there
projectively, i.e. if
there are automorphisms of $V$ satisfying 
\begin{align}
&\tau_+\tau_-^{-1}\tau_+\equal \si=
\tau_-^{-1}\tau_+\tau_-^{-1},\
\tau_\pm=\vep\tau_\mp\vep^{-1},
\label{pglc}
\end{align}
and inducing the  $\vep, \tau_\pm$ from Section 2 on 
End$(V).$
By $PGL_2^c(\Z),$ we mean the central extension of  
$PGL_2(\Z)$ due
to Steinberg introduced  by the relations (\ref{pglc}).

Neither functional nor polynomial representations
are self-dual. Either $\tau_-$ or $\tau_+$ act there
respectively, but not all together.
Actually self-dual  semisimple irreducible modules 
(from $\o_X$)
have to be
finite-dimensional. Indeed, they belong to both 
$\o_X$ and $\o_Y,$ and one can readily check that
they are finite-dimensional.

The
non-cyclic module $\v_1\equal \tga^{-1}\v$   is
self-dual. Namely, the following $\ast$-linear
involution 
\begin{align}
&\psi(\e_b\tga^{-1})\ =\ q^{-(b_\#\, ,\, b_\#)/2+
(\rho_k\, ,\, \rho_k)/2}
\e_b\tga^{-1} \for b\in B 
\label{vepega}
\end{align}
corresponds  to $\vep$ from the  double Hecke algebra. 
It is formula (5.8) from
[C5], which is equivalent to (\ref{epeps}) above.
Remark that the module  $\v_1$ is not 
$PGL_2^c(\Z)$-invariant 
because any extension of the  automorphism $\tau_+$
has to be proportional to  the multiplication by 
$\tga(X),$ which
does not preserve $\v_1.$ 

\vskip 0.2cm
{\it Gaussians.}
We turn to the Gaussians in $V_{\sph}[\xi],$ which will be
called {\it restricted Gaussians}. 
Generally speaking, the problem is to make $\tau_+$ an
inner automorphism. The analysis is
not difficult thanks to the main theorem
of Section 6. The Gaussian commutes with $X$-operators,
so preserves the $X$-eigenspaces. This theorem
provides their description.
The Gaussian, if it exists, has to be an element 
(a function)
in spherical representations, not only an operator,
because the $X$-action makes these representations
algebras. So this case is especially interesting.

According to  
the main theorem,  $\tau_+$ is inner in $V=V_{\sph}[\xi]$ 
if and only if the map 
\begin{align}
&S^\xi_{\hw}\mapsto \tau_+(S^\xi_{\hw}),\where 
\hw\in \hW_*^\flat[\xi]\and  
S^\xi_{\hw}=S_{\hw}(q^\xi),
\label{staup}
\end{align}
is inner (i.e. multiplication
by a matrix) in the eigenspace $V(\xi).$
In the spherical case, $S_{\hw}$ become simply $\hw$
in the realization $\FF'[\xi],$ and the automorphisms
$\tau_+(\hw)$ can be
readily calculated. We come to the following proposition.

\begin{proposition}
The restricted Gaussian defined by the relations 
\begin{align}
&\ga\, H\, \ga^{-1}\ =\ \tau_+(H) \hbox{\ in\ } V
\label{gaussd}
\end{align}
exists in $V=\FF'[\xi]$ if and only if 
\begin{align}
&q^{(b,b+2\xi)/2}\ =\ 1 \hbox{\ whenever\ \ } 
 bw\in \hW_*^\flat[\xi] \for b\in B,\  w\in W. 
\label{gaussv}
\end{align}
Then it is proportional  to
$\ga_*(bw)\equal$ $q^{(b+2w(\xi),b)/2}.$
\label{GAUSSV}
\end{proposition}
{\it Proof.}  Let $v$ be the cyclic generator
of $V$ (of weight $\xi$). 
We identify $\ga$ with the operator
of multiplication by $\ga.$
If it exists then
\begin{align}
&\ga(S_{\hw}(v))\ =\ S_{\hw}(S_{\hw}^{-1}\ga S_{\hw})(v)
\hbox{\ for\ all\ }\hw\in \Up_+[\xi].
\label{gausfd}
\end{align}
Here $\ga'_{\hw}=S_{\hw}^{-1}\ga S_{\hw}\ga^{-1}$ is a function which
can be readily calculated due to the definition
of $\tau_+.$ However it is more convenient
to involve $\FF'[\xi],$ where $S_{\hw}$
is nothing else but  $\hw.$ Recall that $\tau_+$ is interpreted as
the formal conjugation by $\ga(q^z)=q^{(z,z)/2}$ in any
functional spaces. See (\ref{gausseq}) and (\ref{gaussf}).
Therefore  
\begin{align}
&\ga'_{\hw}(\hu)\ =\ 
q^{(\xi'',\xi'')/2-(\xi',\xi')/2}
\for \xi'=\hu\llb\xi\rrb,\ \xi''=\hw\llb\xi\rrb.
\label{gaussdd}
\end{align}
Its action on $v$ is $\ga'_{\hw}(\hbox{id})=
\ga_*(\hw).$

Thanks to the argument which has been used several times,
formula (\ref{gausfd}) defines $\ga$ correctly 
if and only if
$\ga'_{\hw}=1$ for $\hw\in\hW_*^\flat[\xi].$
This is exactly (\ref{gaussv}).
\qed

\vskip 0.2cm
{\bf Perfect representations.}
Following the same lines, let us find out when 
$V=V_{\sph}[\xi]$
is self-dual i.e. $\vep$ can be made an inner 
involution in $V.$
For the sake of simplicity, we use the functional realization of
$\FF'[\xi]$ of $V.$  So $1$ is the spherical vector, 
$\de_0=\de_{\hbox{\small id}}$ is
the $X$-cyclic generator,  
$\de_0=\de_{\hbox{\small id}}.$

Recall that 
\begin{align}
&\chi_{\hw}(\hu)\ =\ \de_{\hw,\hu},\ \
\de_{\hw}(\hu)\ =\ \mu_\bullet(\hw)^{-1}\chi_{\hw} 
\label{chide}
\end{align}
for the Kronecker $\de_{\hw,\hu}.$
We also set $\de^\pi_{b}=\de_{\pi_b}$ in the case
$\xi=-\rho_k.$

\begin{theorem} 
i) Provided (\ref{sphset}),
the finite-dimensional  representation $V_{\sph}[\xi]$ 
is self-dual
if and only if it is $X$-unitary and
$q^{-\rho_k}=q^{\hw\llb \xi\rrb}$ for a certain
$\hw\in \Up_+[\xi].$ The representation 
$V_{\sph}[-\rho_k],$
is $\HH^\flat$-isomorphic to
the irreducible $X$-cospherical $Y$-semisimple nonzero 
quotient of $\i_Y[\rho_k],$
which is unique. 
The $X$-spectrum of  $V_{\sph}[\xi],$ which is 
$\Up_+[-\rho_k],$ belongs to $\pi_B=\{\pi_b\, \mid\, 
b\in B\}$
and is the negative of its $Y$-spectrum.
The $X$-unitary structure of  $V_{\sph}[-\rho_k]$ is 
also $Y$-unitary. 

ii) The module $\FF'[-\rho_k]\simeq V_{\sph}[-\rho_k]$  
has a basis
formed by the discretizations  of the spherical 
polynomials $\e_b:$
\begin{align}
&\e'_b(\hw)\equal \e_b(q^{-\hw\llb \rho_k\rrb})\for 
\hw\in 
\Up_+[-\rho_k].
\label{eprime}
\end{align}
They are well defined 
assuming that $\pi_b\in \Up_+[-\rho_k].$ 
Explicitly,
$$
\FF'[-\rho_k]\ =\ \oplus_{b^\bullet}\,\Q_{q,t}\,
\e'_{b^\bullet}, 
$$
where the summation is over the set of representatives 
$$
\pi_{b^\bullet}\in \Up_+[-\rho_k]/ 
\hW^\flat_*[\rho_k]
\equal \Up'_+[-\rho_k].
$$ 
Here we use that 
$\hW^\flat_*[-\rho_k]=\hW^\flat_*[\rho_k].$ Moreover
\begin{align}
&\e'_b= \e'_c \hbox{\ whenever\ } \pi_b\llb \rho_k\rrb\ =
\pi_c\llb \rho_k\rrb, \hbox{\ i.e.\ } 
\pi_b^{-1}\pi_c\in \hW^\flat_*[\rho_k].
\label{epcoin}
\end{align}

iii) The map
\begin{align}
&\psi': \sum g_{b}\e'_{b}\mapsto 
\sum g_{b}^{*}\de^{\pi}_b
\for g_b\in \Q_{q,t},\ \pi_b\in  \Up_+[-\rho_k],
\label{isoprime}
\end{align}
is an $\HH^\flat$-automorphism of $\FF'[-\rho_k],$
inducing $\vep$ on $\HH^\flat.$ Here $\de^{\pi}_b$
depend only on the images of $\pi_b$ in
$\Up'_+[-\rho_k].$
Equivalently,
\begin{align}
&\lan 1\ran'\,\psi'(f)\,(\pi_b)\ =\ 
\lan \e'_b\, ,\,f\ran'\ =\ 
\lan\, f^*\,\e'_b\,\ran',\notag\\   
&\where \lan f\ran'\equal 
\sum_{\hw\in \Up'_+[-\rho_k]}
\mu_\bullet(\hw) f(\hw).
\label{fours}
\end{align}
\label{SDUAL}
\end{theorem}
{\it Proof.}
The condition  $q^{-\rho_k}=q^{\hw\llb \xi\rrb}$
is  necessary because $\vep$ sends $X_b$ to $Y_b,$  
conjugating the eigenvalues, and  the self-duality 
implies that  Spec$_X= -$Spec$_Y.$ 
So if the involution of $V,$ let us denote it by $\psi',$ 
exists then $\psi'(1)$  has to be proportional to $\de_0.$
Its $X$-weight is $-\rho_k.$ Here we use the 
$\FF'$-realization of $V.$
Recall that the $Y$-weight of $1$ is $\rho_k.$

Now let $\xi=-\rho_k$ and
\begin{align}
&\Up_0[\rho_k]\, \supset\, \Up_-[\rho_k]\, \subset\, 
\Up_+[\rho_k].
\label{sphrho}
\end{align}
It is a reformulation of (\ref{sphset}) using
the obvious relations  $\Up_\pm[\xi]=\Up_\mp[-\xi].$
We know how to describe nonzero
irreducible $Y$-semisimple 
quotients of $\i_Y[\rho_k].$ 
The latter module has such quotients
thanks to (\ref{sphrho}) and Theorem \ref{PUNIT}.
The universal semisimple quotient is 
$U_Y[\rho_k]\equal \h_X^{\flat}/ \j$ for
\begin{align}
&\j\equal\sum \h_X\cdot 
\tau_+\bigl(\Phi_{\hw}(q^{-\rho_k})\bigr), 
\label{uys}   
\end{align}
where the summation is over $\hw\in 
\ddot{\Up}_+[-\rho_k]$
and $1$ is identified with the cyclic generator of 
$\i_Y[\rho_k].$
Using Theorem \ref{SPHSS}, we get that
$U_Y[\rho_k]$ has a unique nonzero irreducible  $X$-cospherical
quotient.
This implies that the latter is a unique
irreducible $X$-cospherical quotient of $\i_Y[\rho_k].$ 
Cf. Lemma 6.2 from [C1]. We will denote it
by $V^Y_{\sph}[\rho_k].$ 

Let us check that the module 
$V_{\sph}[-\rho_k]\simeq \FF'[-\rho_k]$ 
is not only $Y$-co\-she\-rical but also
$X$-co\-sphe\-rical. 
Since it is irreducible, we need
to find a vector $v$ satisfying (\ref{txv}): 
\begin{align}
&T_i(v) =t_i\, v \for i>0,\ 
X_a(v)=q^{-(a\, ,\, \rho_k)}v \for a\in B.
\label{txvs}
\end{align}
Such $v$ is exactly the generator $\de_0.$ 
Note that the
first formula in (\ref{txvs}) holds because
$\{s_1,\ldots,s_n\}\subset$ $\ddot{\Up}_+[-\rho_k].$

Summarizing, the $\HH^\flat$-homomorphism
$\i_Y[\rho_k]\to$  $\FF'[-\rho_k]$ sending
$1\mapsto 1,$ which exists
because the former module is induced,
identifies  the latter
with $V^Y_{\sph}[\rho_k]$ constructed above. There are several
immediate corollaries.

First, the operators $Y_a$ are diagonalizable in  $\FF'[-\rho_k].$
Second, the $Y$-spectrum is simple and coincides with the
negative of the $X$-spectrum. Third, the invariant $Y$-definite 
$\ast$-bilinear form on  $V^Y_{\sph}[\rho_k]$ has to be 
proportional to the $X$-definite invariant 
form on  $\FF'[-\rho_k]$ upon this identification. Fourth,
the polynomials $\e_b\in \v$ are well defined for 
$\pi_b\in \Up_+[-\rho_k].$ Let us check the latter claim.  

Using the standard $\HH^\flat$-homomorphism 
from $\v$ to $\FF'[-\rho_k]$ 
and the simplicity of the  $Y$-spectrum of the target,
we diagonalize the $Y$-operators in the subspaces  
\begin{align}
&\Si(b)\ =\ \oplus_{c\succ b}\,\Q_{q,t}X_c\,
\oplus\Q_{q,t}X_b,
\for B \ni c\succ b,
\label{macdd}
\end{align}
starting with $b=0$ and ``decreasing'' $b$ as far as
$\pi_b \in \Up_+[-\rho_k].$ See (\ref{macd}).  
The leading monomial
$X_b$ always contributes to $E_b.$  Thus the Macdonald polynomials
$\{E_b\}$ are well defined in this range and their images $E'_b$ 
in $\FF'[-\rho_k]$ are nonzero.  
Let us verify that $E_b(q^{-\rho_k})\neq 0$ for such $b.$

Following [C3,C4], we set
\begin{align}
 &[\![f,g]\!]\ =\  \{L_{\imath(f)}(g(x))\}(q^{-\rho_k}) \for
f,g\in \v,
\label{symfg}
\\ 
&\imath(X_b)\ =\ X_{-b}\ =\ X_b^{-1},\ 
\imath(z)\ =\  z \for 
z\in \Q_{q,t}\, ,
\notag \end{align}
where $L_f$ is from Proposition \ref{YONE}.
Let us introduce the following
$\Q_{q,t}$-linear anti-\-involution of $\HH^\flat$
\begin{align}
&\triangledown\equal \vep\,\star=\star\,\vep:\ 
X_b\mapsto Y_b^{-1},\ 
T_i\mapsto T_i\ (1\le i\le n).
\label{triangledown}
\end{align}
The basis of all duality statements is the following
lemma.

\begin{lemma}\label{RADI}
For arbitrary nonzero $q,t_{\sht},t_{\lng},$
\begin{align}
&[\![f,g]\!]=[\![g,f]\!] \and
 [\![H(f),g]\!]=[\![f,H^\triangledown\,
(g)]\!],\ H\in \HH^\flat.
\label{syminv}
\end{align}
The quotient $\v'$ of $\v$ by the radical 
Rad$[\![\, ,\,]\!]$ 
of the pairing $[\![\ ,\ ]\!]$ is an 
$\HH^\flat$-module such that
a) all $Y$-eigenspaces $\v'_Y(\xi)$ are zero or 
one-dimensional,
b) $E(q^{-\rho_k})\neq 0$ if the image $E'$ of 
$E$ in  $\v'$
is a nonzero $Y$-eigenvector.
\end{lemma}
{\it Proof.} Formulas (\ref{syminv}) are from
Theorem 2.2 of [C4]. See also  [C1], Corollary 5.4.
Thus Rad$[\![\, ,\,]\!]$  is a submodule and the 
form  $[\![\ ,\ ]\!]$
is well defined and nondegenerate on $\v'.$
For any pullback $E\in \v$ of $E'\in \v',$
$E(q^{-\rho_k})=[\![E,1]\!]=$ $[\![E',1']\!].$
If $E'$ is a $Y$-eigenvector of weight $\xi$ and 
$E(q^{-\rho_k})$ vanishes then
$$
[\![\y^\flat (E'),\h_Y^\flat(1')]\!]=0= 
[\![E',\v\cdot\h_Y^\flat(1')]\!]
\for \y^\flat=\Q_{q,t}[Y_b].
$$
Therefore $ [\![E',\v']\!]=0,$ which is impossible.

We can use it to check that dim$\,\v'_Y[\xi]\le1.$ 
Indeed, there is always a linear combination of 
two eigenvectors in $\v'$ of the same weight
with zero value at $q^{-\rho_k},$ which has to be
zero identically.
\qed

Let us go back to the proof of the theorem.
Obviously the radical of the  pairing $[\![\ ,\ ]\!]$
belongs to the kernel of
the discretization map $\de':\,$ $\v\to\FF'[-\rho_k].$
Since the latter module is irreducible, the radical and 
the kernel coincide, and we have a nondegenerate pairing 
on $\FF'[-\rho_k]$ inducing $\triangledown$ on 
$\HH^\flat.$
Since the image of $E_b$ in this module
is nonzero for  $\pi_b\in \Up_+[-\rho_k],$ we conclude 
that  
$E_b(q^{-\rho_k})\neq 0$ and the spherical polynomials 
$\e_b=E_b/(E_b(q^{-\rho_k}))$ 
are well defined for such $b.$ 
Moreover, if $\e_b$ and $\e_c$ are well defined and their
eigenvalues coincide then they can be different in 
$\v$ but
their images in  $\FF'[-\rho_k]$ have to coincide thanks
to the normalization $\e_b(q^{-\rho_k})=1.$

So far we have not used the $X$-unitary structure at all. 
By the
way, if it exists then 
\begin{align}
&E_b(q^{-\rho_k})\ =\ 0\ \Rightarrow\ 
(\,E'_b\, ,\, \de_0\,)\ =\ 0\ 
\Rightarrow\notag\\ 
&(\, Y_a (E'_b)\, ,\, T_i X_c(\de_0)\,)=0\for i>0,\ 
a,c\in B \Rightarrow\notag\\ 
&(\, E'_b\, ,\, Y_a^{-1}T_i X_c(\de_0)\,)\ =\ 0\ \Rightarrow\
(\, E'_b\, ,\, \HH^\flat(\de_0)\,)\ =\ 0\ 
\Rightarrow\notag\\ 
&(\, E'_b\, ,\,\FF'[-\rho_k]\,)\ =
\ 0 \Rightarrow E'_b=0.
\label{enonz}
\end{align}

The existence of an $X$-unitary structure, i.e. the relation
$\mu_\bullet(\hw)=1$ for  $\hw\in \hW_*^\flat[\rho_k],$
is equivalent to the existence of $\psi'.$ Indeed,
the form $(f,g)=$ $[\![ \psi(f),\psi(g)]\!]$ is obviously 
invariant pseudo-hermitian,
assuming that $\psi$ induces $\vep,$ and the other way round.
Recall that by invariant pseudo-hermitian forms, we mean 
$\star$-invariant $\ast$-skew-symmetric nondegenerate forms.
As we already used, such a form is automatically 
$X$-definite because the $X$-eigenspaces of 
$\FF'[-\rho_k]$ 
are one-dimensional. 
It is equally applicable to $Y$ instead of $X.$

Recall that we do not suppose  $X$-definite or 
$Y$-definite forms to
be positive (or negative) hermitian forms. 

If $\mu_\bullet=1$ on $\hW_*^\flat[\rho_k],$ then
the delta-functions $\de^\pi_{b}$ depend only on the classes
of $\pi_b$ modulo the latter group. Hence,
the map $\psi'$ from iii) is a well defined restriction 
of $\psi_\circ$ from  formula (\ref{isofour}).
It induces $\vep$ on $\HH^\flat$ 
(Theorem \ref{ISOVD}) 
and
\begin{align}
&\e_b(X)\ =\ \vep( G_{\pi_b})(1)\for b\in B.
\label{isvepp}
\end{align}
The latter holds until the elements 
$\pi_b$ sit in $\Up_+[-\rho_k].$ 
\qed

Let us combine together all the structures discussed above.
Our aim is to describe finite-dimensional irreducible
self-dual semi\-simple pseudo-unitary 
re\-pre\-sen\-ta\-tion 
with the ac\-tion of $PGL_2^c(\Z).$
We call them {\it perfect}.

\begin{corollary}
i) A self-dual spherical $X,Y$-semisimple
pseudo-uni\-tary ir\-reducible
representation of $\HH^\flat$ is 
finite-dimensional and possesses
an $X$-eigenvector of weight $-\rho_k$ and a $Y$-eigenvector
of weight $\rho_k.$ It exists if and only if
\begin{align}
&\Up_0[-\rho_k] \supset \Up_+[-\rho_k] 
\subset \Up_-[-\rho_k],\ 
|\,\Up_+[-\rho_k]\,|\le \infty,\notag\\ 
&\and \mu_\bullet(\hw)=1 \for 
\hw\in \hW_*^\flat[\rho_k]. 
\label{sphsd}
\end{align}
There is only one such representation 
up to isomorphisms, namely,
$\FF'[-\rho_k].$ 

ii) It is $PGL_2^c(\Z)$-invariant if and only 
it contains the
restricted Gaussian 
$$ 
\ga_*(\pi_b)\equal q^{(b-2u_b^{-1}(\rho_k)\, ,\, b)/2}
\for b_{\#}=b-u_b^{-1}(\rho_k),
$$ 
which means the relations
\begin{align}
&q^{(b,b-2\rho_k)/2}\ =\ 1 \hbox{\ for \ all \ } 
\pi_b\in \hW_*^\flat[\rho_k]. 
\label{gaurho}
\end{align}
In this case, $\tau_+$ corresponds to the multiplication
by $\ga_*,$ and $\psi'$
from (\ref{isoprime}) is proportional to the  involution
of $\FF'[-\rho_k]$ sending
\begin{align}
&\e'_b\ga^{-1}\, \mapsto\, \ga_*(\pi_b)^{-1}\, 
\e'_b\ga_*^{-1}, 
\where \pi_b\in \Up_+[-\rho_k]\, .
\label{psega}
\end{align}
\label{SSSG}
\end{corollary}
{\it Proof.} Only the last formula requires some comment.
It is a straightforward  specialization of 
(\ref{vepega}). 
\qed

\vskip 0.2cm
{\bf Main examples.}
The above corollary gives an approach
to the classification of the triples
$\{q,t,B\subset P\}$ such that
$\HH^\flat$ possesses a perfect representation, i.e.
a finite-dimensional nonzero irreducible
self-dual semi\-simple pseudo-unitary 
re\-pre\-sen\-ta\-tion 
with the ac\-tion of $PGL_2^c(\Z).$ If it exists
then this representation is isomorphic to
$\FF'[-\rho_k].$  The importance of these 
representations is obvious. They carry all 
properties of the classical Fourier transform,
directly generalizing the {\it truncated Bessel 
functions}
considered in [CM], the Verlinde algebras, and
irreducible representations of the Weyl algebras
at roots of unity.
Here we will discuss 
only two "main sectors" of perfect representations,
negative (generic $q$) and positive (roots of unity),
where the description of $k$ is reasonably simple.
The case of $A_1$ is considered
in detail in [C7] and [CO] including the 
explicit formulas for the Gauss-Selberg sums.

\begin{theorem}
i) In notation from Theorem \ref{GANEG}, let us impose a
somewhat stronger variant of 
condition (\ref{munonz}),  assuming that $q$
is generic and
\begin{align}
&(\rho_k,\al^\vee)\ \not\in\ \Z \setminus\{0\}
\hbox{\ for\ all\ } \al\in R_+,
\label{notinz}
\\ 
h_\al(k)&=(\rho_k,\al^\vee)+k_\al\ \not\in\ 
\Z 
\hbox{\ for\ all\ extreme\ } \al\neq \vth,\ 
h_{\vth}(k) \in \ -\N.
\notag \end{align}
For instance, we may pick 
$k_{\lng}=k_{\sht}=-e/h,$ where $(e,h)=1$ for 
$e\in \N$ and $h=(\rho,\vth)+1$ is the Coxeter number. 
Then $\FF'[-\rho_k]$ is perfect. 
It exists and remains perfect if 
$q$ is a root of unity,
provided that $h_{\vth}(k)\ge -N$ 
and all fractional powers of $q$ which
may appear in the formulas are primitive
roots of unity of maximal possible order.

ii) Setting $e\equal -h_\vth(k),$ 
$$
\hW_*^\flat[\rho_k]\ =\ 
\{(e\om_r)u_r^{-1},\, \mid\, e\om_r\in B,\, r\in O\}.
$$ 
Identifying $b$ upon the symmetries 
$$
\pi_b\mapsto \pi_b\cdot u_r\cdot(e\om_r),
\hbox{\ i.e.\ } b\mapsto b+eu_b^{-1}(\om_r) \for
(e\om_r)\cdot u_r^{-1}\in \hW_*^\flat[\rho_k],
$$
\begin{align}
& \FF'[-\rho_k]\ =\ \sum_{b\in B} \Q_{q,t} \de_b^\pi, 
\hbox{\ where\ either\ \,} b=0, 
\label{ffstruct}
\\ 
&\hbox{\ or\ } (\vth,\rho_k)+k_\vth-(b_-,\vth)<0,
\hbox{\ or\ }\notag\\ 
&(\vth,\rho_k)+k_\vth-(b_-,\vth)=0 \and
u_b^{-1}(\vth)\in R_-.
\notag \end{align}
For example, let
$\om=\om_1,\ k=k_{\lng},\ e=-2k\in \N,\ B=P$ 
in the case $R=A_1.$ Then 
assumption (\ref{notinz}) means that $e$ is odd and  
$$\hW_*^\flat[\rho_k]=\{\pi_0,\pi_{e\om}\},\
\FF'[-\rho_k]=\sum_{j=0}^{e-1} \Q_{q,t}\, 
\de_{-j\om}^\pi.
$$
\label{FFNEGAT}
\end{theorem}
{\it The proof} is close to that of 
Theorem \ref{GANEG}. The explicit description of
the set $\Up_+[-\rho_k]$ is nothing else but 
the definition (\ref{upset}).

Let $bw\in \hW_*^\flat[\rho_k].$ 
See (\ref{hwflatxi}) and (\ref{unstab}). 
Then
$bw\llb\rho_k\rrb= w(\rho_k)+b=\rho_k,\ $ 
$w\in W,\,b\in B.$
Setting  $\be=w^{-1}(\al_i)$ for $1\le i\le n,$ 
$$(w(\rho_k)-\rho_k\, ,\, \al_i^\vee)=
(\rho_k\, ,\, \be^\vee)-k_\be=-(b,\al_i^\vee)\in \Z.
$$
If $\be>0$ then there exists a simple
root $\al_j$ of the same length as $\be$ such that 
$\be^\vee-\al_j^\vee$ is a positive coroot. Hence 
$\be=\al_j$
thanks to the first condition in (\ref{notinz}). 
Similarly,
$\be$ can be negative only for extreme $-\be.$
However it contradicts the second condition unless
$\be=-\vth.$ We see that $w$ preserves the set 
$\{-\vth,\al_1,\ldots,\al_n\}$ and is $u_r^{-1}$ for
certain $r\in O.$ We get that 
$$
b=\rho_k-u_r^{-1}(\rho_k)
=e\om_r \and \pi_{e\om_r}=(e\om_r)\cdot u_r^{-1}.
$$ 
Indeed,
$$
(b,\al_i)=0\,\Leftrightarrow\, 
u_r(\al_i)\neq -\vth
\,\Leftrightarrow\, i\neq r,\
(b,\al_r)=-h_{\vth}(k).
$$

There is an additional condition to check, namely,
(\ref{ffstruct}) for $b_-=-e\om_{r^*}$
and $u_r:$
$$
(e\om_{r^*},\vth)<e, \hbox{\ or\ }
\{\,(e\om_{r^*},\vth)=e \and
u_r^{-1}(\vth)\in R_-\, \}.
$$
Since $(\om_{r^*},\vth)=1,$  only the second relation
may happen. It 
holds because $u_r^{-1}(\vth)=-\al_r.$ 

Now let us establish the existence of the pseudo-unitary
structure, which is equivalent to the relations
$\mu_\bullet(\pi_{e\om_r})=1.$ Thanks to 
Proposition \ref{SECORD}, only $A_{4l-1}$ and $D_{2l+1}$
have to be examined. In other cases, the orders of the
elements of $bw\in \hW_*^\flat[\rho_k]$ are not divisible
by $4.$  The proof below is actually uniform but 
the calculation is more relaxed in the simply-laced case.

Letting $k=k_{\lng},$ the integer $e=-kh$ 
has to be relatively prime to $h$ and, in particular, 
must be odd.
Recall that $r\mapsto r^*$ for $r\in O$ describes the
inversion in $\Pi$ and corresponds to the automorphism
of the nonaffine Dynkin diagram induced by 
$\varsigma=-w_0.$
We will use (\ref{muhatpi}):
\begin{align}
&\mu_\bullet(\pi_b)\ =\ \prod
\Bigl(
\frac{
t_\al^{-1/2}-t_\al^{1/2} q_\al^{(\al^\vee,\rho_k)+j}}{
t_\al^{1/2}-t_\al^{-1/2} q_\al^{(\al^\vee,\rho_k)+j}
}
\Bigr) \for \al\in R_+,
\label{muhatpii}
\\ 
-( b_-, \al^\vee )&>j> 0 
\iif  u_b^{-1}(\al)\in R_-,\ 
-( b_-, \al^\vee )\ge j > 0 \hbox{\ otherwise}. 
\notag \end{align}
In the simply-laced case, $t_\al=t,\ q_\al=q,\ 
\al^\vee=\al.$
Let $b=e\om_r.$

One has $-(b_-,\al^\vee)=e(\om_{r^*},\al)=e$ 
when  $\al$ contains $\al_{r^*},$ and $=0$ otherwise.
Using that the total number of factors here is even, 
namely  
$2*(e-1)*(\rho,\om_{r^*}),$ we may transpose the
binomials in the numerator:
\begin{align}
&\mu_\bullet(\pi_{e\om_r}) = 
\prod_{(\om_{r^*},\al)=1}^{1<j<e}\, 
q^{2(\al,\rho_k)+2j}\,
\Bigl(
\frac{
t^{1/2}-t^{-1/2} q^{-(\al,\rho_k)-j}}{
t^{1/2}-t^{-1/2} q^{(\al,\rho_k)+j}
}
\Bigr).
\label{muhatpit}
\end{align}
Here $j<e$ because  $j=e$ may appear only under
the condition
$u_r^{-1}(\al)>0,$ which never holds.
Indeed, let $\al=\al_{r^*}+\be$ where
$\be$ doesn't contain $\al_{r^*}.$ Then 
$u_r^{-1}(\al_{r^*}+\be)= -\vth+\be'$ for $\be'$ without
$\al_r.$ Since $\vth$ contains all simple roots, we conclude
that $-\vth+\be'$ is always negative.

As a by-product, we have established that $-u_r^{-1},$ sending
\begin{align}
&\La\equal\{\al>0,\, (\al,\om_{r^*})=1\}\ni \al\,
\mapsto\notag\\ 
&\al'=-u_r^{-1}(\al)\in 
\{\al'>0,\, (\al',\om_{r})=1\}\equal\La',
\label{urminus}
\end{align}
is an isomorphism.
Its inverse takes
$$
(\al,\rho_k)\mapsto (-u_r(\al),\rho_k)=
-(\al, u_r^{-1}(\rho_k))=(\al,e\om_r-\rho_k)=
e-(\al,\rho_k).
$$
Since the sets $\La$ and $\La'$ are isomorphic under
the automorphism $\varsigma$ the sets 
$\{(\al,\rho_k)\}$ and
$\{(\al',\rho_k)\}$ coincide for $\al\in \La$ and
$\al'\in \La'.$ Finally, we conclude that
$(\al,\rho_k)\mapsto e-(\al,\rho_k)$ is a symmetry
of $\La=\La'.$ 

It ensures a complete cancelation of the binomials
in (\ref{muhatpit}) and results in
\begin{align}
&\mu_\bullet(\pi_{e\om_r})\ =\ q^{2\Si},\ 
\Si=\sum_{\al\in \La,1<j<e} ((\al,\rho_k)+j)\notag\\ 
&=\ (e-1)hk(\om_{r^*},\rho)+(e(e-1)/2)
(2(\om_{r^*},\rho))
\ =\ 0.
\label{mubupi}
\end{align}
We have used that $hk=e,$ $|\La|=2(\om_{r^*},\rho),$
and the following general formula 
$$\sum _{\al>0}(a,\al)(b,\al^\vee)\ =\ 
h(a,b),\ a,b\in \C^n.$$
Therefore $\FF'[-\rho_k]$ is pseudo-unitary.

The last check is (\ref{gaurho}), 
ensuring the existence of
the restricted Gaussian in this module:
\begin{align}
&(\om_r,e\om_r-2\rho_k)/2\ =\
 (\om_r , -u_r^{-1}(\rho_k)-\rho_k)/2
\label{gaurhoo}
\\ 
&= -(\om_r+u_r(\om_r) ,\rho_k)/2\ =\
 -(\om_{r^*}-\om_r,\rho_k)/2=0,\notag\\ 
& \where
\{(e\om_r)\cdot u_r^{-1}\}\ =\ \hW_*^\flat[\rho_k]. 
\notag \end{align}
Once again we have used that the automorphism 
$\varsigma=-w_0$
transposes $\om_r$ and $\om_{r^*}$ and preserves 
$\rho_k.$

The analysis of the case of $A_1$ is straightforward, as
well as the statement about roots of unity.
\qed

It is instructional to calculate the dimension of 
$\FF'[-\rho_k].$
For the sake of simplicity let $B=Q.$
We combine (\ref{ffstruct}) with formula (\ref{toaff})
applied to $z\in (1/e)Q,$ which gives the existence and
uniqueness of $\ u\in W,\ a\in Q$ such that
\begin{align}
& z_-\equal ua \llb z \rrb \in (1/e)P_-, \ 
(z_-,\vth)\ge -1, 
\label{toaffe}
\\ 
&(z_-,\vth)=-1  \ \Rightarrow\   
u^{-1}(\vth) \in R_-,
\hbox{\ and\ }\notag\\ 
&(\al_i,z_-)= 0 \ \Rightarrow\   
u^{-1}(\al_i)\in R_+,\ i>0.  
\notag \end{align}
Multiplying $z$ and $a$ by $e,$
we get that the elements $b$ from (\ref{ffstruct}) are 
in one-to-one correspondence with the elements of $Q/eQ.$
So the dimension is $e^n.$

There is an immediate proof of this formula based
on formula (\ref{semibase}). Indeed, it gives that the dimension
is the volume of the domain 
$\{z\in \R^n\}$ such that 
$(z,\al_i)\le 0$ for $i>0$ and  $(z,\vth)\ge -e$ 
divided by the volume
of the affine Weyl chamber.

\vskip 0.2cm
{\it Roots of unity.}
Let $q$ be a $N$-th root of unity under assumption
(\ref{aforbb}). See also (\ref{aforb}). However we do not
assume now that fractional powers of $q$ which appear in
the formulas are primitive roots of unity unless 
otherwise
stated. 

We are going to describe
the main "positive" sector of per\-fect 
$\HH^\flat$-mo\-dules, 
generalizing the Verlinde algebras. In this case 
$k_{\lng}$
and $k_{\sht}$ are positive and rational. 
To simplify considerations,
we will assume that $N$ is greater than the Coxeter
number, but this is actually not necessary. 
Perfect representations
below are well defined for small $N.$ To see one can follow
the proof of the previous theorem instead of 
using the affine Weyl chambers (see below).

There exist other sectors, for
instance, with negative $h_\vth(k)$ (see the above theorem)
and with more special choices of the roots of unity.
The complete list remains unknown, although Corollary   
\ref{SSSG} seems sufficient for the classification.

We pick 
$$\hat{q}=q^{1/(2\hat{m})}\for 
(B,B+2\rho_k)\ =\ \hat{m}^{-1}\Z
$$
setting 
$$q^{(b,c+2\rho_k)}\ =\ 
\hat{q}^{\, (b,c+2\rho_k)\hat{m}}
$$ unless stated otherwise. Note that
$\hat{q}$ is not supposed to be primitive.

\begin{theorem}
i) We assume that $k_{\lng},k_{\sht}>0,$ 
\begin{align}
&\N\ni h_\vth(k)=(\rho_k,\vth)+k_{\sht}\ <\ N_{\sht}=N,
\label{hveen}
\\ 
&\hbox{and\ also \ either\ (a)\ } (N,\nu_{\lng})=1, 
\hbox{\ or\ }\notag\\ 
&\hbox{(b)\ }
\N\ni h_\th(k)=(\rho_k,\th^\vee)+k_{\lng}\ <\ 
N_{\lng}=
N/(N,\nu_{\lng})\notag\\ 
&\hbox{for\ the\ longest
\ root\ } \th\in R_+, \hbox{\ or}\notag\\ 
&\hbox{(c)\ }
(\rho_k,\al^\vee)\ \not\in\ \Z_+
\hbox{\ for\ all\ long\ } \al\in R_+.
\notag \end{align}
In case (c), we suppose that $\hat{q}$ is 
a primitive root of unity.
Then relation (\ref{sphset}) holds, 
which gua\-ran\-tees the ex\-is\-ten\-ce
of the $Y$-co\-sphe\-rical $X$-semisimple
finite-dimensional irreducible 
module $\FF'[-\rho_k].$ It is nonzero if
$h_\vth(k)-N$ $\le \hbox{\, Max\,}(B_-,\vth).$
In case (b), we also add here the counterpart of 
this inequality for $\th.$

ii) We also impose either the condition
\begin{align}
&\{ \rho\in B,\, h_\vth(k)+h-1\le N,\,
 h_\th(k)+h-1<N_{\lng}\}\hbox{\ or\ }\notag\\ 
&\{\rho\not\in B,\, h_\vth(k)+h-1\le N+1,\,
 h_\th(k)+h-1<N_{\lng}+1 \}.
\label{hvevth}
\end{align}
Only the relations with $\vth$ are
necessary under conditions (a) or (c) from i).
Then $\FF'[-\rho_k]$ is pseudo-unitary and 
\begin{align}
&\ \hW_*^\flat[\rho_k]= \{(e\om_r)u_r^{-1}\, \mid \,
e\om_r\in B,\, q^{N(\om_r,B)}=1\}, 
\hbox{\ as\ } e\equal N-h_\vth(k).
\label{hvestab}
\end{align}
It is positive unitary if $q^{1/2}=\exp(\pi i/N).$
The restricted Gaussian exists in $\FF'[-\rho_k]$ 
if and only 
$$
q^{m_r}=\hat{q}^{\, 2m_r\hat{m}}=1 \for
m_r=(N-h_\vth(k))N(\om_r,\om_r)/2
$$ 
and $\om_r$ from (\ref{hvestab}). This may restrict 
the choices for $\hat{q}.$

iii) Provided the conditions from i) and ii),
\begin{align}
& \FF'[-\rho_k]\ =\Q_{q,t} \oplus
\sum_{0\neq b \in B} \Q_{q,t} \de_{b}^\pi, 
\hbox{\ where\ }  
\label{hvebstr}
\\ 
\hbox{(a)}\ &h_\vth(k)-(b_-,\vth)<N,\hbox{\ or}\notag\\ 
&\{h_\vth(k)-(b_-,\vth)=N \and
u_b^{-1}(\vth)\in R_-\}, \and\notag\\ 
\hbox{(b)}\  &h_\th(k)-(b_-,\th^\vee)<N_{\lng},
\hbox{\ or}\notag\\ 
&\{h_\th(k)-(b_-,\th^\vee)=N_{\lng} \and
u_b^{-1}(\th)\in R_-\}.
\notag \end{align}
Here we identify $b$ modulo $\hW_*^\flat[\rho_k]:$ 
$$
b\mapsto b+u_b^{-1}(c),\ b_-\mapsto u_c(b_-)+c_- \for 
\pi_c\in \hW_*^\flat[\rho_k].
$$ 
The space without this identification
is also an $\HH^\flat$-module, which always
contains the restricted Gaussian. 

iv) For example, let
$$\om=\om_1,\ k=k_{\lng},\ e=N-2k\in \N,\ B=P$$ 
in the case $R=A_1.$ Then
$$
\FF'[-\rho_k]=\sum_{j=-e+1}^{e} \Q_{q,t}\, 
\de_{j\om}^\pi
\hbox{\ for \ primitive\ } q^{1/2}
$$
is perfect. It 
is positive unitary for $q^{1/2}=\pm \exp(\pi i/N)$
if $k\in \N;$ the sign has to be plus for 
half-integral $k.$ 
If $k\in \N$ and $N$ is odd we may pick 
$q^{1/2}$ in primitive $N$-th roots of unity. 
Then the summation above is over $1\le j\le e.$
The Gaussian exists in $\FF'[-\rho_k]$ in this case if 
either $N=4l+3$ for integral $k$ or $k$ is
half-integral. This module is positive unitary for
$q^{1/2}=-\exp(\pi i/N).$
\label{HVEB}
\end{theorem}

{\it Proof.} 
The elements $\pi_b$ for $b$ satisfying
(\ref{hvebstr}) constitute
the set $\Up_+[-\rho_k]$ due to the positivity of
$k_\al.$ Here we use that $\rho_k\in P_+,$ 
$\vth$ is the maximal coroot, 
and $\th^\vee$ is maximal among the short coroots. 
The positivity also results in the inequalities
$$
(\rho_k,\al^\vee) -(b_-,\al)-\upsilon k_\al
+\ep^b_\al< N_\al=N/(\nu_\al,N)
$$ 
for the same $b$ and all $\al>0,$ 
where $\upsilon=0,1,\ $
$\ep^b_\al=1$ if $u_b^{-1}(\al)>0$
and $0$ otherwise. These inequalities are necessary and
sufficient for the existence of 
$\FF'[-\rho_k].$ Let us calculate 
the stabilizer $\hW_*^\flat[\rho_k].$

\begin{lemma}
Under the assumptions from ii), the group of
the elements $\pi_c$ such that 
$\Up_+[-\rho_k]\pi_c=\Up_+[-\rho_k]$ is exactly
$\Pi=\{(e\om_r)u_r^{-1},r\in O\}$ for $e=N-h_\vth(k).$
\label{UPLEMMA}
\end{lemma}
{\it Proof.}
By definition,  $\pi_c=c u_c^{-1}$ 
preserve the set of all element $b$ 
satisfying (\ref{hvebstr}) under the mapping 
$$
b\mapsto b',\where \pi_b\pi_c=\pi_{b'},\ 
b'=b+u_b^{-1}(c),
\ u_{b'}=u_c u_b.
$$ 

The corresponding set of
$b_-$ will be denoted by $B_-^N.$
In terms of $b_-:$ 
$$
\pi_b\pi_c=u_b^{-1}b_-u_c^{-1}c_-=
u_b^{-1}u_c^{-1}u_c(b_-)c_-=u_{b'}^{-1}
\cdot(u_c(b_-)+c_-).
$$
The element $b'_-=u_c(b_-)+c_-$ has to be from $B_-^N.$
Using the affine action (\ref{afaction}): 
\begin{align}
&b_-\mapsto b'_-\ =\ c_- u_c\llb b_-\rrb\ =\ 
u_c c\llb b_-\rrb.
\label{hvebmap}
\end{align}

Due to conditions (\ref{hvevth}) from part ii), 
either $-\rho\in B_-^N$ (it happens
when $\rho\in B$) or
$-\rho-\om_r\in B_-^N$ for proper 
minuscule $\om_r$ otherwise.
Inner points of the segment 
connecting $0$ and such $b_-$ in $\R^n$ are also inner
in the polyhedron $\bar{B}_-^N\in \R^n$ which is defined
by the same inequalities but for real $\bar{b}_-.$
The set of its inner points is either $-e\CC_a,$
where $e= N-h_\vth(k),$ $\CC_a$ is the
affine Weyl chamber, or its intersection with the counterpart
for $\th$ and $N_{\lng}-h_\th(k)$ in case (b) from i). 
See (\ref{lamaff}) for the definition of the affine Weyl 
chamber and its basic properties.

We get that some points do not
leave the negative Weyl chamber $-\CC_a$ under
$$
z\mapsto e^{-1}u_c c\llb\, e z\, \rrb.
$$ 
So do all points. The stabilizer of
$-\CC_a=w_0(\CC_a)$ is $\{(-\om_{r^*})u_{r^*}^{-1}\}.$ 
Therefore
\begin{align}
&c_- u_c\ =\ (-e\om_{r^*}) u_{r^*}^{-1} \for r\in O
\and c_-\ =\ -e\om_{r^*},\, u_c\ =\ u_{r^*}^{-1},
\notag\\  
&\pi_c\ =\ u_c^{-1}c_-\ =\ 
u_{r^*}\cdot (-e\om_{r^*})\ =\
(e\om_{r})u_{r}^{-1},\ c\ =\ e\om_{r}.
\notag \end{align}
\qed

The elements $\pi_c\in \hW_*^\flat[\rho_k]$ satisfy the
conditions of the lemma.
See (\ref{hwflatxi}) and (\ref{unstab}). 
The definition is
as follows:  
\begin{align}
&\hW^\flat[\rho]\ =\ 
\{\hw\in \Up_+[-\rho_k],\ 
q^{(\hw\llb\rho_k\rrb-\rho_k\, ,\, b)}= 1 
\hbox{\ for\ all\ } b\in B\}.
\label{hwflatpic}
\end{align}
Setting $\hw=(e\om_{r})u_{r}^{-1},$
$$
\hw\llb\rho_k\rrb-\rho_k=u_r^{-1}(\rho_k)-\rho_k+
e\om_{r}=h_\vth(k)\om_r+e\om_{r}=N\om_r.
$$
Cf. the proof of the previous theorem.
Hence $q^{N(\om_r,b)}$ has to be $1$ for all $b\in B.$
We arrive at the conditions from ii). Now we can simply
follow Theorem \ref{FFNEGAT}, formulas (\ref{mubupi}) and
(\ref{gaurhoo}), to check the existence of the 
pseudo-unitary structure and the restricted Gaussian.

The claim about the positivity of the pairing for
the "smallest" root of unity $q^{1/2}=\pm \exp(\pi i/N)$ 
is straightforward too. We represent the binomials in
formula (\ref{muhatpit}) for $\mu_\bullet$
in the form $(e^{ix}-e^{-ix})$ and use inequalities
(\ref{hvebstr}). In the numerator,
the $x$ do not reach $\pi$
and $e^{ix}-e^{-ix}=ic$ for positive numbers $c.$
Therefore the same holds true in the denominator. 
\qed

As an application, we get 
a counterpart of Theorem \ref{EEJACK}
at roots of unity. Namely
\begin{align}
\ \ \langle
\e'_b\, (\e'_c)^*\,\ga{\mu_\bullet}
\rangle'_\#\ &=\
 q^{-(b_\#,b_\#)/2-(c_\#,c_\#)/2 +(\rho_k,\rho_k)} 
\e_c(q^{b_\#})\langle \ga{\mu_\bullet}\rangle'_\#,
\label{eejacksr}
 \\ 
\ \ \langle \e'_b\, \e'_c \ga\mu_\bullet
\rangle'_\#\ &=\ 
q^{-(b_\#,b_\#)/2-(c_\#,c_\#)/2 +(\rho_k,\rho_k)}
\notag\\  
&\times\prod_{\nu}t_\nu ^{l_\nu(w_0)/2}
T^{-1}_{w_0}(\e_c)(q^{b_\#})
\langle \ga\mu_\bullet\rangle'_\#.
\label{eejacr}
\end{align}
Here $\lan f\mu_\bullet\ran'_\#$ is the sum 
$\sum f\mu_\bullet(\pi_{b})$ over
$\{\pi_{b}\}$ for $b$ satisfying (\ref{hvebstr})
upon the identification modulo $\hW^\flat_*[\rho_k].$
The functions $\e'_b,\e'_c$ are the images of
the polynomials $\e_b,\e_c$ in $\FF'[-\rho_k]$
for $b,c$ satisfying the same constraint. They are
well defined and nonzero.
See Theorem \ref{SDUAL}.

Concerning the generalized Gauss-Selberg sums
$\langle \ga{\mu_\bullet}\rangle'_\#\, ,$
the formulas can be obtained using the shift
operators, which we do not discuss in the paper.
See [C6] about integral $k$ and [C7] with
the complete list of formulas for $A_1.$ 
Generally speaking,
the shift operator can be used as follows.

Given $k,$ first we calculate 
$\langle \ga\mu^\ka_\bullet\rangle'_\#$
for $\varkappa=\{\varkappa_{\sht},\varkappa_{\lng}\}$ 
taken from the sets $k_{\sht}+\Z$
and  $k_{\lng}+\Z$ with the simplest possible 
perfect representations. 
We denote the corresponding $\mu$
by $\mu^\varkappa.$ If $k$ are fractional we take $\varkappa$ 
negative and use the Macdonald
identities from Section 5. Then we may apply 
the counterpart of (\ref{hatmu}):
\begin{align}
& \frac{\langle \ga\mu^k_\bullet\rangle'_\#\ 
}{ \langle \ga\mu^\varkappa_\bullet\rangle'_\# }
=\prod_{\al\in R_+}\prod_{ j=1}^{\infty}
\Bigl(\frac{ 
1- q_\al^{(\rho_k,\al^\vee)+j}}{
1-t_\al^{-1}q_\al^{(\rho_k,\al^\vee)+j} }\Bigr)
\Bigl(\frac{1-t_\al^{-1}q_\al^{(\rho_\varkappa,\al^\vee)+j} 
}{ 1- q_\al^{(\rho_\varkappa,\al^\vee)+j}
 }\Bigr). 
\label{hatmur}  
\end{align}


\begin{thebibliography}{MTV}
%\vfil
\bibitem [Ao] {Ao}
{K.~Aomoto},
{On product formula for Jackson integrals associated with root
systems}, Preprint (1994).

\bibitem [AI] {AI}
{R.~Askey}, and {M.E.H.~Ismail},
{\it A generalization of ultraspherical polynomials},
in {\it Studies in Pure Mathematics} (ed. P.~Erd\"os),
Birkh\"auser (1983), 55--78.

\bibitem [Bi] {Bi}
{J.~Birman},
{\it On Braid groups},
Commun.\ Pure\ Appl.\ Math.,
{\bf 22} (1969), 41--72.

\bibitem [Bo] {Bo}
{N.~Bourbaki},
{\it Groupes et alg\`ebres de Lie}, Ch. {\bf 4--6},
Hermann, Paris (1969).

\bibitem [C1] {C1}
{I.~Cherednik},
{\it Intertwining operators of double affine Hecke algebras},
Selecta Math. New ser. {3} (1997), 459--495

\bibitem [C2] {C2}
\bysame,{\it Double affine Hecke algebras and  Macdonald's
conjectures},
Annals of Mathematics {141} (1995), 191--216.

\bibitem [C3] {C3}
\bysame, 
{\it Macdonald's evaluation conjectures and
difference Fourier transform},
Inventiones Math. {122} (1995),119--145.

\bibitem [C4] {C4}
\bysame, 
{\it Nonsymmetric Macdonald polynomials },
IMRN {10} (1995), 483--515.

\bibitem [C5] {C5}
\bysame,
{\it Difference Macdonald-Mehta conjecture},
IMRN {10} (1997), 449--467.

\bibitem [C6] {C6}
\bysame,
{\it From double Hecke algebra to analysis },
Doc.Math.J.DMV, Extra Volume ICM 1998,II, 527--531.

\bibitem [C7] {C7}
\bysame,
{\it One-dimensional double Hecke algebras and Gaussian
sums }, Duke Math. J. {108}:3  (2001), 511--538.

\bibitem [C8] {C8}
\bysame,
{\it Inverse Harish--Chandra transform and 
difference operators}
IMRN { 15} (1997), 733--750.

\bibitem [C9] {C9}
\bysame,
{\it A new interpretation of Gelfand-Tzetlin bases },
Duke Math J., { 54}:2 (1987), 563--577.

\bibitem [C10] {C10}
\bysame
{\it On $q$-analogues of Riemann's zeta function}
Sel.math., New ser. { 7} (2001), 1--44.

\bibitem [C11] {C11}
\bysame
{\it On special bases of irreducible finite-dimensional 
representations of the degenerate affine Hecke algebra }
Funct. Analysis Appl. { 20} (1986), 87--89.

\bibitem [CM] {CM}
\bysame, and {Y.~Markov},
{\it Hankel transform via double Hecke algebra },
CIME publications, 2001.

\bibitem [CO] {CO}
\bysame, and {V.~Ostrik},
{\it From Double Hecke Algebras to
Fourier Transform}, Selecta Math. (2002). 
\vfil

\bibitem [Ch] {Ch}
{K.~Chandrasekharan},
{\it Elliptic functions}, 
A series of Comprehensive Studies in 
Math. {281},
Springer-Verlag, Berlin, Heidelberg, New York, Tokyo (1985).

\bibitem [DS] {DS}
{J.F.~van Diejen}, and {J.V.~Stokman}, 
{\it Multivariable $q$-Racah polynomials}.
Duke Math. J. { 91} (1998), 89--136.

\bibitem [D] {D}
{C.F.~Dunkl}, 
{\it Hankel transforms associated to finite reflection groups},
Contemp. Math. {138} (1992), 123--138.

\bibitem [Ev] {Ev}
{R.J.~Evans},
{\it The evaluation of Selberg character sums},
L'Enseignement Math.   {37} (1991), 235--248.

\bibitem [J] {J}
{M.F.E.~de Jeu},
{\it The Dunkl transform }, Invent. Math. {113} (1993), 147--162.

\bibitem [HO1] {HO1}
{G.J.~Heckman}, and {E.M.~Opdam},
{\it Root systems and hypergeometric functions I},
Comp. Math. { 64} (1987), 329--352.

\bibitem [HO2] {HO2}
\bysame, and \bysame
{\it Harmonic analysis for affine Hecke algebras},
Preprint (1996).

\bibitem [He] {He}
{S.~Helgason},
{\it Groups and geometric analysis}, 
Academic Press, New York (1984). 

\bibitem [K] {K}
 {V.G.~Kac},
{\it Infinite dimensional Lie algebras },
Cambridge University Press, Cambridge (1990).

\bibitem [KL1] {KL1}
{D.~Kazhdan}, and { G.~Lusztig},
{\it Proof of the Deligne-Langlands conjecture for Hecke algebras},
Invent.Math. {  87} (1987), 153--215.

\bibitem [KL2] {KL2}
{D.~Kazhdan}, and { G.~Lusztig},
{\it Tensor structures arising from affine Lie algebras. III,}
J. of AMS { 7} (1994), 335--381.

\bibitem [Ki] {Ki}
{A.~Kirillov, Jr.},
{\it Inner product on conformal blocks and Macdonald's
polynomials at roots of unity }, Preprint (1995).

\bibitem [KnS] {KnS}
 {F.~Knop}, and {S.~Sahi},
{A recursion  and a combinatorial formula for Jack
polynomials}, Inventions Math., (1997).

\bibitem [KS] {KS}
{E.~Koelink}, and { J.~Stokman},
{\it The big $q$-Jacobi function transform },
Publ. IRMA (Univ. Lois Pasteur), 1999/23 (1999).

\bibitem [KK] {KK}
{B.~Kostant}, and { S.~Kumar},
{\it  T-Equivariant K-theory of generalized flag varieties},
J. Diff. Geometry {  32} (1990), 549--603.

\bibitem [L] {L}
{G.~Lusztig}, 
{Affine Hecke algebras and their graded version},
J. of the AMS { 2}:3(1989), 599--685.

\bibitem [M1] {M1}
{I.G.~Macdonald}, {\it  Some conjectures for root systems},
SIAM J.Math. Anal. { 13}:6 (1982), 988--1007.

\bibitem [M2] {M2}
\bysame, {\it  A new class of symmetric functions },
Publ.I.R.M.A., Strasbourg, Actes 20-e Seminaire Lotharingen,
(1988), 131--171 .

\bibitem [M3] {M3}
\bysame, {\it Orthogonal poly\-no\-mials asso\-cia\-ted with root 
systems}, Preprint (1988).

\bibitem [M4] {M4}
\bysame, {\it Affine Hecke algebras and orthogonal polynomials},
S\'eminaire Bourbaki { 47}:797 (1995), 01--18.

\bibitem [M5] {M5}
\bysame, {\it A formal identity for affine root systems},
Preprint (1996).

\bibitem [M6] {M6}
\bysame, {\it Affine root systems and Dedekind's $\eta$-function},
Inventiones Math. { 15} (1972), 91--143.

\bibitem [Ma] {Ma}
{H.~Matsumoto},
{\it Analyse harmonique dans les systemes de
Tits bornologiques de type affine},
Lecture Notes in Math.  { 590} (1977).

\bibitem [Na] {Na}
{M.~Nazarov},
{\it Young's symmetrizers for projective 
representations of the symmetric group},
Adv. Math. { 127} (1997), 190--257.

\bibitem [O1] {O1}
{E.M.~Opdam}, 
{\it  Some applications of hypergeometric shift
operators}, Invent.Math.{  98} (1989), 1--18.

\bibitem [O2] {O2}
\bysame, {\it Harmonic analysis for certain 
representations of graded Hecke algebras}, 
Acta Math. { 175} (1995), 75--121.

\bibitem [O3] {O3}
\bysame, {\it Dunkl operators, Bessel functions and the 
discriminant of a finite Coxeter group},
Composito Mathematica { 85} (1993), 333--373.

\bibitem [O4] {O4}
\bysame, 
{\it A generating function for the trace of the 
Iwahori -- Hecke algebra}, math. RT/0101006 (2001).

\bibitem [O5] {O5}
\bysame, 
{\it On the spectral decomposition of affine 
Hecke algeb\-ras}, math. RT/0101007 (2001).

\bibitem [St] {St}
{J.V.~Stokman}, 
{\it Difference Fourier transforms for 
non\-reduced root sys\-tems}.
math. CA/\-0105093 (2002).

\bibitem [Va] {Va}
{E.~Vasserot},
{\it Induced and simple modules of 
double affine Hecke algebras}.
Preprint (2002).

\bibitem [Ve] {Ve}
{E.~Verlinde}, {\it Fusion rules and modular 
transformations in $2D$ conformal field theory}, 
Nucl.\ Phys.\ B{ 300}(FS22) (1987),
360--376.

\end{thebibliography}
\end{document}